\documentclass{article}
\usepackage[letterpaper]{geometry}
\usepackage[final]{microtype}

\usepackage{amsfonts,amssymb,amsmath,pdfsync,xcolor,bm,bbm}
\usepackage{amssymb}
\usepackage{amsthm}
\usepackage[shortlabels]{enumitem}
\usepackage[normalem]{ulem}
\usepackage{cancel}

\makeatletter% From gimac.sty
\usepackage{xspace}
\usepackage{graphicx}% Needed so texpresso doesn't bork the initial page
\usepackage{mathtools}
\mathtoolsset{showonlyrefs=true,showmanualtags=true}% Only show referenced equations

\newcommand{\GI@given}[1]{\nonscript\:\mathopen{}#1\vert\nonscript\:\mathopen{}}
\newcommand{\given}[1][]{\GI@given{#1}}
\DeclarePairedDelimiterX\paren[1]{(}{)}{%
\renewcommand{\given}{\GI@given{\delimsize}}#1%
}
\DeclarePairedDelimiterX\brak[1]{[}{]}{%
\renewcommand{\given}{\GI@given{\delimsize}}#1%
}

\newcommand{\GI@st}[1]{\nonscript\:#1\vert\nonscript\:\mathopen{}\allowbreak}
\newcommand{\st}[1][]{\GI@st{#1}}
\DeclarePairedDelimiterX\set[1]\{\}{%
\renewcommand{\st}{\GI@st{\delimsize}}#1%
}

\DeclarePairedDelimiter{\abs}{\lvert}{\rvert}
\DeclarePairedDelimiter{\norm}{\lVert}{\rVert}
\newcommand{\defeq}{\stackrel{\scriptscriptstyle\textup{def}}{=}}

\newcommand{\blfootnote}[1]{%
  \begingroup%
  \renewcommand{\thefootnote}{}%
  \footnotetext{\hspace*{-1.7ex}#1}%
  \endgroup%
}

\makeatother

\usepackage{mathrsfs} %,showlabels}

\newcommand{\Rm}{{\mathbb R}}
\newcommand{\R}{\mathbb R}
\newcommand{\Z}{\mathbb Z}
\newcommand{\eps}{\varepsilon}
\renewcommand{\epsilon}{\varepsilon}
\newcommand{\commentout}[1]{}
\newcommand{\E}{\mathbb{E}}

\newcommand{\N}{\mathbb N}
\newcommand{\Pm}{{\mathbb P}}
\renewcommand{\P}{\mathbb{P}}
\newcommand{\Em}{{\mathbb E}}

\newcommand{\Zm}{{\mathbb Z}}

\newcommand{\cE}{{\cal E}}
\newcommand{\cQ}{\mathcal{Q}}
\newcommand{\cH}{\mathcal{H}}

\usepackage[pdfpagemode=UseNone,breaklinks=true]{hyperref}
\definecolor{link1}{rgb}{0,0,.7}
\definecolor{link2}{rgb}{0,0.25,0.5}
\hypersetup{final,colorlinks,allcolors=link1,citecolor=link2,hypertexnames=false}
\usepackage[numberlinked=false]{footnotebackref}

\theoremstyle{plain}
\newtheorem{theorem}{Theorem}[section]
\newtheorem{lemma}[theorem]{Lemma}
\newtheorem{proposition}[theorem]{Proposition}
\newtheorem{corollary}[theorem]{Corollary}

\newtheorem*{theorem*}{Theorem}
\newtheorem*{lemma*}{Lemma}
\newtheorem*{proposition*}{Proposition}
\newtheorem*{corollary*}{Corollary}

\theoremstyle{definition}

\theoremstyle{remark}
\newtheorem{remark}[theorem]{Remark}
\newtheorem*{remark*}{Remark}

\newcommand{\farc}{\frac}

\newcommand{\one}{\mathbbm{1}}

\newcommand{\cT}{{\mathcal T}}

\newcommand{\cP}{{\mathcal P}}

\newcommand{\bhat}[1]{\expandafter\hat#1}

\DeclareMathOperator*{\argmin}{argmin}

\RequirePackage{seqsplit}
\usepackage{doi}
\definecolor{lightblue}{HTML}{083085}

\numberwithin{equation}{section}

\newcommand{\cG}{{\cal G}}

\newcommand{\dist}{\hbox{dist}}
\newcommand{\grad}{\nabla}
\newcommand{\dv}{\nabla \cdot}
\newcommand{\lap}{\Delta}
\newcommand{\cA}{{\mathcal A}}

\newcommand{\cM}{{\mathcal M}}
\newcommand{\cB}{{\mathcal B}} 
 
\newcommand{\cS}{{\mathcal S}} 
\newcommand{\cU}{{\mathcal U}}
\newcommand{\cZ}{{\mathcal Z}}

\DeclareMathOperator{\conv}{conv}
\DeclareMathOperator{\cut}{cut}

\DeclareMathOperator{\hess}{Hess}
\DeclareMathOperator{\Dom}{Dom}

\newcommand{\Lip}{\mathrm{Lip}}
\newcommand{\BV}{\mathrm{BV}}

\definecolor{GIcolor}{HTML}{349900}
\newcommand{\GIf}[2][]{%
  \ifhmode\unskip\fi%
  \footnote{\textcolor{GIcolor}{%
    \textbf{GI\xspace #1:} #2%
  }}\ignorespaces%
}

\renewcommand{\le}{\leqslant}
\renewcommand{\leq}{\leqslant}
\renewcommand{\ge}{\geqslant}
\renewcommand{\geq}{\geqslant}

\newcommand{\vphi}{\varphi}

\newcommand{\II}{I\!I}

\usepackage{float}

\begin{document}

\author{%
Ayya Alieva%
  \footnote{%
    Institute for Computational Mathematics and Engineering,
    Stanford University,
    Stanford, CA 94305, USA;
    ayya@stanford.edu}%
\and%
Gautam Iyer%
  \footnote{%
    Department of Mathematical Sciences,
    Carnegie Mellon University,
    Pittsburgh,
    PA 15213;
    gautam@math.cmu.edu}%
\and%
Alexei Novikov%
  \footnote{%
    Department of Mathematics,
    The Pennsylvania State University,
    University Park, 16802;
  novikov@psu.edu}%
\and%
Lenya Ryzhik\footnote{%
  Department of Mathematics,
  Stanford University,
  Stanford, CA 94305, USA;
  ryzhik@stanford.edu}}

\title{Support recovery of diffusion models under the manifold hypothesis%
  \blfootnote{%
    \emph{2020 Mathematics Subject Classification.}
       Primary:
	 60H10. %
       Secondary:
	 60J60, %
	 65C30, %
	 68T07. %
    \par
    \emph{Key words and phrases.}
      Diffusion models;
      manifold hypothesis;
      score based generative models;
      cut-locus.
    \par
    GI has been partially supported by the National Science Foundation under grants
    DMS-2406853, %
    DMS-2342349, %
    LR by grants DMS-2205497 and DMS-2510166, and AN by  DMS-2407046.
This work was also supported by the Center for Nonlinear Analysis at Carnegie Mellon University.
  }
}

\maketitle

\begin{abstract}
We study diffusion models under the manifold hypothesis with possibly large score estimation and discretization errors, 
and show that while they may not exactly recover the data distribution, they will exactly recover the support 
(or data manifold). In continuous time, we use a Lyapunov function 
to show that even if the relative score estimation error is fixed (and not necessarily vanishing), 
then the rate at which the process approaches the support is $O(\sqrt{(T-t)|\log(T-t)|})$ as~$t \to T$, where~$T$ is the final time.
In discrete time, we use geometric time steps with a fixed (non-vanishing) ratio and show that with a fixed 
(also not small) relative score 
estimation error, the process approaches the data manifold exponentially with the number of iterations.
The proof relies on controlling a Gaussian measure of a neighborhood of the cut-locus 
of the data manifold, using a polynomial approximation argument combined 
with quantitative semi-algebraic estimates.
\end{abstract}

\tableofcontents

\section{Introduction}
\subsection{Main Results}

The basic paradigm in diffusion models is to use data samples to learn the score of a 
noisy data distribution, and then run a ``denoising'' diffusion with the drift that equals
to the learned score function, in order  to generate new samples from the data 
distribution.
The main result of this paper
is that diffusion models can exactly recover the data manifold even when 
the learned score has a large, 
non-vanishing relative error with respect to the ``true'' score function,
and when the time steps of the discretization are not necessarily small.

Suppose~$\rho_0$ is a probability distribution in~$\R^d$ representing the data.
A commonly used noising schedule uses the regularized family of 
distributions~$\rho_t$, $ 0 < t \leq T$. Here,~$T > 0$ is a fixed terminal time, and
\begin{equation}\label{e:rhoTdef}
\rho_t(x) \defeq \rho_0 * \phi_{2t}(x)\,, \quad
\phi_t(x)
  \defeq \frac{1}{(2\pi t)^{d/2}} e^{-|x|^2/(2t)}
\,,~x\in\R^d.
\end{equation}
The true score (in the sense that will become clear soon) of the noisy data distribution is the vector field%
\begin{equation}\label{e:sStarDef}
s_t(x) \defeq 
  \grad \ln \rho_t(x)
  = \frac{\grad\rho_t(x)}{\rho_t(x)}
  \,.
\end{equation}
The diffusion model simulates the ``denoising'' SDE
\begin{equation}\label{e:dXExactScore}
dX_t = u_t(X_t) \, dt + \sigma \, dB_t,
\end{equation}
where
\begin{equation}\label{e:uDef}
u_t(x)\defeq \gamma s_{T - t} (x)= \frac{\gamma \grad \rho_{T - t}(x)}{\rho_{T - t}(x)},
\end{equation}
and
\begin{equation}\label{e:SigmaDef}
\sigma = \sqrt{2 (\gamma - 1)},
\end{equation}
on the time interval $0<t<T$.
Here~$B_t$ is a standard~$d$-dimensional Brownian motion,  
and~$\gamma >1$ is a parameter.

A direct substitution in the Fokker-Planck equation corresponding to~\eqref{e:dXExactScore}--\eqref{e:uDef}
(observed first by Anderson~\cite{Anderson82}) shows that if~$X_0 \sim \rho_T$ then~$X_T \sim \rho_0$.
In this sense, the SDE~\eqref{e:dXExactScore} is a denoising process that undoes the noising~\eqref{e:rhoTdef} and exactly recovers~$\rho_0$.
This is also why~$s_t(x)$ given by~\eqref{e:sStarDef} is referred to as the true score of the distribution $\rho_0$. 
The most common version of~\eqref{e:dXExactScore} found in the literature is with~$\gamma = 2$ in~\eqref{e:uDef} and~\eqref{e:SigmaDef}.
As discussed below, this choice has certain advantages.   
The restriction $\gamma>1$ keeps $\sigma>0$. The presence of the noise, $\sigma\neq 0$ in~\eqref{e:SigmaDef}, does help us 
to prove convergence in Theorems~\ref{t:relErrorScore} and~\ref{t:supportRecoveryDiscreteIntro} below. 
The noiseless case~$\gamma=1$, known as  the probability flow ODE, is also sometimes used. 
One of its main advantages is the ability to use higher order or even implicit time stepping, 
see~\cite{song2022denoisingdiffusionimplicitmodels,lu2022dpmsolverfastodesolver}. However,
these ODE methods tend to suffer from convergence and mode collapse 
(only recovering some specific parts of the support of the distribution $\rho_0$)
issues and often require SDE simulations in practice.

Of course, it is natural to expect that solutions to a discrete approximation to~\eqref{e:dXExactScore} converge to the true solution, provided the approximation and discretization errors vanish.
To elaborate, there are three principal sources of error: (1) score approximation error, 
(2)~time-discretization error, and (3) initialization error.
The score approximation error arises because the true score~$s_t$ in~\eqref{e:sStarDef} 
is not exactly known but is approximated through a computationally expensive training process.
The time discretization error appears because one can only simulate time discrete approximations of~\eqref{e:dXExactScore}, 
and, as will become apparent below, the drift term is routinely large or even unbounded.
Moreover, in practice, often only a very limited number of time steps can be taken per denoising process, since the score function is approximated by a frequently large neural network that is expensive to evaluate at each time step.
Finally the initialization error arises when one does not start~\eqref{e:dXExactScore} with~$X_0 \sim \rho_T$ (as needed to exactly recover~$X_T \sim \rho_0$), but instead initializes~$X_0$ according to a known distribution (typically, a Gaussian with a large variance). 
Because the approximation to~$s_t$ is trained jointly for all~$t\in(0,T]$, there is a trade-off between model accuracy 
at~$t \approx 0$ and the size of $T$~\cite{zhang2025collapseerrorsinduceddeterministic}. As a result, the initialization error cannot be simply controlled by selecting $T \gg 1$, in which case $\rho_T$ would be close to a Gaussian.
Due to the above considerations, in many situations, each of these errors can not be controlled, and may very well be large.
 
Motivated by this, we study~\eqref{e:dXExactScore}
and obtain ``non-asymptotic results'', in the sense that we do not require the approximation 
errors to be small. When the above errors are not negligible, one can not hope to recover the original data 
distribution~$\rho_0$.
However, in many situations of interest, $\rho_0$ is supported on a low dimensional ``data manifold''
that we denote by $\cM\subset\Rm^d$ \cite{fefferman2013testingmanifoldhypothesis,Gorban_2018,7080175,NIPS2004_74934548}.
In this case, we show that the diffusion models still  exactly
recover the data manifold, even when the score approximation error, time-discretization error and initialization error are all large.
This is important for practical applications, as we guarantee that the model will only produce points on the (low dimensional) data manifold instead of 
terminating at arbitrary points in the (high dimensional) ambient space, even if it may not recover the true data distribution on the manifold. This
``recovery only of support" phenomenon has been observed in practice, where it is usually described as``error in mode amplitude" 
or``model bias towards specific subgroups", see e.g. 
\cite{dennehy2026diffusionmodelsrecoveraccurate,koehler2022statisticalefficiencyscorematching,roos2025imetbiasinvestigating,hakemi2025deeperdiffusionmodelsamplify}. 
Our findings underline the true power of diffusion models: even when the problem is too complex to solve exactly, 
the denoising SDE will still converge onto an acceptable, realistic, sample.

Our main results can be informally summarized as follows:
\begin{enumerate}
  \item 
In continuous time, if~$\gamma \geq 2$, the terminal time $T>0$ is fixed, 
and one learns the score up to a relative error of at most~$1 - 1/\gamma$, then the diffusion model almost surely approaches the manifold~$\mathcal M$ as~$t \to T$ at rate~$O(\sqrt{T - t})$, up to a logarithmic factor.
  While a few similar results are available~\cite{
      DeBortoli22,
      Pidstrigach22,
      LuWangEA24,
      LiShenEA26},
    we provide a new and elegant proof by constructing a Lyapunov function.
    In addition, our approach allows us to prove a stronger (quantitative) result under significantly
    weaker assumptions on the score estimation error.
    We only require the aforementioned 
    bound on the relative score approximation error, which is   important, as it 
    allows for very large errors as~$t$ approaches~$T$ and does not even require the relative error
    to be small.

\item
For the discrete time algorithm, we fix $T>0$ and use geometric time steps with a fixed ratio~$\alpha \in (0, 1)$.
If we learn the score up to a relative error of~$\delta < 1$, then the distance to the data manifold vanishes exponentially with the number of iterations.
To the best of our knowledge, there are no previously known results to this effect that apply to the discrete time diffusion model.
These can not be obtained from the continuous time results since the mesh-size (in time) is fixed, not vanishing, the coefficients in~\eqref{e:dXHat} are not bounded, and the score approximation error need not vanish.
\end{enumerate}
 
We now state our main results precisely.

\subsubsection{Support recovery in continuous time}

We first state our results in the continuous time case.
Let~$\hat s_t(x)$ be an approximate score function obtained in training.
We define the denoising process~$\hat X_t$ 
based on this approximate score by
\begin{equation}\label{e:dXHat}
d\hat X_t = \hat u_t(\hat X_t) \, dt + \sigma \, dB_t,
\end{equation}
where
\begin{equation}\label{e:uHatDef}
\hat u_t(x) \defeq \gamma \hat s_{T - t}(x),
\end{equation}
and the diffusivity $\sigma>0$ is still given by~\eqref{e:SigmaDef}.
We reiterate that the final time $T>0$ is fixed and is not assumed to be large.

Our first result shows that if~$\gamma \geq 2$ and we make a relative error  
of size at most~$1 - 1/\gamma$, then we still exactly recover the support of the initial distribution.
In the commonly used case~$\gamma = 2$, this allows an at most~$50\%$ relative error 
in learning the score.
Moreover, we obtain quantitative bounds on the rate at which~$d(\hat X_t, \mathcal M) \to 0$.
Here, and below, we denote by $d(x,\cM)$ the distance between a point $x\in\Rm^d$ and a sub-manifold $\cM\subset\Rm^d$.
\begin{theorem}\label{t:relErrorScore}
  Let~$m < d$, and~$\mathcal M$ be an~$m$-dimensional closed compact~$C^2$ sub-manifold of~$\R^d$.
  Let~$\rho_0$ be a probability distribution which has a bounded, non-degenerate density~$\bar \rho_0$ with respect to the volume measure on~$\mathcal M$.
  That is, we assume
  \begin{equation}\label{e:rho0Bounded}
    \rho_{0, \inf } \defeq \inf_{z \in \mathcal M} \bar \rho_0(z) > 0
    \quad\text{and}\quad
    \rho_{0, \sup } \defeq \sup_{z \in \mathcal M} \bar \rho_0(z) < \infty
    \,.
  \end{equation}
Suppose that $\gamma\ge 2$, and there exists
  \begin{equation}\label{e:deltaGamma}
    0 \leq \delta < 1 - \frac{1}{\gamma}
    \,,
  \end{equation}
  such that for every~$x \in \R^d$, $t \in (0, T]$ we have
  \begin{equation}\label{e:RelErrorScore}
    \abs{\hat s_t(x) - s_t(x)} \leq \delta \abs{s_t(x)}
    \,.
  \end{equation}
  Let~$\hat u_t(x)$, $\hat X_t$ be as in~\eqref{e:dXHat}--\eqref{e:uHatDef}.
  Then, for any initial distribution of~$\hat X_0$, we have, almost surely
  \begin{equation}
    \label{e:dXtHatMrate}
    \limsup_{t\to T^-}\paren[\Big]{\frac{1}{(T-t)|\log(T-t)|}}^{1/2}
      d(\hat X_t, \mathcal M)
      \leq \sqrt{2 \paren[\big]{(\gamma - 1)d - m}}
    \,,
 \end{equation}
 and
 \begin{equation}
    \label{e:uHatConvRate}
    \limsup_{t\to T^-}
    \paren[\Big]{\frac{T-t}{|\log(T-t)|}}^{1/2} \abs{\hat u_{t}(\hat X_t)}
      \leq \frac{\gamma(1+\delta)}{\sqrt{2}} \sqrt{(\gamma - 1) d - m}
    \,.
  \end{equation}
\end{theorem}

We prove Theorem~\ref{t:relErrorScore} by finding a Lyapunov function and 
using the martingale convergence theorem.
Our proof can be modified to also give quantitative bounds on the rate at 
which~$\E [d(\hat X_t, \mathcal M)] \to 0$, as outlined in Section~\ref{sec:lim-super} below.  

\begin{remark}\label{r:convRateSharp}
The convergence rates~\eqref{e:dXtHatMrate} and~\eqref{e:uHatConvRate} in Theorem~\ref{t:relErrorScore} are at worst sub-optimal by a logarithmic factor, as can be seen by choosing~$\rho_0$ to be a Dirac~$\delta$-measure.
This is discussed further in Section~\ref{s:deltaExample}, below. 
Further, note that the size of the score $s_{T-t}(x)$ blows up as~$(T-t)^{-1}$ for $x$ away from the
manifold $\cM$. Thus, the upper bound on $\hat u_t(\hat X_t)$ in~\eqref{e:uHatConvRate} is another 
way to indicate that $\hat X_t$ is close to $\cM$. It also explains that one need not train the score
$\hat s_t(x)$ in regions where it takes ``huge'' values.
\end{remark}

\begin{remark}\label{r:GammaEq2Conv}
Under the assumptions of Theorem~\ref{t:relErrorScore},
the process~$\hat X_t$ itself converges almost surely
to a point on~$\mathcal M$ as~$t\to T^-$.
This is simply because the velocity bound~\eqref{e:uHatConvRate}
is integrable near~$T$.
\end{remark}

\subsubsection{Support recovery in discrete time}

We now describe our results for the discrete time approximations to~\eqref{e:dXExactScore}.
A common approach in practice~\cite{karras2022elucidatingdesignspacediffusionbased} is to discretize~\eqref{e:dXExactScore} using \emph{geometric time steps} -- that is, the time steps are chosen so that the increments are geometrically decreasing.
Explicitly, fix~$\alpha \in (0, 1)$, and define
\begin{equation}\label{e:tnDef}
  t_n = T (1 - \alpha^n)\,,
  \quad
  \Delta t_n = t_{n+1} - t_n = T(1- \alpha) \alpha^n
  \,.
\end{equation}
An Euler--Maruyama based discretization of~\eqref{e:dXHat} is to set
\begin{equation}\label{e:EMhatX} 
  \hat X_{n+1} = \hat X_n + \hat u_{t_n}(\hat X_n) \Delta t_n + \sigma \sqrt{ \Delta t_n} \, \xi_{n+1}\,,
  \quad
  \sigma = \sqrt{2 (\gamma - 1)}
  \,,
\end{equation}
where~$\xi_n$ are i.i.d.\ standard Gaussians.
The largest time increment (or mesh size in time) is precisely~$T ( 1- \alpha)$, which converges to~$0$ as~$\alpha \to 1$.
However, since the coefficients of~\eqref{e:dXHat} are unbounded, solutions to the Euler--Maruyama discretization~\eqref{e:EMhatX} still may not converge to the solution of~\eqref{e:dXHat} as~$\alpha \to 1$.
 
We also note that there is a trade-off between making~$\alpha \ll 1$, and~$\alpha \approx 1$.
Making~$\alpha \approx 1$ results in smaller time steps, and an approximation that is closer to the continuous model.
On the other hand, choosing~$\alpha \ll 1$ makes the chosen times~$t_n$ converge to~$T$ faster, which is numerically advantageous, 
and often required due to time and cost constraints.
Similarly, while larger $T$ results in smaller initialization error between the commonly used Gaussian distribution for $X_0$ and $\rho_T(x)$,
increasing $T$ can often compromise stability and accuracy of the final result, rendering smaller $T$ a better choice.
 
We fix~$\alpha \in (0, 1)$ (which isn't necessarily close to~$1$), and~$T > 0$ (which isn't necessarily large), and ask whether   
the iterations of~\eqref{e:EMhatX} still converge to the support of~$\rho_0$ as~$n \to \infty$.
Our main result shows that if we use an approximate score~$\hat u_t$ with a relative error coming from~\eqref{e:uHatDef} 
and~\eqref{e:RelErrorScore} for some~$\delta < 1$, then the iterates~$\hat X_n$ approach~$\mathcal M$ exponentially fast as~$n \to \infty$.
For simplicity, we state the result here with~$\gamma = 2$.
The general result for any~$\gamma > 1$ is Theorem~\ref{t:supportRecoveryDiscrete} in Section~\ref{s:discreteConv}, below.

\begin{theorem}\label{t:supportRecoveryDiscreteIntro}
Let~$m < d$, and~$\mathcal M$ be an~$m$-dimensional compact closed~$C^\infty$ sub-manifold of~$\R^d$.
Let~$\rho$ be as in Theorem~\ref{t:relErrorScore}.
Suppose~$\delta \in (0, 1)$ and~$\hat s_t$ is an approximate score function satisfying 
the relative error bound~\eqref{e:RelErrorScore}, and let~$\hat u$ be the associated velocity field defined by~\eqref{e:uHatDef}.
Let $\hat X_{0} \in\R^d$ be an arbitrary initial point and inductively define~$\hat X_n$ by~\eqref{e:tnDef}--\eqref{e:EMhatX} 
with~$\gamma = 2$ and a fixed $\alpha \in (0, 1)$.
Given $b\in(0,1)$, choose any~$\bar \alpha_b$ such that
\begin{equation}\label{e:barAlphaB}
  \max\set{\alpha^b, \alpha + \delta(1 - \alpha)} < \bar \alpha_b < 1
  \,.
\end{equation}
\begin{enumerate}
  \item
For all~$b < 1/6$  
there exists a constant~$C = C(\alpha, b, \gamma, \mathcal M)<+\infty$ such that
    \begin{equation}\label{e:EdXnHatMGamma2}
      \E d(\hat X_n, \mathcal M) \leq
	    (d(\hat X_0, \mathcal M) + C) \bar \alpha_b^n,
      \,.
    \end{equation}

  \item
For all~$b < 1/4$  
there exists a random constant~$D = D(\alpha, \gamma, b, X_0, \mathcal M)<+\infty$ such that,
almost surely, for every~$n \in \N$, we have
    \begin{equation}\label{e:dXnHatMGamma2}
      d(\hat X_n, \mathcal M)
	\leq
	    D \bar \alpha_b^n
	    \,.
    \end{equation}
\end{enumerate}
\end{theorem}

\begin{remark}
When we have access to the exact score (that is, $\delta = 0$), 
then Theorem~\ref{t:supportRecoveryDiscreteIntro}
guarantees that~$d(\hat X_n, \mathcal M)$
converges to~$0$ almost surely with rate almost~$\alpha^{n/4}$, and in~$L^1$ with rate 
almost~$\alpha^{n/6}$.
In fact, if~$\hat X_n$ is confined to a tubular neighborhood around~$\mathcal M$, the proof will 
show that~$d(\hat X_n, \mathcal M)$ vanishes at the faster rate~$\alpha^{n / 2}$.
In this case one can drop the noise entirely, and get convergence at the even faster rate 
of~$\alpha^n$.
\end{remark}

The proof of Theorem~\ref{t:supportRecoveryDiscrete} is significantly more involved than that of Theorem~\ref{t:relErrorScore}. In the continuous case, existence of the Lyapunov function 
tells us that the process $\hat X_t$ will go towards the manifold, albeit with no convergence rate. 
This allows us 
to obtain the rate of convergence to the manifold, as in~\eqref{e:dXtHatMrate}--\eqref{e:uHatConvRate}, 
purely from an asymptotic expression for the score \emph{close} to the manifold.
Explicitly,  the score~$s_t(x)$ satisfies
\begin{equation}\label{e:scoreTube}
  s_t(x) \approx -\frac{x - x^*}{2 t}
  \quad \text{as } t \to 0
  \,.
\end{equation}
Here,~$x^*$ is the unique point in~$\mathcal M$ that is closest to~$x$.
Note that~$x^*$ is well defined in a tubular neighborhood around~$\mathcal M$,
and the Lyapunov function tells us that eventually $\hat X_t$ will enter that neighborhood.
A similar expansion was also used in the related continuous time results~\cite{DeBortoli22,Pidstrigach22,LuWangEA24}, which are discussed 
further in Section~\ref{s:closelyRelated}, below. The asymptotic behavior in~\eqref{e:scoreTube} 
allows us to deduce the rates of convergence in the continuous case after $\hat X_t$ entered this
good tubular neighborhood.

As one moves away from the manifold, the asymptotics~\eqref{e:scoreTube} remains true (with some caveats) until you get to the \emph{cut-locus}.
Recall that the cut-locus of a manifold~$\mathcal M \subseteq \R^d$ is the set of all points in~$\R^d$ which \emph{do not} have a unique closest point in~$\mathcal M$.
That is, the cut-locus\footnote{
This is also called the \emph{medial axis} or \emph{skeleton} of~$\mathcal M$.
We note that some authors define the cut-locus as the closure of the set in~\eqref{e:cutLocus}.
This may additionally include focal points at which the distance to the manifold has a degenerate minimum.
}
of~$\mathcal M$ is defined by
\begin{equation}\label{e:cutLocus}
  \cut(\mathcal M)
    \defeq
      \set[\Big]{x \in \R^d \st \# \argmin_{z \in \mathcal M} |z - x| \geq 2 }
      \,.
\end{equation}

At the cut-locus, the point~$x$ does not have a unique closest point on~$\mathcal M$, and the right side of~\eqref{e:scoreTube} is not even defined.
It is a neighborhood of these points that is a potential bottleneck to convergence of the discrete
algorithm.
This bottleneck is also likely encountered in practice, as
at the initial stage, when~$t$ is far from $T$, the score pushes all points towards the barycentre, which is often close to the cut-locus. 
In the discrete case, the proof of Theorem~\ref{t:supportRecoveryDiscrete} requires understanding the asymptotic behavior of~$s_t$ everywhere, including the ``bottleneck'' points close to the cut-locus.
The proof relies on two main ingredients: first, a quantitative version of how~\eqref{e:scoreTube} 
degenerates near the cut-locus, see Lemma~\ref{l:zetaBdKTau-bis2} below. 
Second, we estimate the probability that the iterates~$\hat X_n$ avoid a small tubular neighborhood 
of the cut-locus.

This last step, the proof of the fact that $X_n$ avoid a small neighborhood of the cut-locus is surprisingly difficult, and the proof is substantial enough that it merits a mention here.
The specific cut-locus avoidance lemma we need is as follows.

\begin{lemma}\label{l:ballVolumeFractions2}
For all~$\epsilon>0$ and $R > 0$ there exists a constant~$C = C(R, \epsilon)$ such that for every~$r > 0$, $\kappa \in (0, 1]$, and $x \in \R^d$, we have
\begin{equation}\label{e:PNinCutRM2}
  \P( d(\mathcal N(x, \kappa I_d), \cut_{R}(\mathcal M)) \leq r )
    \leq \frac{C r}{\kappa^{\frac{1}{2} + \epsilon}}
  \,.
\end{equation}
Here~$\cut_R(\mathcal M) \defeq \cut(\cM) \cap B(0, R)$, $\mathcal I_d$ is the~$d\times d$ identity matrix, and~$\mathcal N(x, \kappa I_d)$ is a Gaussian random variable with mean~$x \in \R^d$ and  covariance matrix~$\kappa I_d$.
\end{lemma}

If the~$\cut_R(\mathcal M)$ is contained in the \emph{graph} of a Lipschitz function of a compact~$(d-1)$ dimensional manifold, then the Gaussian bound~\eqref{e:PNinCutRM2} can be proved quickly even 
with~$\epsilon = 0$.
However, classical results of Itoh and Tanaka~\cite{ItohTanaka01} and Li, Nirenberg~\cite{LiNirenberg05} 
show that the cut-locus (restricted to a large ball) is contained in the \emph{image} (not graph!) of a Lipschitz function of a compact~$(d-1)$ dimensional manifold.
This small and subtle difference makes the proof much harder.

In general, $\cut_R(\mathcal M)$ can be a wild fractal set and even dense~\cite{GluckSinger78,GluckSinger79}, and in such cases one may not be able to prove~\eqref{e:PNinCutRM2} with~$\epsilon = 0$.
This, however, likely only happens in exceptional cases~\cite{Thom72}.
All real-analytic manifolds~\cite{Buchner77,Wall77} have a regular cut-locus of lower dimension and will satisfy~\eqref{e:PNinCutRM2} with~$\epsilon = 0$.
Yomdin~\cite{Yomdin81} (see also~\cite{Mather83}) showed that the cut-locus of a ``generic'' 
smooth sub-manifold of~$\R^d$ is stratifiable, with strata of dimension at most~$d - 1$.
Without quantitative control of the embeddings of the strata, however, one can not guarantee~\eqref{e:PNinCutRM2} holds, even for some~$\epsilon > 0$.
We prove Lemma~\ref{l:ballVolumeFractions2} without appealing to the stratifiable structure, 
using a polynomial approximation and 
a quantitative semi-algebraic approximation~\cite{BasuPollackEA06}.

\subsubsection{A numerical illustration of the effect of the cut-locus} 

\begin{figure}[htb]
    \includegraphics[width=\linewidth]{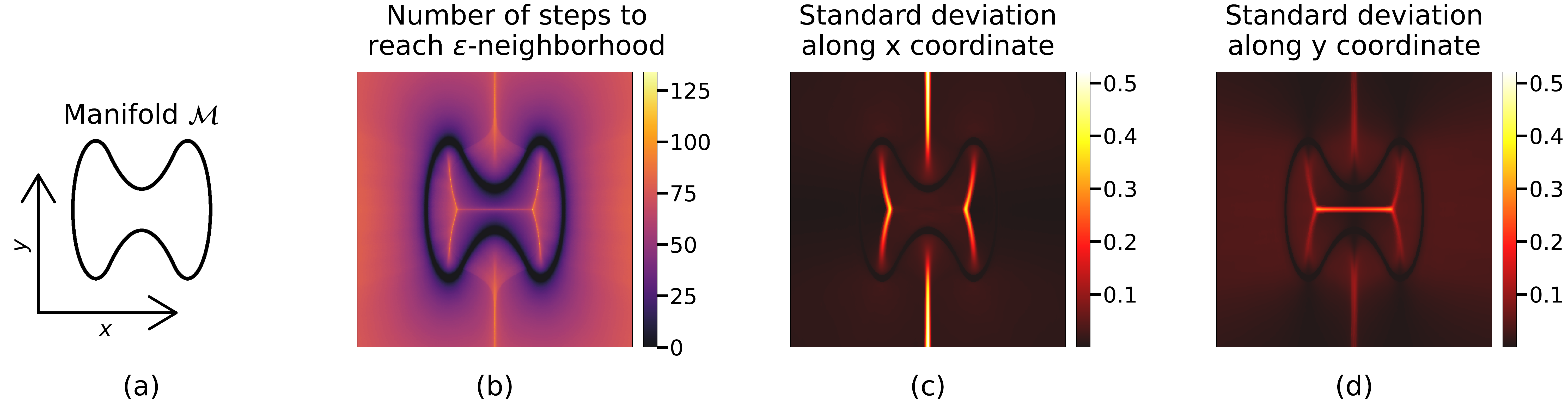}
    \caption{Illustration of the influence of the cut locus on the diffusion model. }
\label{fig:cut_locus}
\end{figure}

Figure~\ref{fig:cut_locus} illustrates the effect  of the cut-locus on a simple two-dimensional example.
The manifold itself is depicted in Figure~\ref{fig:cut_locus}(a). 
In Figures (b)-(d), we approximated the score by a neural network, and generated samples using~\eqref{e:EMhatX}.
Figure~\ref{fig:cut_locus}(b) shows the number of steps needed to reach an~$\epsilon$-neighborhood of~$\mathcal M$ as a function of the starting point.
The figure clearly shows that the number of steps required to reach an $\eps$-neighborhood of the manifold is higher when the starting point is close to the cut-locus. As minimizing the number of time steps at evaluation has significant practical interest (e.g. \cite{li2023autodiffusiontrainingfreeoptimizationtime,pmlr-v235-sabour24a}), the cut-locus is likely to have practical implications.
In the simulations, we choose~$\varepsilon=0.05$, $\alpha=0.95$, and~$T = 10^{-6}$.
These parameters were chosen to emphasize and illustrate the effects of the cut-locus.
In practice one typically uses much larger values of~$T$ (to reduce initialization errors), and smaller values of~$\alpha$ (so that the number of iterations is smaller).
The cut-locus effects are observable both for larger~$T$ (but are blurred due larger initial step), and for smaller~$\alpha$ (but the figures have sharper heatmap transitions).

Figures~\ref{fig:cut_locus}(c) and~\ref{fig:cut_locus}(d) illustrate an effect related to applied work on a phenomenon known as \emph{model variance at intermediate noise scales}~\cite{xu2023stabletargetfieldreduced}.
This terminology refers to an instability in the models when the denoising process~\eqref{e:dXHat} starting at initial conditions which are close produces very different final outcomes. 
These regions of high variance likely correspond to the cut locus, as this is precisely where the direction of the score changes rapidly. 
Figures~\ref{fig:cut_locus}(c) and~\ref{fig:cut_locus}(d) show the standard deviation of the $x$ and $y$-coordinates of $X_T$ as a function of the starting point $X_0$, 
respectively.
Together, they demonstrate that the variance of the outcome of \eqref{e:EMhatX} is large on all of the cut-locus.
For the variance calculations in Figures~\ref{fig:cut_locus}(c) and~\ref{fig:cut_locus}(d), we ran 200 simulations for each $X_0$ on the grid, and used $T=0.001$.
Smaller values of $T$ produces equivalent images, and larger values of~$T$ blurs the cut-locus

\subsection{Literature overview}\label{s:litReview}

\subsubsection{General diffusion models}

Sampling from complex, high-dimensional data distributions is a longstanding problem, 
with many applications ranging from machine learning to molecular dynamics (see for 
instance~\cite{LelievreRoussetEA10,Liu08}).
Since foundational work of Neal~\cite{Neal01}, a commonly used technique is to initially smooth 
out local structures, and gradually reduce the smoothing until the true distribution is recovered.
The advantage of this approach is the reduction in complexity -- sampling the true 
distribution directly is typically exponentially expensive 
in the barrier height between local structures, while  
an ``annealed'' algorithm can reduce this to polynomial
(see, for instance, the book~\cite{Liu08}, or a few recent results~\cite{HanIyerEA26,LeeSantanaGijzen24,PaulinJasraEA19}, and references there-in).

Annealing based strategies typically smooth out local structures by increasing a temperature parameter, 
which effectively reduces the barrier between local structures, and flattens the distribution.
\emph{Diffusion models}, on the other hand, smooth out the local structures by a Gaussian 
regularization, and then run a denoising process to recover the true distribution.
The only unknown in the denoising process is the score function $s_t = \nabla \log \rho_t$ of the 
noised data distribution, which is approximated in practice using neural networks via denoising score 
matching algorithms~\cite{Hyvarinen05, Vincent11}.
This had tremendous success recently in score based generative 
models~\cite{SongErmon19,SongSohlDicksteinEA21}, and denoising diffusion probabilistic 
models~\cite{HoJainEA20}. We refer the reader to the review articles~\cite{Lu26,YangZhangEA25} for a 
comprehensive introduction to this rapidly developing field.

The empirical success of diffusion models led to several mathematical studies establishing 
rigorous convergence bounds. A central line of inquiry studies the discrepancy between the true 
data distribution and the generated distribution.
When the target data distribution~$\rho_0$ is absolutely continuous with respect to the 
Lebesgue measure in the ambient space $\R^d$, an extensive body of literature now exists.
Lee et\ al.\ and Chen et\ al.~\cite{ChenChewiEA23,LeeLuEA22} provide polynomial-time 
convergence guarantees in total variation and Kullback-Leibler divergence, but require strong 
regularity assumptions
(e.g. smoothness, strict log-concavity, or a log-Sobolev/Poincar\'e constant bound)
on the target distribution.
Block et\ al.~\cite{BlockMrouehEA22}, established non-asymptotic bounds in Wasserstein 
distance \cite{BlockMrouehEA22}, accommodating less regular potentials, and Benton et\ al.~\cite{BentonDeBortoliEA24} analyzed the impact of various discretization schemes.

In many practical applications, however, the target distribution is singular, and not regular.
In fact, many modern machine learning applications (such as image or audio generation) postulate the \emph{manifold hypothesis}: that in spite of the high ambient dimension, the data itself is \emph{intrinsically} low dimensional and supported on a lower dimensional manifold~$\mathcal M \subseteq \R^d$.

The analysis when the target distribution is singular is much more involved.
The denoising process is necessarily degenerate at the final time~$T$, and stability estimates 
for such processes are difficult to obtain.
In addition to the analytical challenges, there is also an important practical consideration: due to the size and intricate structure of the data, diffusion models require large neural networks that are notoriously difficult to train, thus making it hard to guarantee that the score approximation errors are small in practice.
Combined with time constraints on evaluation, the denoising process frequently results in non-vanishing, potentially large discretization errors.
Thus, even if good stability estimates are available, it would be hard to verify applicability.

In this paper we take a different approach.
Rather than performing a stability/error analysis, we study support recovery.
We show that even though the denoising process runs in a high dimensional ambient space, it rapidly approaches the (low dimensional) data-manifold, even if the score estimation errors or time discretization errors are not small.
This is significant in practice, as it means that the model will always terminate on the data 
manifold (for example, produce an acceptable image) even though it may not recover the true data 
distribution on the manifold.
There are a few related results in this direction, and we summarize them and compare them to our work in the next section.

\subsubsection{Comparison to closely related work}\label{s:closelyRelated}

Among the results most closely related to ours 
are \cite{DeBortoli22,LiShenEA26,LuWangEA24,Pidstrigach22}.
De Bortoli~\cite{DeBortoli22} establishes~$1$-Wasserstein distance bounds between the generated distribution and the true data distribution supported on a low-dimensional manifold. 
The obtained bounds require vanishing, or an~$L^2$-small absolute score error.
In addition, to avoid singularities of the score, 
an early stopping time   is introduced in~\cite{DeBortoli22},  
which effectively regularizes the data.
De Bortoli also considers time discrete approximations, but uses small 
polynomial time-steps in order to combat the growth in coefficients.

Our results are complementary to De Bortoli's, as we only consider support recovery, 
and not quantitative error bounds.
However, in our results we allow for a (fixed, non-vanishing) relative error in the score, which 
corresponds to an extremely large absolute error.
In the discrete case, we are also able to handle geometric time steps, do not require them
to be small, and show that the distance to the manifold vanishes exponentially 
with the number of iterations even with these large relative errors.

Pidstrigach~\cite{Pidstrigach22} proves a support recovery result for continuous 
time diffusion models, without introducing an early stopping time, 
whereas Lu et\ al.~\cite{LuWangEA24} study the singular behavior of the score near the data manifold.
While Lu et\ al.\ use the exact score, Pidstrigach allows for an absolute error in the score, 
provided the error satisfies an exponential integrability condition.
Both papers use an asymptotic expansion of the score in a tubular neighborhood of~$\mathcal M$, 
as in~\eqref{e:scoreTube}. More specifically, 
Pidstrigach uses it to prove convergence to the support, while Lu et\ al.\ use it to analyze the behavior of the diffusion near the support.
Compared with~\cite{LuWangEA24,Pidstrigach22}, Theorem~\ref{t:relErrorScore} 
allows for a much larger relative error in the score, which allows the absolute error 
to be~$O(1/(T - t))$, which is far from satisfying the exponential integrability condition 
of Pidstrigach~\cite{Pidstrigach22}.
Theorem~\ref{t:relErrorScore} also gives a quantitative~$O((T - t)\abs{\log(T - t)})^{1/2}$ 
rate of convergence as in~\eqref{e:dXtHatMrate}, \eqref{e:uHatConvRate}.
The proof of Theorem~\ref{t:relErrorScore} is also very different from the proofs 
in~\cite{LuWangEA24,Pidstrigach22}.
Our proof constructs a Lyapunov function, whereas~\cite{Pidstrigach22} 
uses the Girsanov theorem and a relative entropy argument, 
and~\cite{LuWangEA24} employs a direct trajectory analysis near the data manifold 
for singular SDEs via asymptotic 
ODE stability.

Finally, we mention the work of Li et\ al.~\cite{LiShenEA26}.
They consider scores~$\hat s_t$ with an absolute local error of the 
size~$o(1/t)$ as~$t \to 0$.
If~$\hat \rho_t$ is the distribution whose score is~$\hat s_{t}$, then Li et\ al.\ 
show that the support of~$\hat \rho_t$ concentrates on~$\mathcal M$ as~$t \to 0$.
This is different from Theorem~\ref{t:relErrorScore}, as \cite{LiShenEA26}
studies the static snap-shot  stationary distributions, and does not analyze the 
dynamics of the iterates.

\subsection*{Plan of this paper}

In Section~\ref{s:continuousConv}, we prove Theorem~\ref{t:relErrorScore} and 
investigate the optimality of Theorem~\ref{t:relErrorScore} 
by comparing to the explicitly computable example when the data distribution is a Dirac measure.
In Section~\ref{s:discreteConv}, we state and prove Theorem~\ref{t:supportRecoveryDiscrete},
a generalization of Theorem~\ref{t:supportRecoveryDiscreteIntro}   
to include all~$\gamma > 1$. 
The remainder of the paper is devoted to lemmas required for the proofs of Theorems~\ref{t:relErrorScore} and~\ref{t:supportRecoveryDiscrete}.
Section~\ref{s:lyapunov} constructs a Lyapunov function (Lemma~\ref{l:dXtHatMVanish}, below) used in the 
proof of  Theorem~\ref{t:relErrorScore}. Section~\ref{s:rhoDecay} studies the decay of~$\rho_t$ away 
from the manifold~$\mathcal M$ (proving Lemma~\ref{l:rhoDecay}), both of which are needed in the proof of Theorem~\ref{t:relErrorScore}.
Section~\ref{s:zetaConv} obtains the asymptotic behavior of Gaussian averaged projections as~$t \to T$ (proving Lemma~\ref{l:zetaBdKTau-bis2}), which is used in both Theorems~\ref{t:relErrorScore} and~\ref{t:supportRecoveryDiscrete}.
Finally, both Sections~\ref{s:ballVolFractions2} and~\ref{s:cutTubeBd} are devoted to the proof of Lemma~\ref{l:ballVolumeFractions2}.

\subsection*{Tool and computational resource disclosure}

The proof of Lemma~\ref{l:ballVolumeFractions2} but not its statement 
was developed with the assistance of ChatGPT.
The proof was completely written and checked by the authors, and the authors are responsible for its correctness.
Proof assistants were not used significantly in the remainder of this paper, except for checking the misprints and minor technical errors.
The numerical experiment was computed using 1 Tesla P100 PCIe 16 GB GPU.

\section{Support recovery in continuous time}\label{s:continuousConv}

In this section, we first give the proof of Theorem \ref{t:relErrorScore},
though we postpone the proof of Lemmas~\ref{l:dXtHatMVanish},~\ref{l:rhoDecay} 
and~\ref{l:zetaBdKTau-bis2} until later in the paper.  In Section~\ref{s:deltaExample},
we consider a special case~$\rho_0=\delta_0$, in which everything can be computed explicit.

\subsection{Proof of Theorem \ref{t:relErrorScore}}\label{s:relErrorScoreProof}

There are three ingredients that go into the proof of Theorem~\ref{t:relErrorScore}.
The first is finding a (time dependent) Lyapunov function that captures the attracting effect of the score towards the manifold.
We will use the true score of the noisy distribution~$\rho_t$ and construct 
a Lyapunov function as follows. Fix $\beta\in(0,1]$ and 
define the function~$V \colon [0, T) \times \R^d  \to (0, \infty)$ by
\begin{equation}\label{e:Vdef}
  V_t(x) = \frac{(T-t)^p}{ [\rho_{T-t}(x)]^{\beta}},
\end{equation}
where
\begin{equation}\label{26sep902}
p = \frac{\beta(\gamma-2)d}{2}.
\end{equation}
We claim that~$V_t(x)$ is the desired Lyapunov function in the following sense.
\begin{lemma}\label{l:dXtHatMVanish}
Let~$V_t(x)$ be defined by~\eqref{e:Vdef}, choose~$\beta$ so that
\begin{equation}\label{e:deltaConditionGamma}
  0 < \beta \leq 1 - \frac{\delta \gamma}{\gamma - 1}
\end{equation}
and set
\begin{equation}\label{e:zHatDef}
\hat Z_t\defeq V_t(\hat X_t)
  \,.
\end{equation}
The process~$\set{Z_t}_{0 \leq t < T}$ is a nonnegative super-martingale, and converges almost surely 
to an integrable limit as~$t \to T$.
Consequently,~$d(\hat X_t, \mathcal M) \to 0$ almost surely as~$t \to T^-$.
\end{lemma}

\begin{remark}
  When~$\delta = 0$, $\gamma = 2$ the process~$\set{\hat Z_t}_{0 \leq t < T}$ is in fact a martingale, and not only a super-martingale.
\end{remark}

The proof of Lemma~\ref{l:dXtHatMVanish} is presented in Section~\ref{s:lyapunov}, below.
Since~$\hat Z_t = V(\hat X_t)$ is a super-martingale, the process $\hat X_t$ 
will move towards regions in~$\R^d$ where~$V$ is smaller.
The next ingredient in the proof of Theorem~\ref{t:relErrorScore} is to obtain a rate at which the Lyapunov function~$V$ grows away from~$\mathcal M$, and use this to get a bound on the rate at which~$d(\hat X_t, \mathcal M)$ vanishes as~$t \to T^-$.
This is our next lemma.

\begin{lemma}\label{l:rhoDecay}
There exists $d_0>0$ and constants $c_1, c_2 >0$ such that
for every $0 < t\le T$, and every $x \in \R^d$ such that $d(x,\cM)<d_0$,
we have
\begin{equation}\label{e:RhoTUpperLower}
\frac{c_1e^{-d(x,\cM)^2/(4t)}}{(4\pi t)^{(d-m)/2}} \le \rho_t(x)\leq  
\frac{c_2e^{-d(x, \cM)^2/(4t)}}{(4\pi t)^{(d-m)/2}}.
\end{equation}
\end{lemma}
We prove Lemma~\ref{l:rhoDecay} using heat kernel estimates in Section~\ref{s:rhoDecay}, below.
An immediate consequence of Lemma~\ref{l:rhoDecay} is that $V_t(x)$ indicates how close
the point $X_t$ is to the sub-manifold~$\cM$, in the following sense.   
\begin{corollary}
We have~$V_t(x) > 0$, and $V_t(x) \to +\infty$ if~$x \not \in \mathcal M$ 
and~$V_t(x) \to 0$ if~$x \in \mathcal M$, both as~$t \to T$. 
\end{corollary}

Lemmas~\ref{l:dXtHatMVanish} and~\ref{l:rhoDecay} quickly prove the quantitative convergence bound~\eqref{e:dXtHatMrate}.
To bound the velocity along trajectories we need to understand the asymptotic behavior of~$u_t(x)$ as~$t \to T^-$.
We do this using an identity that expresses~$u_t(x)$ as the ``Burgers' velocity'' towards a point in the convex hull of~$\mathcal M$.
Differentiating~\eqref{e:rhoTdef} and using~\eqref{e:dXExactScore} and the identity
\begin{equation}
  \grad \phi_{2t}(x) = -\frac{x}{2t} \phi_{2t}(x)\,,
\end{equation}
we see
\begin{equation}\label{e:uZeta}
  u_t(x)
    =-\frac{\gamma}{2 \rho_{T-t}(x)}\int_{\R^d}
      \paren[\bigg]{\frac{x-y}{T-t}}
      \phi_{2(T - t)}(x - y) \,\rho_0(dy)
    = -\frac{\gamma (x-\zeta_t(x))}{2(T-t)}
  \,,
\end{equation}
where~$\zeta_t$ is defined by
\begin{equation}\label{e:zetaDef}
  \zeta_t(x) \defeq
    \frac{1}{\rho_{T - t}(x)}
    \int_{\R^d} y \, \phi_{2(T - t)}(x - y) \, \rho_0(dy)
    \,.
\end{equation}
Identity~\eqref{e:uZeta} is a form of Tweedie's formula;
see, for example,~\cite{Efron11}.

Since $\zeta_t(x)$ is a weighted average of points on~$\mathcal M$, it 
belongs to the convex hull~$\conv(\mathcal M)$ of the sub-manifold $\cM$.
As~$t \to T^-$, this weighted average is dominated by the points on~$\mathcal M$  
that are closest to~$x$.
Thus, as~$t \to T^-$, $\zeta_t(x)$ should converge to the point on~$\mathcal M$ 
that is closest to~$x$, provided this point is unique.
A version of this observation can be found in~\cite{DeBortoli22,LuWangEA24,Pidstrigach22},
when~$x$ is in in a small tubular neighborhood around the manifold.
While this is enough to prove Theorem~\ref{t:relErrorScore}, it will not be 
sufficient to prove Theorem~\ref{t:supportRecoveryDiscrete}.
Indeed, to prove Theorem~\ref{t:supportRecoveryDiscrete} one needs to understand convergence of~$\zeta_t(x)$ as~$t \to T^-$ even when~$x$ is close to the cut-locus of~$\mathcal M$.

We will bound the rate at which~$\zeta_t(x)$ converges as~$t \to T^-$ in a \emph{stable projection zone}, consisting of all points in a bounded region whose distance from the cut-locus stays away from~$0$.
Explicitly, given~$r, R > 0$, define the ``good set''~$\mathcal G(r, R)$ by
\begin{equation}\label{e:GrR}
  \mathcal G(r, R) \defeq \set{x \in B(0, R) \st d(x, \cut(\mathcal M) ) > r }
  \,.
\end{equation}
The condition~$x \in B(0, R)$ is added for technical reasons.
We will eventually choose~$R$ large enough so that~$\mathcal M \subseteq B(0, R)$.
The next lemma provides the rate at which~$\zeta_t$ converges on~$\mathcal G(r, R)$ as~$t \to T$.

\begin{lemma}\label{l:zetaBdKTau-bis2}
Let $K, R>0$, $b\in (0, 1)$, $c \in [0, 1)$ be such that~$2b + c < 1$.
There exist constants~$C, \tau_0 > 0$, depending on~$b$, $c$, $K$, and $R$, 
such that for every~$x\in\cG(K\tau^c,R)$ and all~$\tau\in(0,\tau_0)$, we have
\begin{equation}\label{e:zetaAsTtoBigT}
  |\zeta_{T - \tau}(x) - x^*|
    \leq C \tau^b
  \,.
\end{equation}
Here, %
$x^* = \argmin_{z \in \mathcal M} \abs{z - x}$ is the unique point in~$\mathcal M$ 
that is closest to~$x$.
\end{lemma}
\begin{remark}
The proof of Lemma~\ref{l:zetaBdKTau-bis2} (specifically~\eqref{e:ZetatMinusTCase2}, below) 
shows that, in fact, we have the bound
  \begin{equation}\label{e:zetaLipschitz}
    \abs{ \zeta_{T - \tau}(x) - x^*} \leq C\sqrt{\tau},
  \end{equation}
  uniformly for~$x$ in a small tubular neighborhood around~$\mathcal M$,
  for all sufficiently small $\tau>0$.
  Using~\eqref{e:uZeta} this implies
  \begin{equation}\label{e:scoreTubeNbd}
    u_t(x) =- \frac{\gamma(x - x^*)}{2(T - t)} + v_t(x)\,,
  \end{equation}
for some function~$v_t(x)$ which satisfies
$|v_t(x)|\leq C(T-t)^{-1/2}$ uniformly in a small tubular
neighborhood around~$\mathcal M$, for $t$ sufficiently close to $T$.
A version of~\eqref{e:scoreTubeNbd} can also
be found in~\cite{DeBortoli22,LuWangEA24,Pidstrigach22}.
\end{remark}

We prove Lemma~\ref{l:zetaBdKTau-bis2} in Section~\ref{s:zetaConv}, below.
We now use Lemmas~\ref{l:dXtHatMVanish},~\ref{l:rhoDecay} 
and~\ref{l:zetaBdKTau-bis2}  to prove Theorem~\ref{t:relErrorScore}.

\begin{proof}[Proof of Theorem~\ref{t:relErrorScore}]
We first prove the upper bound on the distance $d(\hat X_t, \mathcal M)$ given 
in~\eqref{e:dXtHatMrate}.
Let~$\hat Z_t = V_t(\hat X_t)$ be as in~\eqref{e:zHatDef}.
By Lemma~\ref{l:dXtHatMVanish}, the process $\hat Z_t$ is a non-negative super-martingale 
and converges almost surely to an integrable limit $\hat Z_T$ as~$t \to T$,
and, in addition,~$d(\hat X_t, \mathcal M) \to 0$ as~$t \to T$.
Consequently, almost surely, there exists a random time~$\tau_1 < T$ 
such that for all $t \in (\tau_1, T)$, we have 
\begin{equation}
\frac{(T-t)^p}{[\rho_{T-t}(\hat X_t)]^{\beta}} \leq 1 + \hat Z_T
  \,.
\end{equation}
This implies
\begin{equation}\label{e:rhoHatLowerBoundGamma}
  \frac{p}{\beta}\log(T-t) - \log \rho_{T-t}(\hat X_t) \leq \frac{1}{\beta} \log(1 + \hat Z_T).
\end{equation}
Using Lemma~\ref{l:rhoDecay}, for any point $x$ in a tubular neighborhood of $\mathcal M$, we have
\begin{equation}\label{e:rhoLogDecayBoundGamma}
  \frac{d(x, \mathcal M)^2}{4(T-t)} + \frac{d-m}{2} \log(4\pi (T-t)) - \log c_2 \leq -\log \rho_{T-t}(x).
\end{equation}
Evaluating~\eqref{e:rhoLogDecayBoundGamma} at $x = \hat X_t$, 
with $t\in(\tau_1,T)$, and applying \eqref{e:rhoHatLowerBoundGamma} yields
\begin{equation}
  \paren[\Big]{\frac{p}{\beta} + \frac{d-m}{2}}\log(T-t)
  + \frac{d(\hat X_t, \mathcal M)^2}{4(T-t)}
  + \frac{d-m}{2} \log(4\pi) - \log c_2 \leq \frac{1}{\beta} \log(1 + \hat Z_T).
\end{equation}
Using~\eqref{26sep902}, this implies
\begin{equation}
\begin{aligned}
  \
  d(\hat X_t, \mathcal M)^2 &\leq 2((\gamma - 1)d - m)(T-t)|\log(T-t)|
    \\
    \label{e:distHatSquareBoundGamma}
       & + 4(T-t)
      \paren[\bigg]{
    \frac{\log(1+\hat Z_T)}{\beta} + \log c_2 - \frac{d-m}{2}\log(4\pi)
      }
      \,,
\end{aligned}
\end{equation}
which immediately yields~\eqref{e:dXtHatMrate}.
\smallskip

We now prove~\eqref{e:uHatConvRate}.
Since~$\mathcal M$ is~$C^2$, there exists~$d_0 > 0$ so that
$\mathcal T \defeq B(\mathcal M,d_0)$ defines a tubular
neighborhood around~$\mathcal M$
(see for instance~\cite[Chapter 4.5]{Hirsch94}).
Taking $d_0$ sufficiently small, we may assume
that~\eqref{e:zetaLipschitz} holds in~$\mathcal T$.
Since $d(\hat X_t,\mathcal M)\to0$ almost surely as $t\to T^-$,
we must have~$\hat X_t\in\mathcal T$
for all $t$ sufficiently close to $T$.
Now, using~\eqref{e:uZeta} and~\eqref{e:zetaLipschitz} gives
\begin{equation}\label{e:uHatPointwiseEstimateGamma}
\begin{aligned}
  \frac{2}{\gamma}|u_t(\hat X_t)|(T-t)
    &=|\zeta_t(\hat X_t)-\hat X_t|\\
    &\leq|\hat X_t^*-\zeta_t(\hat X_t)|
          +|\hat X_t-\hat X_t^*|\\
    &\leq d(\hat X_t,\mathcal M)+C\sqrt{T-t}.
\end{aligned}
\end{equation}
By~\eqref{e:RelErrorScore},
$|\hat u_t(\hat X_t)|\leq(1+\delta)|u_t(\hat X_t)|$.
This implies
\[
  \left(\frac{T-t}{|\log(T-t)|}\right)^{1/2}
  |\hat u_t(\hat X_t)|
  \leq
  \frac{\gamma(1+\delta)}{2}
  \left(
    \frac{d(\hat X_t,\mathcal M)}
         {\sqrt{(T-t)|\log(T-t)|}}
    +\frac{C}{\sqrt{|\log(T-t)|}}
  \right).
\]
Taking the upper limit as $t\to T^-$ and
using~\eqref{e:dXtHatMrate} gives~\eqref{e:uHatConvRate},
finishing the proof of Theorem~\ref{t:relErrorScore}.
\end{proof}

\subsection{The delta function example}\label{s:deltaExample}

We now consider the simple example where~$\rho_0 = \delta_0$ is the Dirac~$\delta$-measure at the origin.

\begin{proposition}\label{p:deltaExample}
  Suppose~$\mathcal M = \set{0}$ and~$\rho_0 = \delta_0$ is the Dirac~$\delta$-measure at the origin.
  Let~$X_t$, $u_t(x)$ be defined by~\eqref{e:dXExactScore}-\eqref{e:uDef}.
  Then, for any~$\gamma > 1$, the processes $X$ and~$u(X)$ are Gaussian with
  \begin{equation}\label{e:deltaExplicit}
    \E X_t = \frac{(T - t)^{\gamma / 2}}{T^{\gamma / 2}}X_0
    \,,\quad
    \lim_{t \to T} \frac{\operatorname{Var}(X_t)}{T - t} = 2
    \,,
  \end{equation}
  and
  \begin{equation}\label{e:deltaExplicitU}
    \lim_{t \to T} (T - t)\operatorname{Var}(u_t(X_t)) = \frac{\gamma^2}{2}
    \,.
  \end{equation}
\end{proposition}

As mentioned in Remark~\ref{r:convRateSharp}, this implies that the bounds~\eqref{e:dXtHatMrate} 
and~\eqref{e:uHatConvRate} in Theorem~\ref{t:relErrorScore} are at worst 
sub-optimal by a~$\abs{\log(T - t)}$ factor.

\begin{proof}[Proof of Proposition~\ref{p:deltaExample}]
When~$\rho_0 = \delta_0$, $\rho_t = \phi_{2t}$ (where~$\phi$ is where~$\phi$ is defined in~\eqref{e:rhoTdef}).
The velocity~$u$ is precisely

\begin{equation}\label{e:uDelta}
u_t(x)=\frac{\gamma\nabla\phi_{2(T-t)}(x)}
             {\phi_{2(T-t)}(x)}
       =-\frac{\gamma x}{2(T-t)}.
       \end{equation}
The denoising diffusion~\eqref{e:dXExactScore} now takes the form 
\begin{equation}\label{25apr2204}
dX_t=-\frac{\gamma}{2(T-t)}X_tdt+ \sigma \, dB_t,
\end{equation}
which is equivalent to
\begin{equation}
d\big((T-t)^{-\gamma/2}X_t\big)= \sigma (T-t)^{-\gamma/2}dB_t
\,.
\end{equation}
Integrating shows
\begin{equation}
(T-t)^{-\gamma/2}X_t-T^{-\gamma/2} X_0= \sigma \int_0^t(T-s)^{-\gamma/2}dB_s
\,,
\end{equation}
and hence
\begin{equation}\label{26sep1104}
  X_t=\frac{(T-t)^{\gamma/2}}{T^{\gamma/2}} X_0+\sqrt{2 (\gamma-1)}Y_t,
  \qquad
  Y_t= 
  (T-t)^{\gamma/2}\int_0^t(T-s)^{-\gamma/2}dB_s
  \,.
\end{equation}
Taking the expectation in~\eqref{26sep1104}, we obtain the first identity in~\eqref{e:deltaExplicit}. 
The variance of the mean-zero Gaussian process $Y_t$ is  
 \begin{equation}
\begin{aligned}
\Em(Y_t^2)
  &=(T-t)^{\gamma}\int_0^t(T-s)^{-\gamma}ds=(T-t)^{\gamma}\int_{T-t}^T s^{-\gamma}ds
\\
  &=\frac{(T -t)^\gamma}{\gamma-1} \paren[\Big]{
      \frac{1}{(T-t)^{\gamma-1}}-\frac{1}{T^{\gamma-1}}
    }
  \sim \frac{T-t}{\gamma - 1},\quad\text{as } t\to T^-\,.
\end{aligned}
 \end{equation}
This gives
\begin{equation}
\frac{\operatorname{Var}(X_t)}{T-t}
  =\frac{2(\gamma-1)}{T-t}\Em(Y_t^2)\to 2, \quad\text{ as } t\to T^-.
\end{equation}
which proves the second equality in~\eqref{e:deltaExplicit}.
Combining this with~\eqref{e:uDelta} yields~\eqref{e:deltaExplicitU}, concluding the proof.
\end{proof}

\iffalse
%
%
%
%
%
%
%
%
%
%
%
%
%
%
%
%
%
%
%
%
%
%
%
%
%
%
%
%
%
\fi

\section{Convergence of the Discrete Algorithm}\label{s:discreteConv}

In this section, we prove Theorem~\ref{t:supportRecoveryDiscrete}, a generalization 
of Theorem~\ref{t:supportRecoveryDiscreteIntro} 
that allows any~$\gamma > 1$, and obtain more precise convergence rate.
Before doing this, we briefly explain the main idea behind the proof.
Suppose for simplicity we are working with the exact score, that is,~$\delta = 0$.
In this case, the discretization~\eqref{e:EMhatX}   becomes
\begin{equation}
  X_{n+1} = X_n + u_{n}(X_n) \Delta t_n + \sqrt{ 2 \Delta t_n} \, \xi_{n+1}\,.
\end{equation}
Using the representation~\eqref{e:uZeta} of $u_t(x)$ with $\gamma=2$, this is equivalent to
\begin{equation}\label{e:XnExact}
    X_{n+1}
      = \alpha X_n + (1 - \alpha) \zeta_{t_n}(X_n) + \sqrt{2 \Delta t_n} \, \xi_{n+1}
    \,.
\end{equation}
In other words, $X_{n+1}$ is a convex combination of~$X_n$ and~$\zeta_{t_n}(X_n)$.
The strategy we use to prove Theorem~\ref{t:supportRecoveryDiscrete} is to estimate the chance that~$X_n$ 
avoids the cut-locus, using Lemma~\ref{l:ballVolumeFractions2},
then use Lemma~\ref{l:zetaBdKTau-bis2} to ensure
that~$d(\zeta_{t_n}(X_n), \mathcal M)$ is small because $X_n$ is not in the cut-locus, 
and finally use~\eqref{e:XnExact} to show
that this implies that~$d(X_n, \mathcal M)$ decreases exponentially in~$n$.
When~$\delta > 0$ and~$\gamma \neq 2$ the details are more involved, but the essential idea is unchanged.
Once again, we mention that, somewhat surprisingly, the most difficult part of the proof is 
Lemma~\ref{l:ballVolumeFractions2}. This proof is postponed until Sections~\ref{s:ballVolFractions2}
and~\ref{s:cutTubeBd}. 
The rest of the proof of Theorem~\ref{t:supportRecoveryDiscrete} is contained in this section,
though it relies on some of the lemmas used in the proof of Theorem~\ref{t:relErrorScore}
that are also proved later in the paper.

\subsection{The convergence result for any \texorpdfstring{$\gamma > 1$}{gamma > 1}}

As before, let~$\delta \in [0, 1)$ and~$\hat s_t$ be an approximate score function that satisfies 
the relative error estimate~\eqref{e:RelErrorScore}.
Let~$\gamma > 1$ and~$\hat u_t$ be the associated velocity field defined in~\eqref{e:uHatDef}.
In this case the~$\alpha$ in~\eqref{e:XnExact} will get replaced with
\begin{equation}\label{e:betaDef}
  \beta \defeq 1 - \frac{\gamma}{2}(1 - \alpha)
  \,.
\end{equation}
Define~$\beta$ by~\eqref{e:betaDef}, and assume
\begin{equation}\label{e:betaTildeDef}
  \tilde \beta \defeq \abs{\beta} + \delta (1 - \beta) < 1
  \,.
\end{equation}

The generalization of Theorem~\ref{t:supportRecoveryDiscreteIntro} is as follows.

\begin{theorem}\label{t:supportRecoveryDiscrete}
Let~$m < d$,~$\mathcal M$ be an~$m$-dimensional closed~$C^\infty$ sub-manifold of~$\R^d$, 
and~$\rho$ be as in Theorem~\ref{t:relErrorScore}.
Assume, as above, that~$\delta \in [0, 1)$ satisfies~\eqref{e:betaDef}-\eqref{e:betaTildeDef},~$\hat s_t$ 
is an approximate score 
function that satisfies~\eqref{e:RelErrorScore}, and~$\hat u_t$ is defined by~\eqref{e:uHatDef}.
Let $\hat X_{0} \in\R^d$ be an arbitrary initial point and inductively define~$\hat X_n$ by~\eqref{e:tnDef}--\eqref{e:EMhatX}.
Given $b\in(0,1)$, define
\begin{equation}\label{e:tildeAlphaBetaDef}
  \tilde \alpha_b = \max \set{\alpha^b, \tilde \beta}\,.
\end{equation}
\begin{enumerate}
  \item
For all~$b < 1/6$  
there exists a constant~$C = C(\alpha, b, \gamma, \mathcal M)<+\infty$ such that
    \begin{equation}\label{e:EdHatXnMExpDecay}
      \E d(\hat X_n, \mathcal M) \leq
	  \begin{dcases}
	    (d(\hat X_0, \mathcal M) + C) \tilde \alpha_b^n,
	      & \text{ if } \alpha^b \neq \tilde \beta\,,
	    \\
	    (d(\hat X_0, \mathcal M) + C n) \tilde \alpha_b^n,
	      & \text{ if } \alpha^b = \tilde \beta
	    \,.
	  \end{dcases}
      \,.
    \end{equation}

  \item
For all~$b < 1/4$  
there exists a random constant~$D = D(\alpha, \gamma, b, X_0, \mathcal M)<+\infty$ such that,
almost surely,
    \begin{equation}\label{e:dHatXnMExpDecay}
      d(\hat X_n, \mathcal M)
	\leq
	  \begin{dcases}
	    D \tilde \alpha_b^n, & \text{ if } \alpha^b \neq \tilde \beta\,,
	    \\
	    D (1 + n) \tilde \alpha_b^n, & \text{ if } \alpha^b = \tilde \beta
	    \,,
	  \end{dcases}
    \end{equation}
    for every~$n \in \N$. 
\end{enumerate}
\end{theorem}

Theorem~\ref{t:supportRecoveryDiscrete} is a generalization of 
Theorem~\ref{t:supportRecoveryDiscreteIntro}.
Indeed, if~$\gamma = 2$, then~$\beta = \alpha \in (0, 1)$, and~\eqref{e:betaTildeDef} 
is automatically satisfied.
Moreover, $\tilde \alpha_b$ in~\eqref{e:tildeAlphaBetaDef} is smaller than~$\bar \alpha_b$ in~\eqref{e:barAlphaB}, and so~\eqref{e:EdHatXnMExpDecay} and~\eqref{e:dHatXnMExpDecay} will immediately imply~\eqref{e:EdXnHatMGamma2} and~\eqref{e:dXnHatMGamma2} respectively.

In general, the assumption~\eqref{e:betaTildeDef} provides a restriction on the relative error~$\delta$, or on the time step ratio~$\alpha$.
Explicitly, if
\begin{equation}
  \alpha > 1 - \frac{2}{\gamma}
  \,,
\end{equation}
then $\beta\in(0,1]$ and~\eqref{e:betaTildeDef} is satisfied for any~$\delta \in [0, 1)$.
If instead we only have
\begin{equation}
  \alpha > 1 - \frac{4}{\gamma}
  \,,
\end{equation}
then~\eqref{e:betaTildeDef} is only satisfied for any~$\delta \in [0, (1 - \abs{\beta}) / (1 - \beta) )$.
For smaller values of~$\alpha$, the assumption~\eqref{e:betaTildeDef} is not satisfied.

\subsection{Boundedness of iterates}

In order to prove Theorem~\ref{t:supportRecoveryDiscrete}, we need to first ensure that~$\hat X_n$ stays in a bounded region.

\begin{lemma}\label{l:XnHatBounded}
Let~$\hat X_n$ be defined by~\eqref{e:EMhatX} and let
$\alpha_* = \max\set{\sqrt{\alpha}, \tilde \beta}$ if $\sqrt{\alpha} \neq \tilde \beta$, or be any number in~$(\tilde \beta, 1)$ if $\sqrt{\alpha} = \tilde \beta$.
  There exist constants~$R, C_1, C_2 > 0$ depending on $\alpha, \tilde \beta, \delta, \sigma, T, d, \alpha_*$, and~$\mathcal M$, such that for any $r > 0$,
  \begin{equation}\label{e:XnHatTailBound}
    \P\paren[\Big]{ \abs{\hat X_n} > R + \tilde \beta^n \abs{\hat X_0} + C_1 \alpha_*^n + r }
      \leq 2 \exp\paren[\bigg]{-\frac{r^2}{C_2 \alpha_*^{2n}}},\quad\text{for all $n\ge 1$.}
  \end{equation}
\end{lemma}
Recall that the parameter $\tilde\beta$ is defined in in~\eqref{e:betaTildeDef} and that $\tilde\beta<1$.  
Since the proof of Theorem~\ref{t:supportRecoveryDiscrete} relies on ideas used in the proof of Lemma~\ref{l:XnHatBounded}, we prove Lemma~\ref{l:XnHatBounded} here.

\begin{proof}[Proof of Lemma~\ref{l:XnHatBounded}]
Define
\begin{equation}\label{e:zetaHatDef}
  \hat \zeta_t(x) \defeq x + 2 (T - t) \hat s_{T - t}(x)
    \,.
\end{equation}
and observe, using the relations~\eqref{e:tnDef} and the definition~\eqref{e:betaDef} of $\beta$:
\begin{equation}
\begin{aligned}
  \hat X_{n+1}
    &= \hat X_n
      + \frac{\gamma \Delta t_n}{2(T - t_n)} (\hat \zeta_n(\hat X_n) - \hat X_n)
      + \sigma \sqrt{\Delta t_n} \, \xi_{n+1}
  \\
  &=\hat X_n
      + \frac{\gamma (1-\alpha)}{2 } (\hat \zeta_n(\hat X_n) - \hat X_n)
      + \sigma \sqrt{\Delta t_n} \, \xi_{n+1}
  \\
    \label{e:hatXnp11}
    &= \beta \hat X_n + (1 - \beta) \hat \zeta_n(\hat X_n) + \sigma \sqrt{\Delta t_n} \, \xi_{n+1}
    \,.
\end{aligned}
\end{equation}
While
\[
1 - \beta=\frac{\gamma}{2}(1-\alpha),
\]
is always strictly positive, $\beta$ itself may be negative if $\gamma$ is large.
Note that this does not happen for~$\gamma=2$, in which case $\beta=\alpha\in(0,1)$.
Thus, in general, the first two terms in the right side of~\eqref{e:hatXnp11}
do not necessarily represent a convex linear combination of $\hat X_n$
and $\hat\zeta_n$.
Nevertheless, this expression will still be
good enough for our purposes, as long as~$\tilde \beta < 1$.

We first replace~$\hat \zeta_n$ in the right side of~\eqref{e:hatXnp11} by~$\zeta_n$.
For this, note that 
\eqref{e:uHatDef}, \eqref{e:RelErrorScore}, \eqref{e:uZeta} and~\eqref{e:zetaHatDef} imply
\begin{equation}
\begin{aligned}
  \abs{\hat \zeta_n(x) - \zeta_n(x)}
    &= 2 (T - t_n) \abs{\hat s_{T - t_n}(x) - s_{T - t_n}(x)}
  \\
  \label{26sep2304}
    &\leq 2 \delta (T - t_n) \abs{s_{T - t_n}(x)}
    \leq \delta \abs{\zeta_n(x) - x}
  \,,\quad\text{for any~$x \in \R^d$.}
\end{aligned}
\end{equation}
Using \eqref{26sep2304} and the triangle inequality in~\eqref{e:hatXnp11} gives
\begin{equation}
\begin{aligned}
  \abs{\hat X_{n+1}}
    &\leq \abs{\beta} \abs{\hat X_n} + (1 - \beta) \paren[\big]{ \abs{\hat \zeta_n(\hat X_n) - \zeta_n(\hat X_n)} + \abs{\zeta_n(\hat X_n)} } + \sigma \sqrt{\Delta t_n} \, \abs{\xi_{n+1}}
  \\
    &\leq \paren[\big]{ \abs{\beta} + \delta(1-\beta) } \abs{\hat X_n} + (1 - \beta)(1 + \delta)\abs{\zeta_n(\hat X_n)} + \sigma \sqrt{\Delta t_n} \, \abs{\xi_{n+1}}
  \\
    \label{e:HatXnTmp1}
    &\leq \tilde \beta \abs{\hat X_n}
      + (1 - \beta) (1 + \delta) D_C
      + \sigma \sqrt{\Delta t_n} \, \abs{\xi_{n+1}}.
\end{aligned}
\end{equation}
Here, we have set 
\begin{equation}
  D_C \defeq \max_{z \in \conv(\mathcal M)} \abs{z},
\end{equation}
and used the fact that~$\zeta_t$ always takes values in~$\conv(\mathcal M)$, as well as
the definition~\eqref{e:tildeAlphaBetaDef} of~$\tilde\beta$. 
Iterating inequality~\eqref{e:HatXnTmp1} from $0$ to $n$ yields
\begin{equation}
\begin{aligned}
  \abs{\hat X_n} &\leq
    \tilde \beta^n \abs{\hat X_0}
    + (1 - \beta)(1 + \delta) D_C \sum_{k=0}^{n-1} \tilde \beta^{n-1-k}
    + \sigma \sum_{k=0}^{n-1} \tilde \beta^{n-1-k} \sqrt{\Delta t_k} \abs{\xi_{k+1}}
  \\
  \label{e:absHatXn1}
     &\leq
    \tilde \beta^n \abs{\hat X_0}
    + R
    + Z_n
    \,,
\end{aligned}
\end{equation}
where
\begin{equation}
  R \defeq \frac{(1 - \beta)(1 + \delta) D_C}{1 - \tilde \beta},
  \end{equation}
and
\begin{equation}\label{26sep2306}
  Z_n \defeq \sigma \sum_{k=0}^{n-1} \tilde \beta^{n-1-k} \sqrt{\Delta t_k} \abs{\xi_{k+1}}
  \,.
\end{equation}

We now bound~$Z_n$ via a standard Gaussian concentration inequality.
Since
\[
\E\abs{\xi_k} \leq \sqrt{d},
\]
we have
\begin{equation}\label{e:EZn}
  \E Z_n \overset{\eqref{e:tnDef}}{\leq} \sigma \sqrt{d T (1 - \alpha)} \sum_{k=0}^{n-1} \tilde \beta^{n-1-k} \alpha^{k/2} \leq C_1 \alpha_*^n
  \,,
\end{equation}
for some constant $C_1 = C_1 (\alpha, \tilde \beta, \sigma, T, d, \alpha_*)$.
Also, note that~\eqref{26sep2306} can be written as 
\[
Z_n = L_n(\xi_1, \dots, \xi_n),
\]
with a Lipschitz function~$L$.
Its Lipschitz constant can be estimated, using~\eqref{e:tnDef}, as 
\begin{equation}
  \norm{L_n}_{\text{Lip}}^2 = \sigma^2 \sum_{k=0}^{n-1} \tilde \beta^{2(n-1-k)} \Delta t_k
    = \sigma^2 T (1 - \alpha) \sum_{k=0}^{n-1} \tilde \beta^{2(n-1-k)} \alpha^k
    \leq C_2 \alpha_*^{2n}
    \,,
\end{equation}
with some constant $C_2 > 0$.
By the Gaussian concentration inequality~\cite[Theorem~5.2.2]{Vershynin18}, and making~$C_2$ larger if necessary this implies
\begin{equation}
  \P( Z_n > \E Z_n + r )
    \leq 2 \exp\paren[\bigg]{-\frac{r^2}{C L_n^2}}
    \leq 2 \exp\paren[\bigg]{-\frac{r^2}{C_2 \alpha_*^{2n}}}
  \,.
\end{equation}
Combining this with~\eqref{e:absHatXn1} and \eqref{e:EZn}, 
yields~\eqref{e:XnHatTailBound}, concluding the proof.
\end{proof}

\subsection{The expectation bound in Theorem~\ref{t:supportRecoveryDiscrete}}

We now prove the expectation bound~\eqref{e:EdHatXnMExpDecay}, which is the first assertion in Theorem~\ref{t:supportRecoveryDiscrete}.

\begin{proof}[Proof of the first assertion in Theorem~\ref{t:supportRecoveryDiscrete}]
Let~$\hat X_n^*$ denote the projection of~$\hat X_n$ onto~$\mathcal M$.
Since~$\hat X_n$ is a continuous random variable and~$\cut(\mathcal M)$ has Lebesgue measure~$0$, $\hat X_n^*$ is defined almost surely.
Define~$\hat \zeta_t$ by~\eqref{e:zetaHatDef}, as in the proof of Lemma~\ref{l:XnHatBounded}.
Subtracting~$\hat X_n^*$ from both sides of~\eqref{e:hatXnp11} yields
\begin{equation}\label{e:XnXnStarGamma}
  \hat X_{n+1} - \hat X_n^*
    = \beta(\hat X_n - \hat X_n^*) + (1-\beta)(\hat \zeta_n(\hat X_n) - \hat X_n^*)
      + \sigma \sqrt{\Delta t_n} \, \xi_{n+1}
    \,.
\end{equation}
Using the triangle inequality and taking expectations gives
\begin{equation}\label{e:EXn1XStar1Gen}
  \E{\abs{\hat X_{n+1} - \hat X_n^*}}
    \leq \abs{\beta} \E\abs{\hat X_n - \hat X_n^*}
      + (1-\beta)\E \abs{\hat \zeta_n(\hat X_n) - \hat X_n^*}
      + C \alpha^{n/2}
    \,.
\end{equation}

We first estimate the second term in the right side of~\eqref{e:EXn1XStar1Gen} that we split as
\begin{equation}\label{26sep2308}
  \E \abs{\hat \zeta_n(\hat X_n) - \hat X_n^*}
    = \E \big[\one_{A_n} \abs{\hat \zeta_n(\hat X_n) - \hat X_n^*}]
      + \E\big[ \one_{A_n^c} \abs{\hat \zeta_n(\hat X_n) - \hat X_n^*}\big],
      \,
\end{equation}
where
\begin{equation}\label{e:AnDef}
  A_n \defeq \set{\hat X_n \not\in \mathcal G( \alpha^{c n}, R_1 ) }
  \,,
\end{equation}
for some~$R_1 \geq 1$ that will be chosen shortly. Recall that the good set $\cG(r,R)$ is defined 
in~\eqref{e:GrR}. That is, on the event $A_n$, the point $\hat X_n$
is very close to the cut-locus.  The parameter~$c$ in the definition of $A_n$ is chosen so that  
\begin{equation}\label{e:bcRel}
  b + \frac{1}{2} < c < 1 - 2b.
\end{equation}
It is the balance between the lower and upper bounds on $c$ in~\eqref{e:bcRel} that necessitates the restriction~$b\in(0,1/6)$.  
We will see why both of these bounds are needed in the proof below.

To bound the first term in the right side of~\eqref{26sep2308},  note that
for any~$x \in \R^d$, \eqref{e:RelErrorScore}, \eqref{e:uZeta} and~\eqref{e:zetaHatDef} imply
\begin{equation}
\begin{aligned}
  \abs{\hat \zeta_n(x) - \zeta_n(x)}
    &= 2 (T - t_n) \abs{\hat s_{T - t_n}(x) - s_{T - t_n}(x)}
  \\
    \label{e:hatZetan12}
    &\leq 2 \delta (T - t_n) \abs{s_{T - t_n}(x)}
    \leq \delta \abs{\zeta_n(x) - x}
  \,,
\end{aligned}
\end{equation}
so that 
\begin{equation}\label{e:hatZetaNXnStar}
  \abs{\hat \zeta_n(\hat X_n) - \hat X_n^*}
    \leq \delta \abs{\hat X_n} + \delta \abs{\zeta_n(\hat X_n)} + \abs{\zeta_n(\hat X_n)-\hat X_n^*}
    \leq C(1 + \abs{\hat X_n})\,.
\end{equation}
The last inequality above was obtained by using the fact 
that~$\hat X_n^* \in \mathcal M$, $\zeta_n(\hat X_n) \in \conv(\mathcal M)$, and that~$\mathcal M$ is bounded.

Using the definitions~\eqref{e:GrR} and \eqref{e:AnDef},  we see that 
\begin{equation}\label{e:OneAn1}
  \E \big[\one_{A_n} \paren{1 + \abs{\hat X_n}\big]}
    \leq 2 \E\paren[\big]{ \one_{\set{\abs{\hat X_n} \ge R_1}} \abs{\hat X_n} }
      + (1 + R_1) \P( d(\hat X_n, \cut_{R_1}(\mathcal M)) \le \alpha^{c n} )
  \,.
\end{equation}
We bound the first term in the right side above using Lemma~\ref{l:XnHatBounded} by taking
$R_1$ sufficiently large.
Let $R, C_1, C_2$ be as in that lemma and 
choose
\[
R_1 = R + \abs{\hat X_0} + C_1 + 1.
\] 
By Lemma~\ref{l:XnHatBounded} we have
\begin{equation}
\begin{aligned}
  \E\paren[\big]{ \one_{\set{\abs{\hat X_n} \ge R_1}} \abs{\hat X_n} }
     &= R_1 \P(\abs{\hat X_n}\ge R_1)
      + \int_0^\infty \P(\abs{\hat X_n} \ge R_1 + r) \, dr
  \\
    &\le 2R_1\exp\paren[\bigg]{-\frac{1}{C_2\alpha_*^{2n}}}
      + 2\int_0^\infty
      \exp\paren[\bigg]{-\frac{(r+1)^2}{C_2\alpha_*^{2n}}}\,dr
  \\
    \label{e:XnHatTailInt}
    &\le C\exp\paren[\bigg]{-\frac{1}{C_2\alpha_*^{2n}}},
    \end{aligned}
\end{equation}
which decays super-exponentially.

For the second term in~\eqref{e:OneAn1}, we use~\eqref{e:bcRel} and choose~$\epsilon$ so that
\begin{equation}\label{e:epsilonDef}
  0 < \epsilon < c - b - \frac{1}{2}
  \,,
\end{equation}
condition on~$\hat X_{n-1}$, and apply Lemma~\ref{l:ballVolumeFractions2} to obtain
\begin{equation}
\begin{aligned}
  \P( d(\hat X_n, \cut_{R_1}(\mathcal M)) \le \alpha^{c n} )
    &= \E \P \paren[\big]{ d(\hat X_n, \cut_{R_1}(\mathcal M)) \le \alpha^{c n} \given \hat X_{n-1} }
  \\
    &\leq \sup_{x \in \R^d} \P\paren[\big]{ d(\mathcal N(x, \paren{\sigma^2 \Delta t_{n-1}} I), \cut_{R_1}(\mathcal M)) \le \alpha^{cn} }
  \\
    \label{e:dHatXnCutM}
    &\leq \frac{C \alpha^{c n}}{ (\sigma^2 \Delta t_{n-1})^{\frac{1}{2} + \epsilon} }
    = C \alpha^{(c - \frac{1}{2} - \epsilon) n}
  \,.
\end{aligned}
\end{equation}
Using~\eqref{e:XnHatTailInt} and~\eqref{e:dHatXnCutM} in~\eqref{e:OneAn1}, combining it with~\eqref{e:hatZetaNXnStar} and dominating the right  
side of~\eqref{e:XnHatTailInt} by that of~\eqref{e:dHatXnCutM} implies
\begin{equation}\label{e:EZetaN1Gen}
  \E \big[\one_{A_n} \abs{\hat \zeta_n(\hat X_n) - \hat X_n^*}\big]
  \leq C \Em \big[\one_{A_n} (1 + \abs{\hat X_n})]
  \leq C \alpha^{(c - \frac{1}{2} - \epsilon) n}
  \,.
\end{equation}

For the second term in the right side of~\eqref{26sep2308}, we note that
on the event~$A_n^c$, Lemma~\ref{l:zetaBdKTau-bis2} applies. We also recall from~\eqref{e:tnDef} that
$T-t_n=T\alpha^n$. Together, these give 
\begin{equation}\label{e:hatZetan11}
\one_{A_n^c} \abs{\zeta_n(\hat X_n) - \hat X_n^*}
  \leq C \one_{A_n^c} \alpha^{b n}\,.
\end{equation}
Let us note that the upper bound on $c$ in~\eqref{e:bcRel} is needed here so that we can use Lemma~\ref{l:zetaBdKTau-bis2}.  
Thus, we have
\begin{equation} 
\begin{aligned}
  \one_{A_n^c}\abs{\hat \zeta_n(\hat X_n) - \hat X_n^*}
    &\leq \one_{A_n^c} \paren[\big]{
	\abs{\hat \zeta_n(\hat X_n) - \zeta_n(\hat X_n)}
	 + \abs{\zeta_n(\hat X_n) - \hat X_n^*}
      }
  \\
    &\overset{\mathclap{\eqref{e:hatZetan12}, \eqref{e:hatZetan11}}}{\leq}\quad~
      \one_{A_n^c} \paren[\big]{
	\delta \abs{\zeta_n(\hat X_n) - \hat X_n}
	 + C \alpha^{b n}
      }
  \\
    \label{e:ZetaHatNTmp1}
    &\overset{\mathclap{\eqref{e:hatZetan11}}}{\leq}\;
      \one_{A_n^c} \paren[\big]{
	\delta \abs{\hat X_n - \hat X_n^*}
	 + C \alpha^{b n}
      }
    \,.
\end{aligned}
\end{equation}
Taking the expectation, we obtain
\begin{equation}\label{e:EZetaN2Gen}
  \E \one_{A_n^c} \abs{\hat \zeta_n(\hat X_n) - \hat X_n^*}
  \leq C \alpha^{bn} + \delta \E \abs{\hat X_n - \hat X_n^*}
  \,.
\end{equation}
Combining~\eqref{e:EZetaN1Gen} and~\eqref{e:EZetaN2Gen} implies the following expectation bound on the second term
in~\eqref{e:EXn1XStar1Gen}:
\begin{equation}\label{e:ZetaNXStarGen}
  \E \abs{\hat \zeta_n(\hat X_n) - \hat X_n^*}
    \leq \delta \E \abs{\hat X_n - \hat X_n^*} + C \paren[\big]{
      \alpha^{(c - \frac{1}{2} - \epsilon)n}
      + \alpha^{bn}
    }
    \overset{\eqref{e:epsilonDef}}{\leq}
    \delta \E \abs{\hat X_n - \hat X_n^*} + C \alpha^{bn}
  \,.
\end{equation}

We now use~\eqref{e:ZetaNXStarGen} in~\eqref{e:EXn1XStar1Gen} to get
\begin{equation}
\begin{aligned}
\E [d(\hat X_{n+1}, \mathcal M)]
  &\leq \E\abs{\hat X_{n+1} - \hat X_n^*}
    \leq (\abs{\beta} + (1 - \beta) \delta) \E \abs{\hat X_n - \hat X_n^*}
      + C (\alpha^{b n} + \alpha^{n/2})
  \\
  \label{e:EdXnHatM1}
  &\overset{\mathclap{\eqref{e:betaTildeDef}}}{\leq}~ \tilde \beta
    \E [d(\hat X_{n}, \mathcal M)]
    + C \alpha^{b n}
  \,.
\end{aligned}
\end{equation}
In the last step above, we used the fact that~$b < 1/6$ to dominate~$\alpha^{n/2}$ by~$\alpha^{b n}$.

Setting 
\[
d_n=\E [d(\hat X_n, \mathcal M)],
\]
we re-write~\eqref{e:EdXnHatM1} as 
\begin{equation}\label{e:dnRec}
  d_{n+1} \leq \tilde \beta d_n + C \alpha^{b n}
  \,.
\end{equation}
Iterating this implies
\begin{equation}
  d_{n} \leq \tilde \beta^n d_0 + C \sum_{k = 0}^{n - 1} \tilde \beta^{n - 1 - k} \alpha^{b k}
  \,.
\end{equation}
If~$\tilde \beta \neq \alpha^b$ this implies
\begin{equation}
  d_n \leq \tilde \beta^n d_0 + \frac{C(\tilde \beta^n - \alpha^{b n})}{\tilde \beta - \alpha^b}
    \,.
\end{equation}
On the other hand, if~$\tilde \beta = \alpha^b$ then we obtain
\begin{equation}
  d_n \leq \tilde \beta^n d_0 + C n \tilde \beta^{n - 1}
    = \paren[\Big]{d_0 + \frac{C n}{\tilde \beta}} \tilde \beta^n
    \,.
\end{equation}
This immediately implies~\eqref{e:EdHatXnMExpDecay}.
\end{proof}
\subsection{The almost sure bound in Theorem~\ref{t:supportRecoveryDiscrete}}

The proof of~\eqref{e:dHatXnMExpDecay} is similar, but uses a Borel--Cantelli argument to almost surely stay away from the cut-locus.

\begin{proof}[Proof of the second assertion in Theorem~\ref{t:supportRecoveryDiscrete}]
In this case know~$0<b < 1/4$.
This allows us to choose~$c$ and~$\epsilon$ so that
\begin{equation}\label{26sep2502}
  \frac{1}{2} < c < 1 - 2b
  \quad\text{and}\quad
  0 < \epsilon < c - \frac{1}{2}
  .
\end{equation}
The upper bound on $c$ in~\eqref{26sep2502} still allows us  to 
apply Lemma~\ref{l:ballVolumeFractions2} to  argue nearly exactly as in~\eqref{e:OneAn1}-\eqref{e:EZetaN1Gen} to arrive at
\begin{equation}\label{26sep2504}
  \Pm \big[A_n\big]
  \leq C \alpha^{(c - \frac{1}{2} - \epsilon) n}
  \,.
\end{equation}
Since~$c > 1/2 + \epsilon$, we deduce that 
\[
\sum \P(A_n) < \infty.
\]
Thus,  by the Borel--Cantelli lemma, there exists a random~$N \in \N$ such that
\begin{equation}\label{e:XnNinBAlphaCn}
  \hat X_n \in \mathcal G(\alpha^{c n}, R)
  \quad\text{for all } n \geq N
  \,,
  ~\text{almost surely}
  \,.
\end{equation}
For all~$n \geq N$, the same argument proving~\eqref{e:ZetaHatNTmp1} shows
\begin{equation}\label{e:BC1}
  \abs{\hat \zeta_n(\hat X_n) - \hat X_n^*}
    \leq
	\delta \abs{\hat X_n - \hat X_n^*}
	 + C \alpha^{b n}
    \,.
\end{equation}
Also, using the Borel--Cantelli lemma again (and making the random~$N$ larger if necessary) we know that almost surely
\begin{equation}\label{e:BC2}
  \abs{\xi_n} \leq C \sqrt{\ln n}
    \quad \text{for all } n \geq N
    \,.
\end{equation}

Using~\eqref{e:BC1} and~\eqref{e:BC2} in~\eqref{e:XnXnStarGamma} implies
\begin{equation}\label{e:dHatXnMBCFinal}
  d(\hat X_n, \mathcal M)
    \leq \tilde \beta d(\hat X_n, \mathcal M)
      + C (\alpha^{b n} +  \sqrt{\ln n} \, \alpha^{n/2})
    \leq \tilde \beta d(\hat X_n, \mathcal M)
      + C \alpha^{b n}
  \,.
\end{equation}
We see that~$d_n = d(\hat X_n, \mathcal M)$  satisfies the exact same recurrence inequality~\eqref{e:dnRec},
except now the relation is only valid for all~$n \geq N$.
Solving this recurrence relation as before, but starting at the random time~$N$, and bounding the first~$N$ terms by a (random) constant yields~\eqref{e:dHatXnMExpDecay} as claimed.
This finishes the proof of the second assertion in~Theorem~\ref{t:supportRecoveryDiscrete}. 
\end{proof}

\section{A Lyapunov function for \texorpdfstring{$\hat X_t$}{hat X t}: 
the proof Lemma~\ref{l:dXtHatMVanish}}
\label{s:lyapunov}

\subsection{Super-solution inequalities} 

We now prove Lemma~\ref{l:dXtHatMVanish} and show that~$V_t(x)$ is a Lyapunov function for~$\hat X_t$.
We begin by proving two inequalities. Let~$\mathcal L_t$, $\hat{\mathcal L}_t$ be the operators whose action on test functions~$f$ is defined by
  \begin{equation}
    \mathcal L_t f = u_t \cdot \grad f + (\gamma - 1) \lap f,\qquad
    \hat{\mathcal L}_t f = \hat u_t \cdot \grad f + (\gamma - 1) \lap f
    \,.
  \end{equation}
\begin{lemma}\label{l:genNegative}
   If~$V_t(x)$ is defined by~\eqref{e:Vdef} then
  \begin{equation}
    \label{e:dtPlusLBound}
    (\partial_t + \mathcal L_t) V_t
      \leq
	- \beta (1 - \beta) (\gamma - 1)
	  \abs{s_{T - t}}^2
	  V_t,
\end{equation}
and
\begin{equation}
    \label{e:dtPlusLhatBound}
    (\partial_t + \hat{\mathcal L}_t) V_t
      \leq
	-\beta \paren[\big]{(1 - \beta) (\gamma - 1) - \gamma \delta} \abs{s_{T-t}}^2 V_t 
	\,.
  \end{equation}
\end{lemma}
\begin{proof}[Proof of Lemma~\ref{l:genNegative}]
Recall that~$\rho_t(x)$ satisfies the heat equation
\begin{equation}\label{e:heat}
\partial_t \rho_t =\Delta\rho_t\,,
\quad x\in\R^d
\,.
\end{equation}
With this, we compute
\begin{equation}
\partial_t V_t
  = -\frac{p(T-t)^{p-1}}{\rho_{T-t}^\beta} - \beta \frac{(T-t)^p}{\rho_{T-t}^{\beta+1}} \partial_t \rho_{T-t}
  \overset{\eqref{e:heat}}{=} -\paren[\bigg]{\frac{p}{T-t} - \frac{\beta \lap \rho_{T-t}}{\rho_{T-t}} } V_t
  \,,
\end{equation}
and  
\begin{equation}
\grad V_t = -\beta \frac{(T-t)^p}{\rho_{T-t}^{\beta+1}} \grad \rho_{T-t} = -\beta V_t s_{T-t},
\end{equation}
so that 
 \begin{equation}
\begin{aligned}
\lap V_t &= \dv (\grad V_t) = -\beta \grad V_t \cdot s_{T-t} - 
\beta V_t \nabla\cdot\Big(\frac{\nabla \rho_{T-t}}{\rho_{T-t}(x)}\Big) \\
&= \beta^2 V_t \abs{s_{T-t}}^2 - \beta V_t \paren[\bigg]{ \frac{\lap \rho_{T-t}}{\rho_{T-t}} - \abs{s_{T-t}}^2 }
= \paren[\Big]{ \beta(1 + \beta) \abs{s_{T-t}}^2 - \beta \frac{\lap \rho_{T-t}}{\rho_{T-t}} } V_t
  \,.
\end{aligned}
\end{equation}
This gives, using the definition of $p$ in~\eqref{26sep902}, 
 \begin{equation}
\begin{aligned}
(\partial_t + \mathcal{L}_t) V_t
  &=\Big(-\frac{p}{T-t}+\frac{\beta\Delta\rho_{T-t}}{\rho_{T-t}}
  -\beta\gamma|s_{T-t}|^2
\\
& 
  +(\gamma-1)\beta(\beta+1)|s_{T-t}|^2
  -(\gamma-1)
\beta \frac{\lap \rho_{T-t}}{\rho_{T-t}}\Big)V_t
 \\
 & = -\biggl( \frac{\beta(\gamma-2)d}{2(T-t)}
    + \beta(\gamma - 2) \frac{\lap \rho_{T-t}}{\rho_{T-t}}
 \label{e:dtPlusL1}
    + \beta \paren[\big]{ \gamma - (\gamma - 1)(1 + \beta) } \abs{s_{T-t}}^2
  \biggr) V_t
  \,.
\end{aligned}
\end{equation}

We now claim that 
\begin{equation}\label{e:LiYau1}
  \frac{\lap \rho_{T-t}}{\rho_{T-t}}
    \geq \abs{\grad \ln \rho_{T-t}}^2 - \frac{d}{2(T-t)}
  \,.
\end{equation}
Inequality~\eqref{e:LiYau1} is the celebrated Li-Yau gradient estimate (see~\cite[Theorem 1.2]{LiYau86}), and is fundamental to geometric analysis.
In our setting, the proof of~\eqref{e:LiYau1} is short, and so we present it here.
Using~\eqref{e:rhoTdef}, we compute
\begin{equation}
  \partial_t \rho_t(x)
    = \int_{\R^d} \paren[\Big]{
	-\frac{d}{2t} + \frac{\abs{x - y}^2}{4 t^2}
      }
      \phi_{2t}(x - y) \, d \rho_0(y)
    \,,
\end{equation}
which gives 
\begin{equation}\label{e:LiYauTmp1}
\paren[\Big]{\partial_t + \frac{d}{2t}} \rho_t(x)
  = \int_{\R^d} \paren[\Big]{ \frac{\abs{x - y}^2}{4 t^2} }
      \phi_{2t}(x - y) \, d \rho_0(y)
      \,.
\end{equation}
Also, we have, using the Cauchy-Schwartz inequality  
 \begin{equation}
\begin{aligned}
  \abs{\grad \rho_t(x)}^2
    &= \abs[\bigg]{ \int_{\R^d} \frac{(x - y)}{2t} \phi_{2t}(x - y) \, \rho_0(dy) }^2
    \\
    &= [\rho_t(x)]^2 \abs[\bigg]{ \int_{\R^d} \frac{(x - y)}{2t} \phi_{2t}(x - y) \, 
    \frac{\rho_0(dy)}{\rho_t(x)} }^2
  \\
    &\leq [\rho_t(x)]^2 \int_{\R^d} \frac{\abs{x - y}^2}{4t^2} \phi_{2t}(x - y) \, \frac{\rho_0(dy)}{\rho_t(x)}
    \overset{\eqref{e:LiYauTmp1}}{=} \rho_t(x) \paren[\Big]{\partial_t + \frac{d}{2t}} \rho_t(x)
    \,.
\end{aligned}
\end{equation}
This, together with~\eqref{e:heat}, immediately implies
\begin{equation}
  \abs{\grad \ln \rho_t(x)}^2
    \leq \frac{\partial_t \rho_t(x)}{\rho_t(x)} + \frac{d}{2t}=
    \frac{\Delta\rho_t(x)}{\rho_t(x)} + \frac{d}{2t}
\end{equation}
Replacing~$t$ with~$T - t$ yields~\eqref{e:LiYau1}, as claimed.
 
Now using the fact that~$\gamma \geq 2$ and~\eqref{e:LiYau1} in~\eqref{e:dtPlusL1} implies
 \begin{equation}
\begin{aligned}
(\partial_t + \mathcal{L}_t) V_t &\le-\Big(
    \frac{\beta(\gamma-2)d}{2(T-t)}+ \beta(\gamma - 2) |s_{T-t}|^2-\frac{d\beta(\gamma-2)}{2(T-t)}
    \\
    &+ \beta \paren[\big]{ \gamma - (\gamma - 1)(1 + \beta) } \abs{s_{T-t}}^2
\Big) V_t
\\
&\le -\beta \paren[\big]{ \gamma - (\gamma - 1)(1+\beta) + (\gamma - 2) } \abs{s_{T-t}}^2 V_t
\\
&=-(\gamma-1)\beta(1-\beta)\abs{s_{T-t}}^2 V_t,
\end{aligned}
\end{equation}
which is~\eqref{e:dtPlusLBound}.

To prove~\eqref{e:dtPlusLhatBound}, we note that
 \begin{equation}
\begin{aligned}
  (\partial_t + \hat{\mathcal{L}}_t) V_t
    &= (\partial_t + \mathcal{L}_t) V_t + (\hat u_t - u_t) \cdot \grad V_t
  \\
  \label{e:dtPlusLHat2}
       &\le -\beta \paren[\big]{
      (1 - \beta) (\gamma - 1) \abs{s_{T-t}}^2
      + \gamma (\hat s_{T-t} - s_{T-t}) \cdot s_{T-t}
  }V_t 
    \,.
\end{aligned}
\end{equation}
By the relative error condition~\eqref{e:RelErrorScore} we know
\begin{equation}
-(\hat s_{T-t} - s_{T-t}) \cdot s_{T-t}
  \le \delta \abs{s_{T-t}}^2
  \,,
\end{equation}
and using this in~\eqref{e:dtPlusLHat2} gives~\eqref{e:dtPlusLhatBound} as desired.
\end{proof}

\subsection{The super-martingale and the proof of Lemma~\ref{l:dXtHatMVanish}}

Once Lemma~\ref{l:genNegative} is established, the proof of Lemma~\ref{l:dXtHatMVanish} is quick.
\begin{proof}[Proof of Lemma~\ref{l:dXtHatMVanish}]
By It\^o's formula and the definition~\eqref{e:Vdef} of $V_t(x)$, we see
  \begin{equation}\label{e:ItodHatZ}
    d \hat Z_t 
      = d(V_t(\hat X_t))
      = \paren[\big]{\partial_t + \hat{\mathcal L}_t} V_t(\hat X_t) \, dt
	- \sigma \beta \hat Z_t s_{T -t}(\hat X_t)
	\, dB_t
    \,.
  \end{equation}
By~\eqref{e:dtPlusLhatBound}, the drift in~\eqref{e:ItodHatZ} 
is negative . Thus,~$\hat Z_t$ is a local super-martingale.
Since~$\hat Z_t \geq 0$, Fatou's lemma implies that~$\hat Z_t$ is a true super-martingale, and Doob's 
martingale convergence theorem implies the existence of an integrable variable~$\hat Z_T$ 
such that~$\hat Z_t \to \hat Z_T$ almost surely as~$t \to T$. Since~$\E \hat Z_T < \infty$, 
$\hat Z_T$ must be finite almost surely.

Finally, since~$\mathcal M$ is compact and~$\rho_0$ is a probability measure 
supported on~$\mathcal M$, the function~$\rho_{T - t}$ converges (as~$t \to T$) to~$0$ locally uniformly on~$\R^d\setminus \mathcal M$. Moreover, Lemma~\ref{l:rhoDecay} gives the estimate
\begin{equation}
\rho_{T-t}(x)\le\frac{C}{(4\pi (T-t))^{(d-m)/2}}e^{-d^2(x,\cM)/[4(T-t)]}.
\end{equation}
It follows that
\begin{equation}\label{26sep904}
\frac{(4\pi)^{\beta (d-m)/2}}{C^\beta}(T-t)^{p+\beta (d-m)/2}e^{\beta d^2(\hat X_t,\cM)/[4(T-t)]}
\le V_t(x)\to \hat Z_T,
\quad\text{ as } t\to T^-,
\end{equation}
and, as we have noted, $\hat Z_T$ is finite almost surely.
This implies~$d(\hat X_t, \mathcal M) \to 0$ almost surely as~$t \to T$.
\end{proof}

\subsection{The limit of the super-martingale}\label{sec:lim-super}

Even though we do not use it in the proof, the limit random 
variable~$\hat Z_T$ above equals~$0$, at least, when~$\beta < 1 - \delta \gamma / (1 - \gamma)$. 
To see this, using~\eqref{e:dtPlusLhatBound} in~\eqref{e:ItodHatZ}, and  
since~$\E \hat Z_T < +\infty$, we deduce that the function~$|s_{T - t}|^2 V_t$ is integrable in time:
\begin{equation}\label{26sep908}
\int_0^T|s_{T - t}(\hat X_t)|^2 V_t(\hat X_t)dt<+\infty \quad\text{a.s.}
\end{equation}
Fix now some~$c_0>0$ and consider the event $\cE_0$ where~$\hat Z_T \geq c_0 > 0$.
Then,~\eqref{26sep908} implies that on $\cE_0$ we have 
\begin{equation}\label{26sep910}
\int_0^T|s_{T - t}(\hat X_t)|^2 dt<+\infty. 
\end{equation}
The lower bound in Lemma~\ref{l:rhoDecay} implies, as in~\eqref{26sep904},
that on the event $\cE_0$ we have 
\begin{equation}\label{26sep906}
C^\beta(4\pi)^{\beta (d-m)/2}(T-t)^{p+\beta (d-m)/2}e^{\beta d^2(\hat X_t,\cM)/[4(T-t)]}
\ge V_t(x)\to \hat Z_T\ge c_0,
\quad\text{ as } t\to T^-.
\end{equation}
Then, on this event we must have
\begin{equation}\label{26sep912}
  d^2(\hat X_t, \mathcal M) \geq C_1 (T - t) \abs{\ln (T - t)} - C_2(T - t),
\end{equation}
with the constants $C_{1,2}$ that depend on $c_0$. 
Now, \eqref{e:uZeta} and~\eqref{e:zetaLipschitz} and~\eqref{26sep912} imply
\begin{equation}
|s_{T-t}(\hat X_t)|=\frac{|\hat X_t-\zeta_{t}(\hat X_t)|}{2(T-t)}
\ge \frac{d(\hat X_t,\cM)-C\sqrt{T-t}}{2(T-t)}
\ge \frac{C\sqrt{|\ln(T-t)|}}{\sqrt{T-t}}.
\end{equation}
Going back to~\eqref{26sep910}, we deduce that 
\begin{equation}
  \int_\tau^T \frac{\abs{\ln(T -t)}}{T - t} \, dt <+ \infty
  \,,
\end{equation}
which is clearly false.
This is a contradiction and hence~$\hat Z_T = 0$ almost surely.

\section{Decay of the density away from \texorpdfstring{$\mathcal M$}{M}: the
proof of Lemma~\ref{l:rhoDecay}}\label{s:rhoDecay}

We now prove Lemma~\ref{l:rhoDecay} using heat kernel estimates.
For this, we first need a lemma bounding the distance to the manifold in coordinate charts.
\begin{lemma}\label{l:chartDistance}
  There exists~$d_0 > 0$ and two constants~$0 < \lambda < \Lambda < \infty$ such that the following holds.
  For every~$x \in \R^d$ with~$d(x, \mathcal M) < d_0$, there is a unique point~$x^* \in \mathcal M$ which is closest to~$x$.
  Moreover, there exists~$U \subseteq \R^m$ open and a chart~$\varphi \colon U \to \R^d$ such that
  \begin{equation}\label{e:xMinusGammaY}
  |x-x^*|^2+\lambda |y-y^*|^2\le |x-\varphi(y)|^2\le |x-x^*|^2+\Lambda |y-y^*|^2\,,
    \quad\text{ for all }
    y\in U
    \,.
  \end{equation}
  and~$y^*$ is the unique point in~$U$ such that~$\varphi(y^*) = x^*$.
\end{lemma}

Momentarily postponing the proof of Lemma~\ref{l:chartDistance}, we now prove Lemma~\ref{l:rhoDecay}.

\begin{proof}[Proof of Lemma~\ref{l:rhoDecay}]
As $x$ lies in a bounded set, it suffices to show that~\eqref{e:RhoTUpperLower} holds for all~$0<t\le t_0$, with some fixed $t_0$. 
Since the density~$\bar\rho_0$ is bounded above and below, it suffices to prove~\eqref{e:RhoTUpperLower} assuming~$\bar \rho_0 = 1$.
In this case, let~$R > 0$ be sufficiently small and write
\begin{equation}\label{26jul2016}
\rho_t(x)
  =\frac{1}{(4\pi t)^{d/2}}\int_\cM e^{-(x-z)^2/(4t)}d\mathcal H^m(z)
  =\frac{1}{(4\pi t)^{d/2}}I_R+\frac{1}{(4\pi t)^{d/2}} \II_R
  \,,
\end{equation}
where
\begin{equation}\label{26jul2008}
I_R \defeq \int\limits_{\cM\cap B(x^*,R)}e^{-(x-z)^2/(4t)} \, d \mathcal H^m(z) \,,
\quad\text{and}\quad
\II_R \defeq \int\limits_{\cM\cap B(x^*,R)^c}e^{-(x-z)^2/(4t)} \, d\mathcal H^m(z),
\end{equation}
and~$x^*$ is unique point in~$\mathcal M$ that is closest to $x$.
Choosing $R$ sufficiently small, we may assume that there is one chart that covers
$\cM\cap B(x^*,R)$.
Let $U\subset\R^{m}$ be an open set such that there exists a~$C^{2,1}$ embedding~$\varphi:U\to \cM\cap B(x^*,R)$ for which
\begin{equation}\label{e:IRGamma}
I_R=\int_U e^{-(x-\varphi(y))^2/(4t)} \sqrt{\det(G)} \, dy
\,.
\end{equation}
Here
\begin{equation}
  G \defeq (D \varphi)^T D\varphi\,,
\end{equation}
is the Gram matrix.
Making $R$ smaller if necessary, let~$0 < \lambda < \Lambda$ be the constants provided by Lemma~\ref{l:chartDistance} such that~\eqref{e:xMinusGammaY} holds.
Since~$G$ is bounded above, replacing~$U$ with~$\R^m$ in~\eqref{e:IRGamma} and using~\eqref{e:xMinusGammaY} immediately implies the upper bound
\begin{equation}\label{e:IRUpper}
    I_R
    \le c_2 t^{m/2} e^{-|x-x^*|^2/(4t)}
    \,,
\end{equation}
for some finite constant~$c_2$.
Since~$\det(G)$ is bounded below away from~$0$, using~\eqref{e:xMinusGammaY} in~\eqref{e:IRGamma} and explicitly integrating yields the lower bound
\begin{equation}\label{e:IRLower}
    I_R
    \geq c_1 t^{m/2} e^{-|x-x^*|^2/(4t)}
    \,,
\end{equation}
for some constant~$c_1 > 0$, provided~$t \leq R^2/C$.

To finish the proof, we will show that the term~$\II_R$ in~\eqref{26jul2016} is exponentially small.
We do this using Lemma~\ref{l:chartDistance}.
By covering~$\mathcal M$ by finitely many charts for which~\eqref{e:xMinusGammaY} holds, we see that for all sufficiently small~$R > 0$, we have
\begin{equation}\label{26jul2010}
  |x-z|^2  -|x-x^*|^2\geq \lambda R^2
  \,,\quad
  \text{provided }
  z \in \mathcal M\setminus B(x^*, R)
  \,.
\end{equation}
This implies
\begin{equation}
  \II_R\le Ce^{- \lambda R^2/(4t)}e^{-|x-x^*|^2/(4t)}.
\end{equation}
Thus, if we choose $R \geq N\sqrt{t}|\log t|$ with a sufficiently large $N$, we can ensure
\begin{equation}\label{26jul2014}
\II_R\le\frac{I_R}{10}.
\end{equation}
Now~\eqref{e:RhoTUpperLower} follows from~\eqref{26jul2016},~\eqref{e:IRUpper}, \eqref{e:IRLower} and~\eqref{26jul2014}, concluding the proof.
\end{proof}

It remains to prove Lemma~\ref{l:chartDistance}.
\begin{proof}[Proof of Lemma~\ref{l:chartDistance}]
  First, the existence of~$d_0 > 0$ so that every point~$x \in B(\mathcal M, R)$ has a unique projection is an immediate consequence of the fact that~$\mathcal M$ is~$C^2$ and tubular neighborhood theorem (see for instance~\cite[Chapter 4.5]{Hirsch94}).
  To verify~\eqref{e:xMinusGammaY} uniformly in the tubular neighborhood, suppose~$U \subseteq \R^m$ is open and~$\varphi \colon U \to \mathcal M$ is a chart.
  Let~$S_x \colon U \to \R$ be defined by~$S_x(y) = \abs{x - \varphi(y)}^2$.
  We compute
  \begin{equation}
    \hess(S_x)(y) = 2 G(y) - 2 R(x, y)\,,
  \end{equation}
  where the \emph{Gram matrix} $G$ and \emph{normal component curvature matrix} $R$ are defined by
  \begin{equation}
    G_{i,j}(y) \defeq \partial_i \varphi(y) \cdot \partial_j \varphi(y)
    \quad\text{and}\quad
    R_{i,j}(x, y) \defeq (x - \varphi(y)) \cdot \partial_i \partial_j \varphi(y)
    \,.
  \end{equation}
  For any non-zero~$v \in \R^m$, linear independence of~$\set{\partial_1 \varphi(y), \dots, \partial_m \varphi(y)}$ implies 
  \begin{equation}
    v \cdot G(y) v = \abs[\Big]{\sum_{i = 1}^m v_i \partial_i \varphi (y)}^2 > 0
    \,.
  \end{equation}
  Replacing~$U$ with a compactly contained neighborhood if necessary, we can now find~$0 < \lambda < \Lambda$ such that
  \begin{equation}
    \lambda \abs{v}^2 \leq v \cdot G(y) v \leq \Lambda \abs{v}^2
    \quad\text{for all } y \in U\,, v \in \R^m
    \,.
  \end{equation}
  Since~$\varphi$ is~$C^{2}$ we know~$\abs{v \cdot R(x, y) v } \leq C \abs{x - \varphi(y)}$ and so
  \begin{equation}
    2(\lambda - C\abs{\varphi(x) - y})) \abs{v}^2
      \leq v \cdot \hess(S_x)(y) v
      \leq 2(\Lambda + C\abs{\varphi(x) - y})) \abs{v}^2
  \end{equation}

  Now restricting to the neighborhood where~$\lambda > 2 C\max_y \abs{\varphi(x) - y}$ gives a family of charts for which we have
  \begin{equation}
    \lambda \abs{v}^2 \leq v \cdot \hess(S_x)(y)v \leq \Lambda \abs{v}^2
  \end{equation}
  for every~$x \in B(\mathcal M, R)$, and every~$y$ close to~$y^*$.
  Taylor expanding~$S_x(y)$ near~$y = y^*$ proves~\eqref{e:xMinusGammaY} in a small neighborhood around~$y^*$.
  Using compactness of~$M$ and covering~$M$ by finitely many such charts finishes the proof.
\end{proof}

\section{Convergence of \texorpdfstring{$\zeta_t$}{zeta t} as \texorpdfstring{$t \to T$}{t to T}: the proof of Lemma~\ref{l:zetaBdKTau-bis2}}\label{s:zetaConv}

We now prove Lemma~\ref{l:zetaBdKTau-bis2} by using the Laplace method.

\begin{proof}[Proof of Lemma~\ref{l:zetaBdKTau-bis2}]
Using~\eqref{e:rhoTdef} and~\eqref{e:zetaDef}, we compute
\begin{equation}
\zeta_{T-\tau}(x)- x^* =\frac{\displaystyle
	\int_{\R^d} (z-x^*) \exp\paren[\Big]{ -\frac{ |x - z|^2 }{4 \tau} } \, d\rho_0(z)  
    }
    {\displaystyle
	\int_{\R^d}\exp\paren[\Big]{ -\frac{ |x - z|^2 }{4 \tau}} \, d\rho_0(z)  
    }\,.
\end{equation}
Take $R_1 > 0$ to be determined later %
and write 
\begin{equation}\label{25oct802bis}
|\zeta_{T-\tau}(x)-x^*|\le I_{R_1} + \II_{R_1}
\,,
\end{equation}
where
\begin{equation}\label{e:IR1def}
I_{R_1} \defeq \frac{\displaystyle
	\abs[\bigg]{ \int_{B(x^*,R_1)} (z-x^*) \exp \paren[\Big]{ -\frac{ |x - z|^2 }{4 \tau} } \, d\rho_0(z)   }
    }
    {\displaystyle
	\int_{\R^d}\exp \paren[\Big]{ -\frac{ |x - z|^2 }{4 \tau}} \, d\rho_0(z) 
    }\, ,
\end{equation}
and
\begin{equation}
\II_{R_1} = \frac{\displaystyle
	\int_{B^c(x^*,R_1)} |z-x^*|\exp \paren[\Big]{ -\frac{ |x - z|^2 }{4 \tau} } \, d\rho_0(z)  
    }
    {\displaystyle
	\int_{\R^d}\exp \paren{ -\frac{ |x - z|^2 }{4 \tau}} \, d\rho_0(z)  
    }\,.
\end{equation}

We will ensure that~$R_1$ is sufficiently small so that~$B(x^*, R_1) \cap \cM$ is covered by a single coordinate chart~$\varphi \colon U \to \cM$ for some domain~$U \subseteq \R^m$.
For notational simplicity, and without loss of generality, we assume~$\varphi(0) = x^*$.
We can write
\begin{equation}
\label{e:IR1Coordinates}
I_{R_1} \leq \frac{\displaystyle
	\abs[\bigg]{ \int_{U \cap \varphi^{-1}(B(x^*,R_1))} F(y) \exp \paren[\Big]{ -\frac{ S_x(y) }{4 \tau} } \, d\rho_0(y)   }
    }
    {\displaystyle
	\int_{U \cap \varphi^{-1}(B(x^*,R_1))} \exp \paren[\Big]{ -\frac{ S_x(y) }{4 \tau} } \, d\rho_0(y)
    }\, ,
\end{equation}
where
\begin{equation}
  \label{e:Sdef}
  S_x(y) \defeq \abs{x - \varphi(y)}^2
    \,,
\end{equation}
and
\begin{equation}
  \label{e:Fdef}
  F(y) = (\varphi(y) - x^*) \rho_0 \circ \varphi(y) \sqrt{\det (D \varphi)^T D\varphi} \,.
\end{equation}

Since~$\mathcal M$ is closed and~$C^{2, 1}$, there exists~$d_0 > 0$ such that
\begin{equation}
  \mathcal T(\mathcal M) \defeq \set{ x \in \R^d \st d(x, \mathcal M) < d_0 }
    \,,
\end{equation}
is contained in a tubular neighborhood (see for instance~\cite[Chapter 4.5]{Hirsch94}).
All points in~$\mathcal T(\mathcal M)$ have unique projections onto~$\mathcal M$, and,
as in Lemma~\ref{l:chartDistance}, the Hessian of~$S_x$ is uniformly bounded (above and below) in this neighborhood.
In this case, one can estimate~\eqref{e:IR1Coordinates} by approximating the numerator and denominator by Gaussian integrals.
If~$x \not\in \mathcal T(\mathcal M)$, it may be close to the cut-locus or points where the Hessian of~$S_x$ is degenerate.
In this case, we just estimate~$I_{R_1}$ by the maximum of~$\abs{F}$, and choose~$R_1$ appropriately.
We now analyze each of these cases separately.

\textbf{Case I: $x \not \in \mathcal T(\mathcal M)$.}
In this case, we use~$\abs{z - x^*} < R_1$ in~\eqref{e:IR1def} to obtain
\begin{equation}\label{e:IR1case1}
  I_{R_1} \leq R_1
  \,.
\end{equation}
For the second term~$\II_{R_1}$, we use the fact that~$\cM$ is bounded to obtain
\begin{equation}\label{e:IIR1Case1}
 \II_{R_1}\le 2D_{\mathcal M} \frac{\displaystyle
	\int_{B(x^*,R_1)^c}\exp \paren[\Big]{ -\frac{ |x - z|^2 }{4 \tau} }
	\, d\rho_0(z) 
    }
    {\displaystyle
	\int_{\R^d}\exp \paren[\Big]{ -\frac{ |x - z|^2 }{4 \tau} } \,  d\rho_0(z)  
    }
    \,.
\end{equation}
where~$D_{\mathcal M}$ is defined by
\begin{equation}\label{e:DMdef}
  D_{\mathcal M} \defeq \sup_{z\in\cM}|z|
  \,.
\end{equation}

To bound the denominator in \eqref{e:IIR1Case1} from below, we apply Lemma~\ref{l:chartDistance} and find~$R_1, \Lambda > 0$ (independent of~$x$) such that~\eqref{e:xMinusGammaY} holds.
Together with the lower bound in \eqref{e:rho0Bounded}, this gives the estimate
\begin{equation}\label{25oct818} %
\begin{aligned}
\MoveEqLeft
\int_{\R^d}\exp\paren[\Big]{ -\frac{|x - z|^2 }{4 \tau}} \, d\rho_0(z)
  \ge \rho_{0, \inf}\int_{B(x^*,R_1)\cap\cM} \exp\paren[\Big]{ -\frac{|x - z|^2 }{4 \tau} }  \, d\mathcal H^m(z)
\\
&\ge \rho_{0, \inf}\int_{B(x^*,R_1)\cap\cM} \exp \paren[\Big]{ -\frac{ |x -x^*|^2 + \Lambda |z-x^*|^2}{4 \tau} } \, d\mathcal H^m(z)\\
&\ge \rho_{0, \inf}\tau^{m/2}\exp \paren[\Big]{ -\frac{|x -x^*|^2 }{4 \tau} }\int_{B(0,R_1/\sqrt{\tau})} 
e^{- \Lambda |y|^2/4} \, dy\\
&=\rho_{0, \inf }\tau^{m/2}\exp \paren[\Big]{ -\frac{|x -x^*|^2 }{4 \tau} }
\varPhi\Big(\frac{R_1}{\sqrt{\tau}}\Big) \,.
\end{aligned}
\end{equation} 
Here,~$\varPhi$ is defined by
\begin{equation}\label{26jul2022}
\varPhi(s) \defeq \int_{B(0,s) }e^{- \Lambda |y|^2/4} \, dy
\,.
\end{equation}
We clarify that the integrals in~\eqref{25oct818} and \eqref{26jul2022} are over $m$-dimensional balls.  

Next, to bound the numerator in~\eqref{e:IIR1Case1} from above, we need a lower bound on~$|x - z|$.
We claim that if~$d(x, \cut(\mathcal M)) > r$ (which is certainly true in~$\mathcal G(r, R)$), then we have the lower bound
\begin{equation}\label{e:xMinusZLower}
  \abs{x - z}^2 \geq \abs{x - x^*}^2 + \abs{z - x^*}^2 \paren[\Big]{ \frac{r}{r + \abs{x - x^*}} }
  \,.
\end{equation}
To avoid breaking continuity, we briefly postpone the proof of~\eqref{e:xMinusZLower} to Lemma~\ref{l:distanceGap}, below.

Using~\eqref{e:xMinusZLower} and the fact that~$|z - x^*| \geq R_1$ on~$B(x^*, R_1)^c$, we see
\begin{equation}
  \abs{x - z}^2 \geq \abs{x - x^*}^2 + R_1^2 \paren[\Big]{ \frac{r}{r + \abs{x - x^*}} }
\end{equation}
provided~$x \in \cG(r, R)$. 
Moreover, since~$x \in \mathcal G(r, R) \subseteq B(0, R)$, we know~$\abs{x - x^*} \leq R + D_{\mathcal M}$.
Thus, there exists a constant~$C > 0$ (depending only on~$R$ and~$\mathcal M$) such that
\begin{equation}
  \abs{x - z}^2 \geq \abs{x - x^*}^2 + C_R r R_1^2
  \,,
\end{equation}
for all~$r \leq 1$.
Since~$\rho_0$ is a probability measure, this implies
\begin{align}
\MoveEqLeft
\int_{B^c(x^*,R_1)}\exp \paren[\Big]{ -\frac{ |x - z|^2 }{4 \tau} } \, d\rho_0(z) 
\le \int_{B^c(x^*,R_1)}\exp\paren[\Big]{ -\frac{ |x -x^*|^2 + C_R r R_1^2}{4 \tau} } \, d\rho_0(z)  \\
\label{25oct820}
&\le \exp\paren[\Big]{ -\frac{ |x -x^*|^2}{4\tau}}
  \exp\paren[\Big]{-\frac{C_R r R_1^2}{4\tau}} 
  \,.
\end{align}
Using~\eqref{25oct818} and~\eqref{25oct820} in~\eqref{e:IIR1Case1} implies
\begin{equation}\label{e:IIR1Bound}
\II_{R_1} \le \frac{1}{2 D_{\mathcal M} \rho_{0, \inf}\tau^{m/2} } \exp\paren[\Big]{ -\frac{ C_R r R_1^2}{4 \tau}}
\,,
\end{equation}
provided~$0<\tau<R_1^2$.
Combining~\eqref{e:IR1case1} and~\eqref{e:IIR1Bound} gives
\begin{equation}\label{e:ZetatMinusTCase1}
|\zeta_{T-\tau}(x)-x^*|\le 
R_1+\frac{C}{\tau^{m/2}} \exp\paren[\Big]{-\frac{C_R r R_1^2}{4 \tau}}
\,,
\end{equation}
provided we again have~$0<\tau<R_1^2$.

Now choose
\begin{equation}
r=K \tau^c,
\quad
R_1=\tau^b d_0,     %
\end{equation}
with $c+2b<1$. 
Since in this case we assumed~$x \notin \mathcal T(\mathcal M)$, we know~$\abs{x - x^*} \ge d_0$.
Thus if
\begin{equation}
  \tau^{1 - 2b} < d_0^2
  \quad\text{and}\quad
  \tau^b < \frac{1}{4}
  \,,
\end{equation}
then~$0 < \tau < R_1^2$ and~$R_1 < d_0/4 \le |x - x^*|/4$, so~\eqref{e:ZetatMinusTCase1} holds.
This implies
\begin{equation}
\begin{aligned}
|\zeta_{T-\tau}(x)-x^*|&\le 
  C\tau^b
  + \frac{C}{\tau^{m/2}} 
      \exp \paren[\Big]{-\frac{ C_R K d_0^2}{4\tau^{1-c-2b}}}
\\
\label{e:zetaTmp2}
&\le C\tau^b 
\,,
\end{aligned}
\end{equation}
provided $\tau$ is sufficiently small, and $x\in\cG(K \tau^c, R)$.
\medskip

\textbf{Case II: $x \in \mathcal T(\mathcal M)$}
In this case we will choose~$R_1 > 0$, independent of~$x, \tau$, to be small enough so that so that~$B(x^*, R_1) \cap \cM$ is always covered by a single normal coordinate chart~$\varphi \colon U \to \cM$ with~$\varphi(0) = x^*$.
Here~$U \subseteq \R^m$ is an open domain.
Since~$R_1$ is fixed, we will have to estimate~$I_{R_1}$ differently and obtain an estimate that is better than~\eqref{e:IR1case1}.

In order to do this, we first use
Lemma~\ref{l:chartDistance} to obtain
\[
  |x-\varphi(y)|^2
  \geq |x-x^*|^2+\lambda|y|^2
\]
in the coordinate chart, with $\lambda>0$ independent of $x$.
Moreover, $|\varphi(y)-x^*|\leq C|y|$, and the density
and the coordinate Jacobian are bounded uniformly.

Then, the numerator of~$I_{R_1}$ satisfies
\begin{equation}
\begin{aligned}
  \operatorname{Num}(I_{R_1})
  &\defeq
  \abs[\bigg]{
    \int_{B(x^*,R_1)}(z-x^*)
    \exp\paren[\Big]{-\frac{|x-z|^2}{4\tau}}
    \,d\rho_0(z)
  }\\
  &\leq
  C\exp\paren[\Big]{-\frac{|x-x^*|^2}{4\tau}}
  \int_{\mathbb R^m}|y|
    \exp\paren[\Big]{-\frac{\lambda|y|^2}{4\tau}}\,dy\\
  &\leq
  C\tau^{(m+1)/2}
  \exp\paren[\Big]{-\frac{|x-x^*|^2}{4\tau}}.
\end{aligned}
\end{equation}
Using the same lower bound~\eqref{25oct818} for the denominator
as in the previous case, we see
\begin{equation}\label{e:IR1FinalCase2}
  I_{R_1}\leq C\sqrt{\tau}.
\end{equation}

For the term~$\II_{R_1}$ we note there exists~$r_0 > 0$ such that~$\mathcal T(\mathcal M) \subseteq \mathcal G(r_0, R)$.
We now bound~$\II_{R_1}$ with the same bound~\eqref{e:IIR1Bound}, but use~$r_0$ instead of~$r$.
Combining this with~\eqref{e:IR1FinalCase2} gives
\begin{equation}\label{e:ZetatMinusTCase2}
|\zeta_{T-\tau}(x)-x^*|\le
  C \paren[\bigg]{ \sqrt{\tau}
  +\frac{1}{\rho_{0,\inf}\tau^{m/2}}
  \exp\paren[\Big]{-\frac{r_0 R_1^2}{8 \tau}} }
 \leq C\sqrt{\tau}
 \,.
\end{equation}
\medskip

Thus using~\eqref{e:ZetatMinusTCase1} outside~$\mathcal T(\mathcal M)$, or~\eqref{e:ZetatMinusTCase2} inside~$\mathcal T(\mathcal M)$ we see, since $b<1/2$, that~\eqref{e:zetaAsTtoBigT} holds in either case.
This concludes the proof of Lemma~\ref{l:zetaBdKTau-bis2}.
\end{proof}

It remains to prove the distance gap estimate~\eqref{e:xMinusZLower}. 
We state and prove this as our next lemma.

\begin{lemma}\label{l:distanceGap}
Suppose~$x \in \R^d$ is such that~$d( x, \cut(\mathcal M)) > r$, and~$x^* \in \mathcal M$ is the projection of~$x$.
Then for any~$z \in \mathcal M$ we have~\eqref{e:xMinusZLower}.
\end{lemma}
\begin{proof}
Since~$d(x, \cut(\mathcal M)) > r$, we can extend the ray from~$x^*$ through~$x$ by a distance~$r$ without crossing the cut-locus.
Let
\begin{equation}\label{e:xPrimeDef}
x' \defeq x + r \frac{x - x^*}{\abs{x - x^*}}
\end{equation}
be the point on this ray such that~$\abs{x' - x} = r$.
Clearly~$|x' - x^*| = r_* + r$, where~$r_* \defeq |x - x^*|$.
Also, $x^*$ is the point on~$\mathcal M$ that is closest to~$x'$, and so
\begin{equation}
  \abs{z - x'}^2 \geq \abs{x^* - x'}^2 = (r_* + r)^2
  \,.
\end{equation}
Expanding the left side around~$x^*$ yields
\begin{equation}
\begin{aligned}
  (r_* + r)^2 &\leq \abs{z - x'}^2
   = \abs[\big]{(z - x^*) - (x' - x^*)}^2
  \\
    &= \abs{z - x^*}^2 + (r_* + r)^2 - 2 (z - x^*) \cdot( x' - x^*)
  \,,
\end{aligned}
\end{equation}
and hence
\begin{equation}\label{e:innerProdBound1}
  2 (z - x^*) \cdot( x' - x^*) \leq \abs{z - x^*}^2
  \,.
\end{equation}
Using~\eqref{e:xPrimeDef}, this implies
\begin{equation}\label{e:innerProdBound2}
  2 (z - x^*) \cdot( x - x^*) \leq \frac{r_*}{r_* + r} \abs{z - x^*}^2
  \,.
\end{equation}
Finally, expanding~$\abs{x - z}^2$ around~$x^*$ gives
\begin{equation}
\begin{aligned}
  \abs{x - z}^2
  &= \abs{z - x^*}^2 + \abs{x - x^*}^2 - 2 (z - x^*) \cdot (x - x^*)
  \\
  &\overset{\mathclap{\eqref{e:innerProdBound2}}}{\geq}~ \abs{x - x^*}^2 + \abs{z - x^*}^2 \paren[\Big]{ 1 - \frac{r_*}{r_* + r} }
  \,,
\end{aligned}
\end{equation}
which implies~\eqref{e:xMinusZLower} as desired.
\end{proof}

\section{Gaussian measure of the fattened cut-locus}\label{s:ballVolFractions2}

\subsection{Proof of Lemma~\ref{l:ballVolumeFractions2}}

In this section, we prove Lemma~\ref{l:ballVolumeFractions2} except for the auxiliary Lemma \ref{l:cutTubeBd}. The proof 
of that lemma is, in turn, postponed until Section~\ref{s:cutTubeBd}.
We first show that the desired Gaussian measure bound in Lemma~\ref{l:ballVolumeFractions2} is equivalent to a tube volume estimate.
\begin{lemma}[Tube volume equivalence]\label{l:tubeEquivalence}
  Let~$E \subseteq \R^d$ be a bounded measurable set, and fix~$\epsilon \in [0, d-1]$.
  Then
  \begin{equation}\label{e:gaussianBound}
    \sup_{\substack{
      x \in \R^d\,,\\
      r > 0,~
      \kappa \in (0, 1]
    }}
      \frac{\kappa^{\frac{1 + \epsilon}{2}}}{r}
      \P\paren[\big]{
	\mathcal N(x, \kappa I_d) \in B(E, r)
      }
    < \infty
  \end{equation}
  if and only if there exists~$r_0 > 0$ such that
  \begin{equation}\label{e:tubeBound}
    \sup_{\substack{
      x \in \R^d\,,
      \\
      r \in (0, r_0),~
      s > 0
    }} \frac{ \abs{B(E, r) \cap B(x, s)}}{r s^{d - 1 - \epsilon}}
    < \infty
    \,.
  \end{equation}
\end{lemma}

Here,~$B(E, r)$ denotes an~$r$-neighborhood of~$E$.
If~$E$ is a compact~$d-1$-dimensional~$C^1$-sub-manifold, or is bounded and contained in a graph of a Lipschitz function of~$d - 1$ variables, then~\eqref{e:tubeBound} is certainly satisfied with~$\epsilon = 0$.
On the other hand, if~$E$ is a non-compact~$d-1$ dimensional~$C^1$-manifold, then it could cluster 
at a point, making~\eqref{e:tubeBound} fail for~$\epsilon = 0$.
This happens, for instance, with the spiral~$\set{t^{-2} (\cos t, \sin t) \st t \in (1, \infty)}$.
Nevertheless, such sets may still satisfy~\eqref{e:tubeBound} for some~$\epsilon>0$, and this weaker estimate is sufficient for our purposes.

As we have mentioned in the introduction, the regularity of the cut-locus has been studied since the~1970's.
While the cut-locus can be a wild fractal set~\cite{GluckSinger78,GluckSinger79}, it is typically much more regular.
Indeed, classical results~\cite{
  Wall77,
  Buchner78,
  Yomdin81,
  Mather83} 
show that the cut-locus of generic smooth manifolds is stratifiable, with top dimension at most~$d - 1$.
More recent results in~\cite{ItohTanaka01,LiNirenberg05} show that $\cut_R(\mathcal M)$ is contained 
in the Lipschitz image of the unit normal bundle -- a compact~$d-1$ dimensional manifold.
However, this does not directly imply~$\cut_R(\mathcal M)$ satisfies~\eqref{e:tubeBound}, even with~$\epsilon > 0$.
We will instead use a polynomial approximation argument combined with quantitative semialgebraic estimates, in the spirit of Yomdin's~\cite{Yomdin15} methods, 
with the ideas also coming from~\cite{BasuPollackEA06}, to prove that~$\cut_R(\mathcal M)$ 
satisfies~\eqref{e:tubeBound} for all small~$\epsilon > 0$.

We begin by defining the \emph{truncated normal map} which allows us to contain~$\cut_R(\mathcal M)$ in a Lipschitz image of the unit normal bundle.
Let
\begin{equation}\label{26oct514}
  \mathcal U=\{(p,\nu):p\in\mathcal M,\ \nu\perp T_p\mathcal M,
      \ |\nu|=1\}
\end{equation}
be the unit normal bundle to $\cM$. This is a compact smooth manifold of
dimension~$d-1$.
Define the normal cut time by
\begin{equation}\label{26oct512}
t_*(p,\nu)=\sup\{t\geq 0:p\text{ is a closest point on }\mathcal M
\text{ to }x=p+t\nu\}\,.
\end{equation}
We also recall the definition of the reach of a manifold:
\[
\operatorname{reach}(\mathcal M)=\sup\{r> 0:~\text{every $x\in\R^d$ s.t. $\dist(x,\cM)<r$
has a unique closest point on $\cM$}\}. 
\]
Fix~$0<a<\operatorname{reach}(\mathcal M)$ and
\begin{equation}\label{26oct516}
L>\max\{a,R+\sup_{p\in\mathcal M}|p|\}.
\end{equation}
The \emph{truncated normal map}, denoted by~$F_L \colon \mathcal U \to \R^d$ is defined by
\begin{equation}\label{e:truncatedNormalMap}
  F_L(p,\nu)=p+t_L(p,\nu)\nu\in\R^d
  \,,
  \quad\text{where}\quad
  t_L(p,\nu)=\min\{t_*(p,\nu),L\}
  \,,~~(p,\nu)\in\cU. 
\end{equation}
We now note that the image of the truncated normal map contains~$\cut_R(\mathcal M)$ (see for instance~\cite{ItohTanaka01}, or the short proof in Section~\ref{s:endPointMap}, below).
\begin{lemma}\label{l:endPointMap}
The image~$F_L(\cU)$ contains~$\cut_R(\mathcal M)$.
\end{lemma}

In light of Lemma~\ref{l:endPointMap}, one can use the manifold structure of~$\mathcal U$ and parametrization~$F_L$ to obtain tube neighborhood bounds of~$\cut_R(\mathcal M)$.
If~$F_L$ were bi-Lipschitz, then~$F_L(\mathcal U)$ would satisfy~\eqref{e:tubeBound} with~$\epsilon=0$, immediately implying Lemma~\ref{l:ballVolumeFractions2}.
In general, however, $F_L$ need not be bi-Lipschitz.
Itoh and Tanaka~\cite{ItohTanaka01} (see also~\cite{LiNirenberg05}) show that~$F_L$ is Lipschitz, and using a semi-concavity argument~\cite{CannarsaSinestrari04} one can readily show~$D F_L \in \BV$.
We will now use polynomial approximation and quantitative semialgebraic estimates to establish~\eqref{e:tubeBound}, provided~$s^{d'}\leq r\leq s$ for some~$d' \geq 1$.
This is our next lemma.

\begin{lemma}[Large tube bounds]\label{l:cutTubeBd}
Let $F_L:\mathcal U\to\mathbb R^d$ be the truncated normal map defined in $\eqref{e:truncatedNormalMap}$.
Then $F_L$ is Lipschitz, $DF_L\in BV(\mathcal U)$, and $\cut_R(\mathcal M)\subseteq F_L(\mathcal U)$.
Moreover, for all~$\eps>0$ and~$d' \geq 1$, there are constants $C>0$, $s_0 \in (0, 1]$ such that the following hold.
\begin{enumerate}
\item
  For every $x\in\R^d$ and $s \in [0, s_0]$, the parametrized~$(d-1)$-dimensional area of~$F_L$ over $B(x, s)$ satisfies
  \begin{equation}\label{e:Jdm1FL}
  \int_{F_L^{-1}(B(x,s))} J_{d-1} F_L \,d\mathcal H^{d-1}
  \leq Cs^{d-1-\eps},
  \quad\text{where}\quad
    J_{d-1}F_L=\sqrt{\det(DF_L^T DF_L)}
    \,.
  \end{equation}

\item For every $r, s$ such that~$0 \leq s^{d'}\leq r\leq s \leq s_0$,
\begin{equation}\label{e:tubeVolumeRlarge}
  |B(\cut_R(\mathcal M),r)\cap B(x,s)|
  \leq Crs^{d-1-\eps}.
\end{equation}
\end{enumerate}
\end{lemma}

The restriction~$r \geq s^{d'}$ in~\eqref{e:tubeVolumeRlarge} arises because the proof of Lemma~\ref{l:cutTubeBd} uses a polynomial approximation argument, and
when~$r$ is small the approximation error may exceed the tube width.
For narrow tubes, the classical theorem of Federer on Minkowski content (see for instance~\cite{EvansGariepy92}) gives an asymptotic tube-volume formula for compact Lipschitz images.
For our purposes, we need a localized, quantitative version, that provides an explicit error bound when the derivative of the parametrization is additionally~$\BV$.
This is our next lemma, and the last ingredient needed in the proof of Lemma~\ref{l:ballVolumeFractions2}.

\begin{lemma}[Tube bound for images]\label{l:BVTube}
Let~$Q\Subset Q'$ be cubes in~$\R^{d-1}$.
Suppose that~$f:Q'\to\R^{d}$ is Lipschitz and~$Df\in\BV(Q')$.
There exist constants~$C,r_0>0$, depending only on the cubes, the
Lipschitz constant of $f$, and~$\|Df\|_{BV(Q')}$, 
such that for every~$x\in\R^{d}$, $s>0$, and~$0<r\leq r_0$ we have
\begin{equation}\label{e:BVTube}
  |B(f(Q),r)\cap B(x,s)|
  \leq C r\int_{Q'\cap f^{-1}(B(x,s+Cr))}J_{d-1}f(y)\,dy
       +Cr^{4/3}
    \,.
\end{equation}
Here (as in~\eqref{e:Jdm1FL}) the notation~$J_{d-1} f$ refers to the~$(d-1)$-dimensional Jacobian defined by
\begin{equation}\label{e:tanjentialJacobianDef}
  J_{d-1} f \defeq \sqrt{\det(Df^T Df)}
  \,,
\end{equation}
\end{lemma}

Momentarily postponing the proofs of the above lemmas, we now prove Lemma~\ref{l:ballVolumeFractions2}.

\begin{proof}[Proof of Lemma~\ref{l:ballVolumeFractions2}]
For any~$\eps\in(0,1)$, we will show that~$E=\cut_R(\mathcal M)$ satisfies~\eqref{e:tubeBound}. 
Once this is established, Lemma~\ref{l:tubeEquivalence} immediately implies Lemma~\ref{l:ballVolumeFractions2} for~$\epsilon \in (0, 1)$.
Since~$\kappa \leq 1$ this implies~\eqref{e:PNinCutRM2} also holds for every~$\epsilon \geq 1$, proving Lemma~\ref{l:ballVolumeFractions2}.
\smallskip

To prove that~$E$ satisfies~\eqref{e:tubeBound}, choose~$r_0\leq1$ sufficiently small for Lemma~\ref{l:BVTube} let~$\set{\psi_i}$ be a finite atlas of~$\mathcal U$, chosen so that there exist cubes~$Q_i \Subset Q_i' \Subset \Dom(\psi_i)$ such that~$\mathcal U = \cup_i \psi_i(Q_i)$.
Let~$r_0$ be small enough so that Lemma~\ref{l:BVTube} applies to each of the coordinate cubes~$Q_i \Subset Q_i'$, and fix~$r \in [0, r_0]$, $s > 0$.
Let $s_0>0$ be a small constant that will be chosen later.
We now divide the analysis into cases.
\medskip

\emph{Case 1: $0 \leq s \leq r$.}
In this case we have the simple bound 
\begin{equation}\label{26oct1302}
  |B(E,r)\cap B(x,s)|
  \leq C_ds^d
  \leq C_drs^{d-1-\eps}.
\end{equation}
\medskip

\emph{Case 2: $0 \leq r \leq s_0 < s$.}
By Lemma~\ref{l:cutTubeBd} we know~$\cut_R(\mathcal M) \subseteq F_L(\mathcal U)$, and that~$F_L$ is Lipschitz.
Since~$\mathcal U$ is a $(d-1)$-dimensional manifold this implies
\begin{equation}\label{26oct1304}
|B(E,r)\cap B(x,s)|\le |B(E,r)|\le|B(F_L(\cU),r)|\leq C_d r\le \frac{C_d}{s_0^{d - 1 - \epsilon}} r s^{d-1-\eps}.
\end{equation}
\medskip

\emph{Case 3: $s^{3d} \leq r\leq s\leq s_0$.}
In this case the bound~\eqref{e:tubeVolumeRlarge} in Lemma~\ref{l:cutTubeBd} with~$d' = 3d$ immediately implies~$E$ satisfies~\eqref{e:tubeBound}.
\medskip

\emph{Case 4: $r < s^{3d} \leq s\leq s_0$.}
By Lemma~\ref{l:cutTubeBd} we know~$E \subseteq F_L(\mathcal U)$,  $F_L$ is Lipschitz and~$DF_L$ is~$\BV$.
Applying Lemma~\ref{l:BVTube} to~$F_L$ in every coordinate cube~$Q_i$, and using Lemma~\ref{l:cutTubeBd} implies
\begin{equation}\label{26oct1306}
\begin{aligned}
  |B(E,r)\cap B(x,s)|
  &\overset{\eqref{e:BVTube}}{\leq} Cr\int_{F_L^{-1}(B(x,s+Cr))}
      J_{d-1}F_L\,d\mathcal H^{d-1}
      +Cr^{4/3}\\
  &\overset{\eqref{e:Jdm1FL}}{\leq}
    Cr(s+Cr)^{d-1-\eps}+Cr^{4/3}\leq Crs^{d-1-\eps}.
\end{aligned}
\end{equation}
Here we used the inequalities~$s+Cr\leq C's$ and $r^{1/3}\leq s^d\leq s^{d-1-\eps}$,
which hold because $r<s^{3d}$ and $0<s\leq s_0\leq1$.
\medskip

Thus in all cases, the set~$E=\cut_R(\mathcal M)$ satisfies~\eqref{e:tubeBound}, and now Lemma~\ref{l:tubeEquivalence} implies~\eqref{e:gaussianBound}, concluding the proof.
\end{proof}

It remains to prove Lemmas~\ref{l:tubeEquivalence}--\ref{l:BVTube}.
We prove Lemma~\ref{l:tubeEquivalence} in Section~\ref{s:tubeEquivalence}, and Lemma~\ref{l:BVTube} in Section~\ref{s:BVTube}.
Lemma~\ref{l:endPointMap} is well known, but since the proof is short we present the proof in Section~\ref{s:endPointMap}.
Finally, the proof of Lemma~\ref{l:cutTubeBd} is a little more involved, and is presented in Section~\ref{s:cutTubeBd}.

\subsection{The equivalent tube volume bound: the proof of Lemma~\ref{l:tubeEquivalence}}\label{s:tubeEquivalence}

We now prove Lemma~\ref{l:tubeEquivalence} using Gaussian density bounds and a standard covering argument.

\begin{proof}[Proof of Lemma~\ref{l:tubeEquivalence}]
First suppose~\eqref{e:gaussianBound} holds, and set
\[
C_G \defeq \sup_{\substack{
      x \in \R^d\,,\\
      r > 0,~
      \kappa \in (0, 1]
    }}
      \frac{\kappa^{\frac{1 + \epsilon}{2}}}{r}
      \P\paren[\big]{
	\mathcal N(x, \kappa I_d) \in B(E, r)
      }.
\]
Fix some $x\in\R^d$, let~$s = \sqrt{\kappa}$ and divide the analysis into two cases.

\emph{Case 1: $0<s\leq1$.}
Then, the density $p(y;x)$ of~$\mathcal N(x,s^2 I_d)$  
obeys a lower bound
\[
p(y,x)\ge (2\pi)^{-d/2}e^{-1/2}s^{-d},~~y\in B(x,s).
\]
This gives an upper bound 
\begin{equation}\label{26oct502}
  |B(E,r)\cap B(x,s)|
    \leq C_d s^d\P\bigl(\mathcal N(x,s^2I_d)\in B(E,r)\bigr)
    \leq C_dC_G r s^{d-1-\epsilon}
  \,,
\end{equation}
which is~\eqref{e:tubeBound}. 

\emph{Case 2: $s > 1$.}
Since $d-1-\eps\ge 0$, to prove that~\eqref{e:tubeBound} also holds
for $s\ge 1$, it suffices to show that
\begin{equation}\label{26oct504}
|B(E,r)|\leq Cr,~~\text{for~$0<r<1$.} 
\end{equation}
Observe that   
since the set~$E$ is bounded, the tube neighborhood
$B(E,1)$ can be covered by finitely many unit balls $B(x_k,1)$, $k=1,\dots,N$, 
with some $x_k\in\R^d$.
Applying the already proved bound~\eqref{26oct502} with~$s=1$ to the balls $B(x_k,1)$ gives
the volume estimate 
\[
|B(E,r)|\leq \sum_{k=1}^N |B(E,r)\cap B(x_k,1)|
\leq NC_dC_Gr,~~\text{for~$0<r<1$,} 
\]
which is~\eqref{26oct504}. 
Thus~\eqref{e:tubeBound} holds with~$r_0=1$.

This concludes the proof that~\eqref{e:gaussianBound} implies~\eqref{e:tubeBound}.
For the converse, we now suppose~\eqref{e:tubeBound} holds for some $r_0>0$.
Let~$X\sim\mathcal N(x,\kappa I_d)$, with~$0<\kappa\leq1$.
We again divide the analysis into two cases.

\emph{Case 1: $r\geq r_0$.}
In this case we simply use the bound~$\P(X\in B(E,r))\leq1$ to see that
\begin{equation}\label{e:GBoundRSmall}
\frac{\kappa^{\frac{1 + \epsilon}{2}}}{r}
      \P\paren[\big]{
	\mathcal N(x, \kappa I_d) \in B(E, r)}\le\frac{1}{r_0}.
\end{equation}

\emph{Case 2: $0<r<r_0$.}
Partition~$\R^d$ into the spherical shells
\[
  S_j(x)=\set{y \st (j-1)\sqrt\kappa\leq|y-x|<j\sqrt\kappa}\,,
\quad j\geq1
  \,.
\]
Assumption~\eqref{e:tubeBound} shows that on each shell we have
\begin{equation}\label{26oct506}
\abs{B(E,r)\cap S_j(x)}\le { \abs{B(E, r) \cap B(x,j\sqrt{\kappa})}}\le 
C{r j^{(d - 1 - \epsilon)}\kappa^{(d - 1 - \epsilon)/2}}.
\end{equation}
Using the upper Gaussian density bound on each shell and~\eqref{26oct506},
we obtain
\begin{align}
  \P(X\in B(E,r))
    &\leq (2\pi\kappa)^{-d/2}
       \sum_{j=1}^\infty |B(E,r)\cap S_j(x)|e^{-(j-1)^2/2}\\
    &\leq Cr\kappa^{-d/2}
       \sum_{j=1}^\infty(j\sqrt\kappa)^{d-1-\epsilon}e^{-(j-1)^2/2}
     \leq \frac{Cr}{\kappa^{(1+\epsilon)/2}}
  \,,
\end{align}
and hence
\begin{equation}\label{e:GBoundRLarge}
\frac{\kappa^{\frac{1 + \epsilon}{2}}}{r}
      \P\paren[\big]{
	\mathcal N(x, \kappa I_d) \in B(E, r)}\le C
	\,.
\end{equation}

Using~\eqref{e:GBoundRSmall} when~$r \leq r_0$, and~\eqref{e:GBoundRLarge} otherwise we obtain~\eqref{e:gaussianBound} as desired.
\end{proof}

\subsection{A tube bound for images (Lemma \ref{l:BVTube})}\label{s:BVTube}

Since the proof of Lemma~\ref{l:BVTube} is shorter than that of Lemma~\ref{l:cutTubeBd}, we present it first.
Before carrying out the details, we briefly describe the main idea.
Let~$h > 0$ and approximate~$f$ by a function~$g_h$ which is piecewise affine on simplices of diameter~$O(h)$.
Standard mollification and piecewise affine interpolation (analogous to the classical Bramble and Hilbert~\cite{BrambleHilbert70} estimates) show that
\begin{equation}\label{e:fMinusGh1}
  \|f-g_h\|_{L^1(Q)} \leq C h^2 \norm{Df}_{\BV(Q')}
  \quad\text{and}\quad
  \|Df-Dg_h\|_{L^1(Q)} \leq C h \norm{Df}_{\BV(Q')}
  \,,
\end{equation}
for a slightly enlarged cube~$Q'$.

Now we compute the tube volume~$\abs{B(g_h(Q), r) \cap B(x, s)}$ by using the area formula on the~$(d-1)$-dimensional term, and bounding the contribution from the lower dimensional faces by~$r^{d - k} h^k $ for~$k \in \set{0, \dots, d-2}$.
Since there are~$O(1/h^{d-1})$ such simplices, the total contribution from lower dimensional faces is bounded by~$C r^2 / h$ when~$r \leq h$.

To bound the tube volume~$\abs{B(f(Q), r) \cap B(x, s)}$, we need to control the approximation error.
We do this by partitioning the domain into cubes of side length~$r$, and use~\eqref{e:fMinusGh1} to show that the cubes where~$f$ and~$g_h$ differ by more than~$r$ at their centers contribute at most~$Ch^2$ to the tube volume.
On the remaining cubes, the Lipschitz bounds allow us to estimate the tube volume using the piecewise affine map~$g_h$.
Balancing~$h^2$ and~$r^2 / h$ errors gives~$h = r^{2/3}$, and the~$O(r^{4/3})$ error term in~\eqref{e:BVTube}.
We now carry out the details.

\iffalse
%
%
%
%
%
%
%
%
%
%
%
%
%
%
%
%
%
%
%
%
%
%
%
%
%
%
%
%
%
%
%
%
%
%
%
%
%
\fi

\begin{proof}[Proof of Lemma~\ref{l:BVTube}]
All approximations below are made on fixed intermediate cubes 
between~$Q$ and~$Q'$. In particular, the convolution and the meshes used below
stay inside the larger cube~$Q'$ when their scales are sufficiently small.
Let
\begin{equation}\label{26oct539}
r\leq h\ll1,
\end{equation}
choose a nonnegative radially symmetric
mollifier~$\zeta\in C_c^\infty(\R^{d-1})$, supported in~$B(0,1)$, with integral one and 
set
\begin{equation}\label{e:mollifiedMap}
  \zeta_h(y)= \frac{1}{h^{d-1}} \zeta\paren[\Big]{ \frac{y}{h}}
  \,,\qquad f_h=f*\zeta_h
  \,.
\end{equation}
Interpolate the smoothed function~$f_h$ on a regular simplicial mesh of size~$h$, and denote
the resulting piecewise affine map by~$g_h$.
This approximation satisfies the Lipschitz bound 
\begin{equation}\label{e:BVInterpolation}
  \operatorname{Lip}(g_h)\leq C\operatorname{Lip}(f_h)\le C\operatorname{Lip}(f)
  \,,
\end{equation}
and the gradient bounds (see for instance~\cite{EvansGariepy92} Sections 4.6, 5.2.2)
\begin{equation}\label{e:DfDfhBv}
 \|Df-Df_h\|_{L^1(Q)}\leq Ch\|Df\|_{BV(Q')},~~ \|Df_h-Dg_h\|_{L^1(Q)}\leq Ch\|Df\|_{BV(Q')},
\end{equation}
Indeed, we have 
\begin{equation}\label{e:DfBVBound1}
\|Df-Df_h\|_{L^1(Q)}\leq Ch\|Df\|_{BV(Q')},
\end{equation}
and the second bound in~\eqref{e:DfDfhBv} is proved similarly as well. 
We will also use the quadratic bound 
\begin{equation}\label{26oct526}
  \|f-g_h\|_{L^1(Q)}\leq Ch^2\|Df\|_{BV(Q')}
  \,.
\end{equation}
To see that~\eqref{26oct526} holds, note that the radial symmetry of $\zeta(x)$ gives
\[
  f_h(x)-f(x)=\frac12\int_{\R^{d-1}}
    \bigl(f(x-hy)+f(x+hy)-2f(x)\bigr)\zeta(y)\,dy
  \,.
\]
Now, we note
on each fixed intermediate cube, that
\begin{align}
  \|f(\cdot+h)+f(\cdot-h)-2f\|_{L^1(Q)}
  &=\norm[\Big]{ \int_0^1
    \bigl(Df(x+th)-Df(x-th)\bigr)\cdot h\,dt
    }_{L^1(Q)}
  \\
    &\leq C|h|^2\|Df\|_{BV(Q')}
  \,.
\end{align}
This gives the bound
\begin{equation}\label{26oct527}
  \|f-f_h\|_{L^1(Q)}\leq Ch^2\|Df\|_{BV(Q')}
  \,.
\end{equation}
On the other hand, we also have a standard bound for the simplicial approximations:
\begin{equation}\label{e:ghMinusfhL1}
  \|g_h-f_h\|_{L^1}\leq Ch^2\|Df_h\|_{BV(Q')}\le Ch^2\|Df\|_{BV(Q')}. 
  \,.
\end{equation}
Combining~\eqref{26oct527} and~\eqref{e:ghMinusfhL1} gives~\eqref{26oct526}.

Next cover~$Q$ by cubes $Q_j$ of side length~$r$ from a translated grid, with centers $u_j$
Averaging over translations and using the quadratic estimate~\eqref{26oct526}
allows us to choose the grid so that
\begin{equation}\label{26oct528} 
r^{d-1}\sum_j|f(u_j)-g_h(u_j)|\leq Ch^2\|Df\|_{BV(Q')}
  \,.
\end{equation} 
Let $\cB$ be the collection of centers $u_j$ such that
\begin{equation}\label{26oct532}
|f(u_j)-g_h(u_j)|>r.
\end{equation}
It follows from~\eqref{26oct528} that
\begin{equation}\label{26oct529}
N_B \defeq \#\{j:~Q_j\in\cB\}\le \frac{Ch^2}{r^{d}}\|Df\|_{BV(Q')}.
\end{equation}
As $f$ is Lipschitz, we know that 
the image $f(Q_j)$ of each cube in $\cB$ is contained in a ball of radius~$Cr$.  That is, we have
\begin{equation}\label{26oct530}
|B(f(Q_j),r)|\le Cr^d,~~\hbox{ for each $Q_j\in\cB$.} 
\end{equation}
Together,~\eqref{26oct529} and~\eqref{26oct530} give 
\begin{equation}\label{26oct531}
\sum_{Q_j\in\cB}|B(f(Q_j),r)|\le Cr^d N_B\le Ch^2. 
\end{equation}

Let now $\cG$ be the collection of the remaining cubes, not in $\cB$. The centers of $Q_j\in\cG$
satisfy the converse of~\eqref{26oct532}:
\begin{equation}\label{26oct533}
|f(u_j)-g_h(u_j)|\le r.
\end{equation}
The uniform Lipschitz bound on $f$ implies that then 
\begin{equation}\label{26oct534}
|f(x)-g_h(u_j)|\le Cr,~~\text{ for all $x\in Q_j$.}
\end{equation}
Therefore, if the tubular neighborhood $B(f(Q_j),r)$ of the image $f(Q_j)$ of 
a cube $Q_j\in\cG$ has a non-empty intersection with a ball~$B(x,s)$, 
then its center~$u_j$ satisfies
\begin{equation}\label{26oct535}
g_h(u_j)\in B(x,s+Cr).
\end{equation}
Let $\cG_S$ be the collection of all cubes $Q_j\in\cG$ whose center $u_j$ lies in a given
simplex~$S$ of the~$h$-mesh, and let 
\begin{equation}\label{26oct542}
P_S(r)=\bigcup_{Q_j\in\cG_S}B(f(Q_j),r).
\end{equation}
It follows from~\eqref{26oct534} and~\eqref{26oct535} that
\begin{equation}\label{26oct536}
P_S(r)\cap B(x,s) \subset B(K_S,Cr).
\end{equation}
Here, we have set
\[
K_S=g_h(S)\cap B(x,s+Cr)
  \,.
\]
Indeed, if $y\in P_S(r)\cap B(x,s)$, then there exists $z\in Q_j\in\cG_S$ such that $|y-f(z)|<r$ and also~$|y-x|<s$. 
It follows from~\eqref{26oct534} that $|y-g_h(u_j)|\le Cr$ and also from~\eqref{26oct535} that~$|g_h(u_j)-x|\le s+Cr$. 
This gives~\eqref{26oct536}.

 Since $g_h$ is an affine function on the simplex element $S$ of size $O(h)$, the set $K_S$ is a convex set of dimension at most~$d-1$ and diameter at
most~$Ch$. Since each side of the simplex has at most size $h\ge r$ because of our choice in~\eqref{26oct539}, the tube formula for convex sets, including its
lower dimensional terms, gives
\begin{equation}\label{26oct537}
\begin{aligned}
  |B(K_S,Cr)|
    &\leq Cr\mathcal H^{d-1}(K_S)
      +C\sum_{j=0}^{d-2}r^{d-j}h^j=Cr\mathcal H^{d-1}(K_S)
      +Cr^d\frac{h^{d-2}}{r^{d-2}}\\
      &=Cr\mathcal H^{d-1}(K_S)+r^2h^{d-2}
  \,.
\end{aligned}
\end{equation}
The first term in the right side of~\eqref{26oct537}  is bounded using the area formula by
\begin{equation}\label{26oct543}
\mathcal H^{d-1}(K_S)
  \leq \int_{\{y\in S:~g_h(y)\in B(x,s+Cr)\}}J_{d-1} g_h(y).
\end{equation}
Here, the Jacobian $J_{d-1}g_h$ is defined as in~\eqref{e:tanjentialJacobianDef}:
\begin{equation}\label{26oct540}
  J_{d-1}g_h=\sqrt{\det(Dg_h^T Dg_h)}
  \,,
\end{equation}
We now estimate, using~\eqref{26oct542}--\eqref{26oct543}
\begin{equation}\label{26oct545}
\begin{aligned}
\abs[\bigg]{
\bigcup_{Q_j\in\mathcal G}
\bigl(B(f(Q_j),r)\cap B(x,s)\bigr)
}
&\leq \sum_S |P_S(r)\cap B(x,s)|
\le \sum_{S}|B(K_S,Cr)|\\
&\le r \int_{\{y\in Q':~g_h(y)\in B(x,s+Cr)\}}J_{d-1} g_h(y)+r^2h^{d-2}\#S.
\end{aligned}
\end{equation}
The total number of simplices is bounded by
\begin{equation}\label{26oct544}
\#S\le C h^{-(d-1)}.
\end{equation}
Using this bound in~\eqref{26oct545} and recalling the estimate~\eqref{26oct531} for the cubes in the collection $\cB$, we arrive at 
 \begin{equation}\label{e:BVTubeIntermediate}
  |B(f(Q),r)\cap B(x,s)|
    \leq Cr\int_{\{y\in Q':g_h(y)\in B(x,s+Cr)\}}J_{d-1} g_h(y)
       +Ch^2+\frac{Cr^2}{h}
  \,.
\end{equation}
It remains to replace $J_{d-1}g_h$ by $J_{d-1}f$  in~\eqref{e:BVTubeIntermediate}, and also $g_h$ by $f$ in the domain of integration,  
in order to obtain~\eqref{e:BVTube}.
However, we know from~\eqref{e:DfDfhBv} that 
\begin{equation}\label{e:Jdm1fg1}
\|J_{d-1} f-J_{d-1} g_h\|_{L^1}\leq Ch.
\end{equation}
Furthermore, the quadratic bound~\eqref{26oct526} implies that, except on a set of measure at most~$Ch^2/r$, the inequality
$g_h(u)\in B(x,s+Cr)$ implies~$f(u)\in B(x,s+C'r)$. In addition, the Jacobian $J_{d-1}f$ is uniformly bounded because of~\eqref{e:BVInterpolation}. 
Together, these estimates, combined with~\eqref{e:Jdm1fg1},  show that 
\[
  r\int_{\{y\in Q':g_h(y)\in B(x,s+Cr)\}}J_{d-1} g_h
    \leq r\int_{\{y\in Q':f(y)\in B(x,s+C'r)\}}J_{d-1} f
       +Crh+Ch^2
  \,.
\]
Substituting this in~\eqref{e:BVTubeIntermediate} gives
\begin{equation}\label{26oct549}
  |B(f(Q),r)\cap B(x,s)|
    \leq Cr\int_{\{y\in Q':f(y)\in B(x,s+Cr)\}}J_{d-1} f(y)
       +Ch^2+\frac{Cr^2}{h}+Crh
  \,.
\end{equation}
As $r\le h$ because of the restriction~\eqref{26oct539} on $h$, the last term in the right side above is dominated by $Ch^2$.
Now, choosing
$h=r^{2/3}$, which satisfies the restriction~\eqref{26oct539} on $h$,  proves~\eqref{e:BVTube}.
\end{proof}

\section{Large tube bounds (Lemma \ref{l:cutTubeBd})}\label{s:cutTubeBd}

\subsection{Proof of Lemma~\ref{l:cutTubeBd}}

We devote this section to proving Lemma~\ref{l:cutTubeBd}.
If~$H \colon \R^d \to \R$ is a~$C^1$ function, then the level set~$\set{H = 0}$ satisfies the tube bound~\eqref{e:tubeBound} with~$\epsilon =0$ \emph{provided} one has a non-degeneracy condition on~$\grad H$.
The key insight that leads to the proof of Lemma~\ref{l:cutTubeBd} is that for non-constant polynomials, no such non-degeneracy condition is required.
That is, if~$H$ is a non-constant polynomial, then~$\set{H = 0}$ satisfies the tube volume estimate~\eqref{e:tubeBound} with~$\epsilon = 0$, and a prefactor that only depends on~$d$, $\deg(H)$, \emph{without} any dependence on the coefficients or requiring any non-degeneracy condition on~$\grad H$.
This is our first lemma.

\begin{lemma}\label{l:polynomialZeroTube}
There exists a constant $C_d$ such that if $H$ is a nonzero polynomial on $\R^d$ of degree $D\geq1$, 
then for $0<r\leq s$, we have 
\begin{equation}\label{e:elementaryPolynomialTube}
  |B(\{H=0\},r)\cap B(x,s)|\leq C_dD^d r s^{d-1}.
\end{equation}
\end{lemma}

In light of Lemma~\ref{l:polynomialZeroTube}, the strategy that will be used to prove Lemma~\ref{l:cutTubeBd} is as follows.
Fix~$s > 0$, choose~$N = 1/s^{O(\epsilon)}$, and find a non-constant polynomial~$H_N$ with degree that
is polynomial in~$N$, so that~$\cut_R(\mathcal M)$ is contained in a set that is well approximated by~$\set{H_N = 0}$.
Then control the approximation errors in the area formula to obtain~\eqref{e:Jdm1FL}, and use~\eqref{e:elementaryPolynomialTube} to obtain~\eqref{e:tubeVolumeRlarge}.

To carry out this strategy, we need a few preliminary lemmas on polynomial approximation.
We begin by briefly introducing some terminology, and refer the reader to~\cite{BasuPollackEA06} for a more comprehensive description.
A set is called a \emph{real algebraic set} if it is the zero set of a polynomial.
A set~$E\subseteq\R^n$ is called \emph{semialgebraic} if it is a finite union of sets defined by finitely many polynomial equalities and inequalities.
A map~$f:E\to\R^m$ is called \emph{semialgebraic} if its graph is semialgebraic.
A \emph{smooth semialgebraic manifold} is a smooth embedded manifold that is also a semialgebraic set.
Finally, a \emph{smooth semialgebraic chart} is a smooth diffeomorphism from an open subset of~$\R^k$ onto an open subset of the manifold, which is also a semialgebraic map.

The first step in the aforementioned construction is to approximate the image $F_L(\cU)$ by semi-algebraic manifolds. 
This is accomplished by the next lemma. 
\begin{lemma}\label{l:cutAlgebraicApproximation}
There is a compact smooth semialgebraic~$(d-1)$-dimensional manifold~$\mathcal B$,
$N_0 \in \N$,
and  maps~$F:\cB\to\R^d$ and $F_N:\mathcal B\to\R^d$, $N\ge N_0$, 
with the following properties. 
\begin{enumerate}
\item
There exists a smooth diffeomorphism~$G:\cB\to\cU$ such that~$F=F_L\circ G$, and hence $F(\mathcal B)=F_L(\mathcal U)$.

\item
There exists~$K \in \N$, such that for every~$i \in \set{1, \dots, K}$ there exist open cubes~$Q_i \Subset Q_i' \Subset \R^{d-1}$ and smooth semialgebraic charts~$\psi_i:Q_i'\to\mathcal B$ such that~$\mathcal B=\bigcup_{i=1}^K\psi_i(Q_i)$.
Moreover, there exists a constant $C > 0$ such that for all~$1\le i\le K$ and $N\ge N_0$ the maps~$F_i=F\circ\psi_i$ and~$F_{N,i}=F_N\circ\psi_i$ obey uniform bounds
\begin{align}\label{26oct702}
\|F_i\|_{\Lip(Q_i')} &\le C,
&\|F_{N,i}\|_{\Lip(Q_i')} &\le C,
\\
\label{26oct704}
\|DF_i\|_{BV(Q_i')}&\le C,
&\|DF_{N,i}\|_{BV(Q_i')}&\le C
\,.
\end{align}

\item
For every~$j\geq1$, there is~$C_j>0$ such that for all $N\ge N_0$ we have 
\begin{equation}\label{e:cutApproximationError}
  \|F-F_N\|_{L^\infty(\mathcal B)}\leq \frac{C_j}{N^j}, %
  \qquad
  \sum_{i=1}^K\|DF_i-DF_{N,i}\|_{L^1(Q_i')}
       \leq \frac{C_j}{N^{j/2}} %
  \,.
\end{equation}
\end{enumerate}
\end{lemma}
With the result of Lemma~\ref{l:cutAlgebraicApproximation} in hand, we will think of $F_L(\cU)=F(\cB)$ as approximately~$F_N(\cB)$. 
In addition, as we will see in the course of the proof of that lemma, the degrees of the polynomials involved in the definition of $F_N$ will be controlled
by $CN^P$ with some fixed $P\ge 1$.
This is important because of the claim of Lemma~\ref{l:polynomialZeroTube} and the following.  
\begin{lemma}\label{l:endpointPolynomial}
There exists $P_2 \in \N$ such that for all~$N \geq N_0$ there exists a non-constant polynomial $H_N:\R^{d}\to\R$, of degree at most $CN^{P_2}$, such that 
\begin{equation}
F_N(\cB)\subseteq \set{H_N = 0}
\,.
\end{equation}
\end{lemma}
We will also need the following estimate that will lead to the claim~\eqref{e:Jdm1FL} of Lemma~\ref{l:cutTubeBd}. 
\begin{lemma}\label{l:algebraicCutArea}
There exist $C,P_3>0$ such that for all $x\in\R^d$ and $s>0$,
\begin{equation}\label{e:algebraicCutArea}
  \sum_{i=1}^K\int_{Q_i'\cap F_{N,i}^{-1}(B(x,s))}
       J_{d-1}F_{N,i}(y)dy\leq CN^{P_3}s^{d-1}.
\end{equation}
\end{lemma}

Momentarily postponing the proofs of Lemmas~\ref{l:polynomialZeroTube}-\ref{l:algebraicCutArea}, we now prove Lemma~\ref{l:cutTubeBd}.
\begin{proof}[Proof of Lemma~\ref{l:cutTubeBd}]
By Lemmas~\ref{l:endPointMap} and~\ref{l:cutAlgebraicApproximation},
\begin{equation}\label{26oct1610-2}
\cut_R(\mathcal M)\subseteq F_L(\cU)=F(\cB).
\end{equation}
Fix~$\eps>0$ and~$d'\ge 1$, and %
choose~$\alpha>0$ sufficiently small that
\begin{equation}\label{26oct1606}
\alpha\max\{dP_2,P_3\}<\farc{\eps}{2}. 
\end{equation}
and, for~$0<s\leq s_0$, put~
\begin{equation}\label{26oct1604}
N=\lceil s^{-\alpha}\rceil.
\end{equation}
Here, $P_2$ is as in Lemma~\ref{l:endpointPolynomial} and $P_3$ is as in Lemma~\ref{l:algebraicCutArea}. 
Choose~$j$ sufficiently large that
\[
\alpha j>d',
\qquad
\alpha j/2>d-1.
\]
By~\eqref{e:cutApproximationError}, we have
\begin{equation}\label{e:cutApproximationSmall}
\delta_N\defeq\|F-F_N\|_{L^\infty(\cB)}=o(s^{d'}),
\qquad
\sum_{i=1}^K\|DF_i-DF_{N,i}\|_{L^1(Q_i')}
=o(s^{d-1}).
\end{equation}

We first estimate the parametrized area of~$F_L$.
The uniform Lipschitz bounds in~\eqref{26oct702} imply
\begin{equation}\label{26oct1602}
|J_{d-1}F_i-J_{d-1}F_{N,i}|
\leq C|DF_i-DF_{N,i}|.
\end{equation}
Moreover, we have, by~\eqref{e:cutApproximationSmall},  
\begin{equation}\label{26oct1608}
F_i^{-1}(B(x,s))
\subseteq F_{N,i}^{-1}(B(x,s+\delta_N)).
\end{equation}
Using the finite atlas from Lemma~\ref{l:cutAlgebraicApproximation},
the area bound~\eqref{e:algebraicCutArea}, and the approximations~\eqref{e:cutApproximationSmall} and~\eqref{26oct1602}, we obtain
\begin{align*}
\int_{F_L^{-1}(B(x,s))}J_{d-1}F_L\,d\cH^{d-1}
&\leq C\sum_{i=1}^K
\int_{Q_i'\cap F_i^{-1}(B(x,s))}
J_{d-1}F_i(y)\,dy\\
&\le C\sum_{i=1}^K
\int_{Q_i'\cap F_{N,i}^{-1}(B(x,s+\delta_N))}
J_{d-1}F_{N,i}(y)\,dy+o(s^{d-1})\\
&\leq CN^{P_3}(s+\delta_N)^{d-1}+o(s^{d-1})\leq Cs^{d-1-\eps}.
\end{align*}
Here, we have also used the uniform bounds on the chart Jacobians, the choice of $N$ in~\eqref{26oct1604} and that of $\alpha$ in~\eqref{26oct1606},
as well as the inclusion~\eqref{26oct1608}. This proves~\eqref{e:Jdm1FL}. 

We next obtain the tube volume estimate~\eqref{e:tubeVolumeRlarge}. We still choose $N$ as in~\eqref{26oct1606}-\eqref{26oct1604}. 
Suppose that
\begin{equation}\label{26oct1612}
s^{d'}\leq r\leq s,
\end{equation}
so that, in particular
\begin{equation}\label{26oct1610}
\delta_N\le r.
\end{equation}
By~\eqref{26oct1610} and~\eqref{e:cutApproximationSmall}, we have %
\begin{equation}\label{26oct1618}
B(\cut_R(\mathcal M),r)\subseteq B(F(\cB),r)
\subseteq B(F_N(\cB),r+\delta_N)
\subseteq B(\{H_N=0\},r+\delta_N).
\end{equation}
Here,~$H_N$ is the polynomial from Lemma~\ref{l:endpointPolynomial}, with
\begin{equation}\label{26oct1614}
\deg H_N\le CN^{P_2}\le Cs^{-\alpha P_2}.
\end{equation}
Lemma~\ref{l:polynomialZeroTube} implies that
\begin{equation}
\begin{aligned}
  |B(\{H=0\},r+\delta_N)\cap B(x,s+r+\delta_N)|&\leq C_d(\deg H_N)^d (r+\delta_N) (s+r+\delta_N)^{d-1}\\
  &\le C_d s^{-\alpha dP_2}rs^{d-1}\le C_d rs^{d-1-\eps}.
\end{aligned}
\end{equation}
We used above~\eqref{26oct1606},~\eqref{26oct1612},~\eqref{26oct1610}, and~\eqref{26oct1614},  
Together with~\eqref{26oct1618}, this gives
\begin{align*}
|B(\cut_R(\mathcal M),r)\cap B(x,s)|\leq Crs^{d-1-\eps}.
\end{align*}
This proves~\eqref{e:tubeVolumeRlarge} and the lemma.
\end{proof}

\subsection{Tube neighborhoods of real algebraic sets (Lemma~\ref{l:polynomialZeroTube})}\label{s:polynomialZeroTube}

We now prove Lemma~\ref{l:polynomialZeroTube}.
The idea behind the proof is to divide~$\R^d$ into a grid of cubes of side length~$r$, and note that the tube estimate~\eqref{e:elementaryPolynomialTube} can be obtained by counting cubes intersecting the zero-level set.
We use induction on the dimension to bound the number of cubes whose boundaries intersect the zero level set, and B\'ezout's theorem and a polynomial perturbation argument to bound the remainder.

\begin{proof}[Proof of Lemma~\ref{l:polynomialZeroTube}]
Let~$r > 0$ and~$a \in \R^d$ and consider the grid~$a + r \Z^d$.
For almost every~$a \in \R^d$, for every~$k \in \set{0, \dots, d}$, the polynomial~$H$ does not vanish identically on any~$k$-dimensional face of any grid cube.
For simplicity, and without loss of generality we now assume~$a = 0$ and fix the grid to be~$r\Zm^d$ . 

Consider any block $\cZ_M$  of~$M^d$ grid cubes and let $N(d, D, M)$ be the number of those cubes that intersect the level set~$\{H=0\}$.
We claim that  there exists a dimensional constant~$C_d$ such that 
\begin{equation}\label{e:polynomialGridCount}
N(d,D,M) \leq C_d D^d M^{d-1}.
\end{equation}
 
We prove \eqref{e:polynomialGridCount} by induction
on the dimension $d\ge 1$. Once again, without loss of generality we may assume that $\cZ_M=[0,rM]^d$.
For~$d=1$, the bound~\eqref{e:polynomialGridCount}
is
\[
N(1, D, M) \le C_1D,
\]
and holds with $C_1=1$ simply because a polynomial of degree
$D$ can have at most $D$ roots in one dimension. 

Suppose the claim holds in dimension~$d-1$.
Of the cubes in $\cZ_M$ that intersect~$\set{H = 0}$, first consider cubes whose boundaries also intersect~$\set{H = 0}$.
Fix $1\le k\le d$, and a coordinate hyperplane $\Pi_{jk} \defeq \{x_k=jr\}$, with some~$0\le j\le M$.
Consider the cubes in $\cZ_M$ such that
  (1)~one of their faces is contained in the hyperplane $\Pi_{jk}$, and
  (2) that particular face intersects the level set~$\set{H = 0}$.
These faces are themselves $(d-1)$-dimensional cubes on the grid~$r \Z^{d-1}$.
By the induction hypothesis,  for each $1\le k\le d$ and $0\le j\le M$ fixed, there are at most $N(d-1,D,M)$ such faces. 
Each of these faces corresponds to two cubes in~$\cZ_M$ whose faces intersect $\set{H = 0}$, while the total
number of hyperplanes $\Pi_{jk}$ is $(M+1)d$.
Thus the total number of cubes in $\cZ_M$ whose faces intersect $\set{H = 0}$ is at most
\begin{equation}\label{e:NFaceIntersect}
2d(M+1)N(d-1,D,M)\le 
2d(M+1)C_{d-1}D^{d-1}M^{d-2}
      \leq C_d D^{d-1}M^{d-1}
      \,.
\end{equation}

Now consider cubes in~$\mathcal Z_M$ whose interiors intersect~$\set{H = 0}$ but faces do not.
Let $Q_m$ be any such cube.
The polynomial $H^2(x)$ has a zero in the interior and is strictly positive on the boundary.
That is, $H^2(x)$ is a polynomial of $d$ variables of degree $D$ that is positive at the boundary of $Q_m$ and has a minimum $x_0$ inside $Q_m$.
Therefore, there exists~$\eps>0$ such that
the same is true for the perturbed polynomial
\begin{equation}
H_\eps(x)=H^2(x)+\eps \sum_{k=1}^d x_k^{2D+1}.
\end{equation}
Hence, in each of $Q_m$ there is an isolated zero of the system
\begin{equation}
\partial_k [H^2(x)] + \eps (2D+1)x_k^{2D}=0,~~1\le k\le d. 
\end{equation}
This is a system of $d$ unknowns of polynomial equations of degree $2D$.
By B\'ezout's theorem it can have no more than $(2D)^d$ isolated solutions.
Therefore, the number of cubes in~$\mathcal Z_M$ whose interiors intersect~$\set{H = 0}$ but faces do not is at most~$(2D)^d$.
Combined with with~\eqref{e:NFaceIntersect}, this means that the total number of cubes in~$\mathcal Z_M$ that intersect $\set{H = 0}$ is at most
\begin{equation}
N(d, D, M) \le C_d D^{d-1}M^{d-1}+(2D)^d\le C_dD^{d}M^{d-1}
\,,
\end{equation}
which proves~\eqref{e:polynomialGridCount}.
\smallskip

Let us now consider the intersection $B(\{H=0\},r)\cap B(x,s)$ with $0<r\le s$.  
Let us apply the bound~\eqref{e:polynomialGridCount} in a grid block of the mesh $r\Zm^d$, that contains
the larger ball~$B(x,s+r)$ and has 
at most~$M_{r,s}=C_d s/r$ cubes on each side.
Every point in~$B(x,s)$ at distance less than~$r$ from the level set $\set{H = 0}$ 
lies in a cube within a bounded number of grid
steps of one of the cubes intersecting $\set{H = 0}$. Each  cube has
at most~$C_d$ such neighbors. The volume of each cube is $r^d$. Thus, we have
\begin{equation}
|B(\{H=0\},r)\cap B(x,s)|\le C_d r^dN(d,D,{M_{r,s}})\le C_d r^d D^d\Big(\frac{s}{r}\Big)^{d-1}
=C_dD^drs^{d-1}\,,
\end{equation}
proving \eqref{e:elementaryPolynomialTube} as desired.
\end{proof}

\subsection{The image of the truncated normal map}\label{s:endPointMap}

We now prove Lemma~\ref{l:endPointMap}.
Studying regularity of~$\cut_R(\mathcal M)$ through the image of the truncated normal map is a well known technique (see for instance Itoh and Tanaka~\cite{ItohTanaka01}).
While the proof of Lemma~\ref{l:endPointMap} is contained in~\cite{ItohTanaka01}, the proof itself is short, and we present it here for completeness.
 
\begin{proof}[Proof of Lemma~\ref{l:endPointMap}]
Take some $x\in\cut_R(\cM)$. Let $x_1$ be one of the closest points to $x$ on $\cM$ and
set~$\nu_1=(x-x_1)/|x-x_1|$ and $t_1=|x-x_1|$. Then, the point 
$(x_1,\nu_1)$ is in $\cU$ and~$t_1=t_*(x_1,\nu_1)$, with $t_*$ defined in~\eqref{26oct512}. 
 We also note
that, as $x\in\cut_R(\cM)$, we know that $|x|\le R$. It follows that  
\[
t_1=|x-x_1|\le R+|x_1|\le R+\sup_{p\in\mathcal M}|p|<L.
\]
We used assumption~\eqref{26oct516} on $L$ in the last step above. 
As a consequence, 
\[
t_L(x_1,\nu_1)=t_*(x_1,\nu_1)=t_1,
\]
which means that
\[
F_L(x_1,\nu_1)=x_1+t_L(x_1,\nu_1)\nu=x_1+t_1\nu =x.
\]
Therefore, $x$ lies in the image $F_L(\cU)$. 
\end{proof}

\subsection{An alternative representation for the truncated normal map} 
 
In order to study the regularity of  the truncated normal map, we will express it in terms of a supremum of a family of smooth scalar-valued functions.

First, we define
\begin{equation}\label{26oct202}
q(z,p,\nu)=\frac{2\nu\cdot(z-p)}{|z-p|^2},~~z\in\cM,~~z\neq p,~(p,\nu)\in\cU.
\end{equation}
Given some $(p,\nu)\in\cU$,
we note that 
\begin{equation}
|p+t\nu-z|^2=t^2+|z-p|^2-2t\nu\cdot(z-p)=t^2+|z-p|^2(1-tq(z,p,\nu)),~~\hbox{for any $z\in\cM$, $z\neq p$.} 
\end{equation}
Therefore,   $p$ is the closest point on $\cM$ to all points of the form
$p+t\nu$, with $t\ge 0$, as long as, 
\[
t\bar q(p,\nu)\leq1,
\]
with 
\begin{equation}\label{e:reciprocalCutTime}
  \bar q(p,\nu)=\max\big\{0,\bar q_1(p,\nu)\big\},~~
  \bar q_1(p,\nu)=\sup_{z\in\cM\setminus\{p\}}q(z,p,\nu). 
\end{equation}
We also set
\begin{equation}\label{26oct210}
\bar q_L(p,\nu)=\max\{\bar q(p,\nu),1/L\}=\max\big\{1/L,\bar q_1(p,\nu)\big\}.
\end{equation}
It follows that
\begin{equation}\label{e:truncatedReciprocal}
  t_*(p,\nu) =\frac{1}{\bar q(p,\nu)},~~~ t_L(p,\nu)=\min\{t_*(p,\nu),L\}=
  \frac{1}{\max\{\bar q(p,\nu),1/L\}}=\frac{1}{\bar q_L(p,\nu)},
 \end{equation}
 and
\begin{equation}\label{26oct208}
 F_L(p,\nu)=p+\frac{1}{\bar q_L(p,\nu)}\nu. %
\end{equation}
Thus, the regularity of $F_L(p,\nu)$ comes from the regularity of $\bar q_L(p,\nu)$, which, in turn, comes from the regularity 
of $q(z,p,\nu)$ defined in~\eqref{26oct202}.  

At a first glance, $q(z,p,\nu)$ seems to have a $O(|z-p|^{-1})$ singularity at $z=p$. This, however, 
cancels out because $\nu\perp T_p(\cM)$. 
To see that, consider a sufficiently small coordinate chart $\cQ\subset\R^m$ of $\cM$ near the point $p$, 
and represent $p=\psi(u)$ and~$z(s;\omega)=\psi(u+s\omega)$, with~$u\in\cQ$,~$|\omega|=1$, and $s\in\R$ 
is sufficiently small. Here, $\psi$ is the coordinate map on $\cQ$. 
Then, we have 
\begin{equation}\label{26oct110}
\begin{aligned}
z(s;\omega)-p=\int_0^1 s D\psi(u+st\omega)\omega\,dt.
\end{aligned}
\end{equation}
Since $\nu\perp T_p(\cM)$, we have 
\begin{equation}\label{26oct204}
\nu\cdot D\psi(u)\omega=0.
\end{equation}
Using~\eqref{26oct110} and~\eqref{26oct204}, the inner product that appears in the numerator 
in \eqref{26oct202} can be written as 
\begin{equation}\label{26oct108}
\begin{aligned}
\nu\cdot(z(s;\omega)-p)&=\int_0^1 s \nu\cdot D\psi(u+st\omega)\omega\,dt=
\int_0^1 s \nu\cdot [D\psi(u+st\omega)-D\psi(u)]\omega\,dt\\
&=\sum_{p,k=1}^m \int_0^1 s \nu_p \big[\partial_k\psi_p(u+st\omega) - \partial_k{\psi_p(u)} \big]\omega_k \, dt\\
&=
\sum_{p,k,j=1}^m \int_0^1 s \nu_p \omega_k\int_0^1 \partial_k \partial_j \psi_p(u+stt'\omega) \, s t \omega_j \, dt \, dt'
\\
&=s^2\sum_{p,k,j=1}^m \int_0^1  \nu_p \omega_k\int_0^t \partial_k \partial_j \psi_p(u+st'\omega) \omega_j \, dt' \, dt\\
&=s^2\sum_{p,k,j=1}^m \nu_p \omega_k\omega_j\int_0^1 (1-t) \partial_k \partial_j \psi_p(u+st\omega) \, dt.
\end{aligned}
\end{equation}
Recalling~\eqref{26oct110}, we see that  
\begin{equation}\label{e:resolvedCutKernel}
q(z(s;\omega),p,\nu)=\frac{2\nu\cdot(z(s)-p)}{|z(s)-p|^2}=  \frac{2\displaystyle\int_0^1(1-t)
       \nu\cdot D^2\psi(u+st\omega)[\omega,\omega]\,dt}
       {\Big|\displaystyle\int_0^1
          D\psi(u+st\omega)\omega\,dt\Big|^2},
\end{equation}
with the notation that is clear by a
comparison of~\eqref{e:resolvedCutKernel} to~\eqref{26oct108}.
The denominator in~\eqref{e:resolvedCutKernel} is uniformly positive for $s$ close to $s=0$, and the 
expression in the right side of~\eqref{e:resolvedCutKernel}
can be extended to~$s=0$. In particular, we see from~\eqref{e:resolvedCutKernel} that,  as long as the chart map~$\psi:\cQ\to\cM$ is $C^4$, the function 
\[
w(s,\omega,p,\nu)=q(z(s;\omega),p,\nu)
\]
is $C^2$ in $s$ and $\omega$. As the function $q(z,p,\nu)$ is also smooth for $z\in\cM$ away from $p$, 
we conclude that
$q(z,p,\nu)$ is $C^2$ in all its arguments if $\cM$ is a $C^4$-manifold. 

\subsection{Stability of the truncated normal map}

The proof of Lemma~\ref{l:cutAlgebraicApproximation} requires some regularity and stability of the truncated normal map~$F_L$, which we now describe.
Let~$\iota \colon \mathcal M \to \R^d$ be the inclusion map, and~$\varphi \colon \mathcal M \to \R^d$ be a smooth embedding.
When~$\norm{\varphi-\iota}_{C^1}$ is sufficiently small, set
\begin{equation}\label{26oct520}
  \nu^\varphi(p,\nu)
    =\frac{P^\perp_{\varphi,p}\nu}{|P^\perp_{\varphi,p}\nu|},
  \qquad
  I_\varphi(p,\nu)=\bigl(\varphi(p),\nu^\varphi(p,\nu)\bigr).
\end{equation}
Here,~$P^\perp_{\varphi,p}$ is the orthogonal projection onto the normal space 
of the image $\varphi(\mathcal M)$ at the point~$\varphi(p)$. 
For $\vphi$ close to $\iota$, the denominator is uniformly positive and $I_\varphi$
identifies the unit normal bundles.
Using this identification, we may view the truncated normal map of the embedded 
manifold~$\varphi(\mathcal M)$ as a map on the unit normal bundle of~$\mathcal M$.
We let~$F^\varphi_L \colon \mathcal U \to \R^d$ denote this map.
Our next lemma concerns stability of the truncated normal map under small~$C^4$ perturbations of~$\mathcal M$.

\begin{lemma}\label{l:truncatedNormalRegularity}
There exist constants~$C,\eta_0>0$ such that for every~$\eta\in(0,\eta_0)$ and every smooth embedding~$\varphi\colon\mathcal M\to\R^d$ satisfying
\begin{equation}\label{e:C4distance}
  \|\varphi-\iota\|_{C^4}<\eta,
\end{equation}
the corresponding truncated normal map~$F_L^\varphi$ is Lipschitz with $DF_L^\varphi\in BV$, and
\begin{align}
  \label{e:FLLipBV}
  \norm{F_L^\varphi}_{\Lip(\mathcal U)} & \leq C
  &
  \norm{D F_L^\varphi}_{\BV(\mathcal U)}
  &\leq C
  \,,
  \\
  \label{e:normalMapStability}
  \|F_L^\varphi-F_L^\iota\|_{L^\infty(\mathcal U)}
    &\leq C\eta,
  &
  \|DF_L^\varphi-DF_L^\iota\|_{L^1(\mathcal U)}
    &\leq C\sqrt{\eta}.
\end{align}
\end{lemma}

We clarify the derivative norms in~\eqref{e:FLLipBV} and~\eqref{e:normalMapStability} are computed locally, in charts.
Explicitly, fix finitely many smooth coordinate charts~$\set{\psi_i}$ of~$\mathcal U$, and choose cubes~$Q_i \Subset \Dom(\psi_i)$ so that~$\mathcal U = \cup_i \psi_i(Q_i)$.
Then
\begin{equation}
  \norm{DF_L^\varphi}_{\BV(\mathcal U)}
    \defeq \sum_i 
      \norm{D (F_L^\varphi \circ \psi_i)}_{\BV(Q_i)}
\end{equation}
and
\begin{equation}
  \|DF_L^\varphi-DF_L^\iota\|_{L^1(\mathcal U)}
    \defeq \sum_i 
      \norm{
	D (F_L^\varphi \circ \psi_i)
	- D (F_L^\iota \circ \psi_i)
      }_{L^1(Q_i)}
      \,.
\end{equation}

\begin{proof}
We will use the representation~\eqref{26oct208} for $F_L(\cU)$.
As explained below~\eqref{e:resolvedCutKernel}, the function
$q(z,p,\nu)$ is $C^2$ in all its arguments if $\cM$ is a $C^4$-manifold.
In fixed local coordinates on~$\cU$, the function~$\bar q_1(p,\nu)$
defined in~\eqref{e:reciprocalCutTime} is the supremum of scalar
functions with uniformly bounded first and second derivatives.
Indeed, for $z$ close to~$p$, this follows from~\eqref{e:resolvedCutKernel}, and for $z$ bounded away from~$p$, it follows directly from~\eqref{26oct202}.
Thus, $\bar q_1$ is Lipschitz and semi-convex: adding a common quadratic makes each function, and hence their supremum, convex.
The same is true for~$\bar q_L(p,\nu)$, related to~$\bar q_1(p,\nu)$
by~\eqref{26oct210}.
Since~$\bar q_L\geq 1/L$, its reciprocal~$t_L$ is Lipschitz and locally uniformly semi-concave.
A Lipschitz semi-concave function has a distributional Hessian which is a finite matrix-valued measure on compact interior cubes.
We see then from~\eqref{26oct208} that~$F_L$ is Lipschitz and~$DF_L\in\mathrm{BV}$.
By~\eqref{e:resolvedCutKernel}, the $C^2$-bounds persist under sufficiently small~$C^4$ perturbations of the inclusion.
Consequently, the preceding estimates are uniform for all such perturbations, proving~\eqref{e:FLLipBV}.
\medskip

Now for $(p,\nu)\in\cU$, define~$\nu^\varphi(p,\nu)$
by~\eqref{26oct520} and note
\[
  F_L^\varphi(p,\nu)=\varphi(p)+t_L^\varphi\nu^\varphi(p,\nu)
  \,,\qquad
  F_L^\iota(p,\nu)=p+t_L^\iota\nu
  \,.
\]
Here~$t_L^\varphi$ denotes the cut time pulled back by the
identification of the unit normal bundle of~$\varphi(\mathcal M)$
and~$\mathcal U$.
Similarly, we define the function~$q^\vphi(z,p,\nu)$.
Now, equation~\eqref{e:resolvedCutKernel} applies to the function
\[
q^\vphi(z,p,\nu)
  =q(\vphi(z),\vphi(p),\nu^\vphi(p,\nu)),
\]
except that the chart mapping~$\psi$ in~\eqref{e:resolvedCutKernel}
needs to be replaced by the composition~$\vphi\circ\psi$.
It follows from the assumption~\eqref{e:C4distance} that for
each fixed coordinate cube~$Q$ we have
\begin{equation}\label{26oct521}
\|\vphi\circ\psi-\psi\|_{C^2(Q)}\le C\eta.
\end{equation}
Moreover, the normal fields~$\nu^\vphi$ and~$\nu$ are
$C^2$-close by~$O(\eta)$ in fixed charts.

We claim that
\[
  \sup_{\substack{z\in\cM\setminus\{p\}\\(p,\nu)\in\cU}}
  |q^\vphi(z,p,\nu)-q^\iota(z,p,\nu)|
  \le C\eta.
\]
To see this, cover~$\cM$ by finitely many coordinate charts.
For~$z$ sufficiently close to~$p$, write $z=\psi(u+s\omega)$ and~$p=\psi(u)$.
The formula~\eqref{e:resolvedCutKernel}, applied to
$\psi$ and~$\vphi\circ\psi$ with their respective normal fields,
shows that the corresponding kernels differ by~$O(\eta)$,
uniformly down to~$s=0$.
Indeed, their numerators and denominators differ by~$O(\eta)$,
and the denominators are uniformly bounded away from zero.
For~$z$ bounded away from~$p$, the same estimate follows
directly from~\eqref{26oct202}, since~$|z-p|$ is uniformly
bounded below. This proves the claim.

Taking suprema over~$z$ and then the maximum with~$1/L$
preserves this estimate. Since~$\bar q_L\geq 1/L$,
reciprocation is uniformly Lipschitz, and we obtain
\begin{equation}\label{e:cutTimeUniformStability}
  \|t_L^\varphi-t_L^\iota\|_{L^\infty(\cU)}
  \le C\eta
  \,.
\end{equation}
Consequently,
\begin{equation}\label{e:normalMapUniformStability}
  \|F_L^\varphi-F_L^\iota\|_{L^\infty(\cU)}
  \le C\eta
  \,.
\end{equation}
Thus, the first estimate in~\eqref{e:normalMapStability} holds.

For the derivative estimate in~\eqref{e:normalMapStability},
we follow a standard argument that can, for instance, be found in~\cite{EvansGariepy92}.
First choose enlarged cubes~$Q_i'$ so that
$Q_i\Subset Q_i'\Subset\Dom(\psi_i)$, and 
set
\[
  w_i =(F_L^\iota-F_L^\vphi)\circ\psi_i.
\]
By~\eqref{e:FLLipBV}, the derivative~$Dw_i$ belongs
to~$\BV(Q'_i)$, with
\[
  \|Dw_i\|_{\BV(Q'_i)}
  \le C.
\]
For sufficiently small~$h>0$, the difference quotient satisfies
\[
  \frac{w_i(x+he_j)-w_i(x)}{h}
  =\int_0^1\partial_j w_i(x+the_j)\,dt
\]
for almost every~$x\in Q_i$.
Using the translation estimate for BV functions, we obtain
\begin{equation}\label{e:wBVbound2}
\begin{aligned}
\int_{Q_i} \Big|
  \partial_j w_i(x)
  -\frac{w(x+he_j)-w(x)}{h}
\Big|\,dx
&\le
\int_0^1\int_Q
  |\partial_j w_i(x+the_j)-\partial_j w_i(x)|
\,dx\,dt\\
&\le
\int_0^1th\|Dw_i\|_{\BV(Q_i')}\,dt\\
&\le Ch.
\end{aligned}
\end{equation}
On the other hand, by~\eqref{e:normalMapUniformStability},
\begin{equation}\label{26oct523}
  \|w_i\|_{L^1(Q_i')}\le C\eta.
\end{equation}
Combining~\eqref{e:wBVbound2} and~\eqref{26oct523},
we obtain
\[
  \|\partial_j w_i\|_{L^1(Q_i)}
  \le Ch+\frac{C}{h}\eta.
\]
Taking~$h=\sqrt{\eta}$, with~$\eta_0$ sufficiently small,
and summing over the coordinate directions and the finitely
many charts gives the second bound in~\eqref{e:normalMapStability}.
\end{proof}

\subsection{The semialgebraic approximation of the normal bundle and cut-locus} 

We now prove Lemma~\ref{l:cutAlgebraicApproximation}.
The main idea is to use the Nash--Tognoli theorem to write~$\mathcal M$ as the smooth embedding of a real algebraic source manifold~$\mathcal S$, and then approximate the embedding by polynomials.
 
\begin{proof}[Proof of Lemma~\ref{l:cutAlgebraicApproximation}]
\textbf{Step 1. Construction of $\cB$.}
By the Nash--Tognoli theorem~\cite{Tognoli73}, there is a compact nonsingular real algebraic source manifold~$\cS\subseteq\R^a$, with some $a\ge d$, and a smooth embedding~$p:\cS\to\R^d$ such that
\begin{equation}
\cM=p(\cS)
\,.
\end{equation} 
Extend~$p$ smoothly to a box $\cQ_a\subset\R^a$ containing~$\cS$.
We now choose a sequence of polynomial maps~$p_N:\R^a\to\R^d$ of degree
at most~$aN$ such that, for every integer~$j>0$ there exists a constant~$C_j$ such that
\begin{equation}\label{e:embeddingPolynomialApproximation}
  \delta_N\defeq\|p-p_N\|_{C^4(\cS)}\leq \frac{C_j}{N^{j}},
  \quad\text{for all }
  N\ge 1
  \,.
\end{equation}
The existence of such a family is well known, and can be found in~\cite[Chapter~4]{Trefethen00}, for instance.
For convenience, we give a short account of this construction here.
Assume without loss of generality that $\cQ_a=[-1,1]^a$ and compose the
extension of~$p$ with the map 
\[
 \Phi(t)=(\cos t_1,\ldots,\cos t_a),~~t\in\R^a. 
\]
The resulting map $p_\Phi \defeq p\circ \Phi \colon \R^a \to \R^a$ is smooth, periodic, and even in each variable.
Truncating the Fourier series of $p_\Phi$ at frequency~$N$ in each variable gives
a trigonometric polynomial~$q_{N}(t)$ in~$\cos t_1,\ldots,\cos t_a$, of total degree at
most~$aN$. Using the Chebyshev polynomials 
\[
T_n(x) \defeq \cos(n\arccos(x)),
\]
we recover from the trigonometric polynomial
$q_{N}(t)$ a polynomial $p_N(x)$ defined on $Q_a$.   
The Fourier coefficients of~$p_\Phi(t)$ decay faster than any power of~$N$,
whereas the first four derivatives of the Chebyshev polynomials are
bounded on~$[-1,1]$ by
\begin{equation}
\|T_N\|_{C^4[-1,1]}\le CN^8.
\end{equation}
Summing the differentiated tails of the Fourier series 
therefore gives~\eqref{e:embeddingPolynomialApproximation}.

Since $p$ is an embedding of $\cS$, it follows from~\eqref{e:embeddingPolynomialApproximation} that there exists~$N_0 \in \N$ such that for all~$N \geq N_0$, the map~$p_N$ is also an embedding of $\cS$.
The images 
\begin{equation}
\cM_N=p_N(\cS)\subset\R^d,
\end{equation}
are compact semialgebraic manifolds because $\cS$ is an algebraic manifold and $p_N$ is a polynomial
map. 
Set
\begin{equation}\label{26oct1102}
  \mathcal B=\{(u,\omega)\in \cS\times\R^d:
       |\omega|=1,\quad
       \omega\perp Dp_{N_0}(u)\cT_u\cS\}.
\end{equation}
Thus,~$\mathcal B$ is the unit normal bundle of~$\cM_{N_0}=p_{N_0}(\cS)$,
with base points labelled by~$\cS$. Since~$\cB$ is the unit normal bundle of 
a compact semialgebraic manifold in $\R^d$, it is itself a compact smooth
semialgebraic manifold of dimension~$d-1$. 
\medskip

\textbf{Step 2: Construction of the maps~$F$ and~$F_N$.}
Let~$P_p^\perp(u)$ and $P_{p_N}^\perp(u)$ be the orthogonal
projections onto the normal space to~$\cM$ and~$\cM_N$, respectively, 
at the corresponding points~$p(u)\in\cM$ and $p_N(u)\in\cM_N$, and set
\begin{equation}\label{26oct806}
  \nu(u,\omega)
     =\frac{P_p^\perp(u)\omega}{|P_p^\perp(u)\omega|},~~
     \nu_{N}(u,\omega)
     =\frac{P_{p_N}^\perp(u)\omega}{|P_{p_N}^\perp(u)\omega|},~~(u,\omega)\in\cB.
\end{equation}
We can choose~$N_0$ sufficiently large so that
the denominators in~\eqref{26oct806} do not vanish.

Let~$t_L(u,\omega)$ and $t_{L,N}(u,\omega)$ be the respective
normal cut time of the manifolds~$\cM$ and~$\cM_N$ at
the points $(p(u),\nu(u,\omega))$ and $(p_N(u),\nu_N(u,\omega)$, truncated at~$L$, and define
the maps
\begin{equation}\label{26oct710}
\begin{aligned}
  &F(u,\omega)=p(u)+t_{L}(u,\omega)\nu(u,\omega),\\
  &F_{N}(u,\omega)=p_N(u)+t_{L,N}(u,\omega)\nu_{N}(u,\omega),~~(u,\omega)\in\cB.
  \end{aligned}
\end{equation}
In particular, each cut time in~\eqref{26oct710} is computed on its own manifold; only
the truncation value~$L$ is shared.
\medskip

 \textbf{Step 3. Uniform bounds on $F$ and $F_N$.}
We now prove the uniform bounds~\eqref{26oct702} and~\eqref{26oct704}. 
Choose a fixed finite smooth semialgebraic atlas~$\psi_i$ on $\mathcal B$ and finitely many cubes
\[
  Q_i\Subset Q_i'\Subset Q_i'', ~~i=1,\dots,K,
\]
compactly contained in the chart domains, with the union~$\bigcup_{i=1}^K\psi_i(Q_i)$ covering~$\mathcal B$.

Define the maps~$G, G_N$ by
\begin{equation}\label{26oct714}
  G(u,\omega)=(p(u),\nu(u,\omega)),~~
  G_{N}(u,\omega)=(p_N(u),\nu_{N}(u,\omega)),~~(u,\omega)\in\cB, 
  \,,
\end{equation}
and note they are diffeomorphisms between~$\mathcal B$ and the corresponding unit normal bundles of $\cM=p(\cS)$ and~$\cM_N=p_N(\cS)$, respectively. 
Let~$I_N$ be the identification between~$\mathcal U$ and the unit normal
bundle of~$p_N(\cS)$ via~\eqref{26oct520}.
Then, we can represent $F$ and $F_N$ as 
\[
  F_N=F_L^{\varphi_N}\circ(I_N^{-1}\circ G_{N}),
  \qquad F=F_L\circ G.
\]
Here, we use the notation $F^{\varphi}$ as in Lemma~\ref{l:truncatedNormalRegularity}.
Note this representation proves the first assertion of Lemma~\ref{l:cutAlgebraicApproximation}.

Now consider the embeddings $\varphi_N \colon \mathcal M \to \R^d$ defined by~$\varphi_N \defeq p_N\circ p^{-1}$.
By~\eqref{e:embeddingPolynomialApproximation} we note
\begin{equation}\label{26oct716}
\|\varphi_N-\iota\|_{C^4}\leq C\delta_N
\,,
\end{equation}
and by making~$N_0$ larger if necessary, we can ensure Lemma~\ref{l:truncatedNormalRegularity} applies for every~$N \geq N_0$.
Since the coordinate changes~$I_N^{-1}\circ G_{N}$ and~$G$, and their inverses, have~$C^2$ norms which are uniformly bounded in~$N$, inequality~\eqref{e:FLLipBV} in Lemma~\ref{l:truncatedNormalRegularity} implies
\[
  \|F\|_{\Lip(Q_i'')}+\|F_{N}\|_{\Lip(Q_i'')}
  +\|DF\|_{\BV(Q_i'')}+\|DF_{N}\|_{\BV(Q_i'')}
  \leq C.
\]
This proves the bounds~\eqref{26oct702} and~\eqref{26oct704} on larger cubes.
\medskip

\textbf{Step 4. Proof of approximation bounds~\eqref{e:cutApproximationError}.} 
We next compare~$F_N$ and~$F$ at the same point
$b=(u,\omega)\in\mathcal B$. The normal fields satisfy
\begin{equation}
\|\nu_N-\nu\|_{C^2(\cB)}\leq C\delta_N.
\end{equation}
For~$v\in \cS\setminus\{u\}$ and $b\in\cB$, define
\[
q(v,b)=
  \frac{2\nu(b)\cdot(p(v)-p(u))}{|p(v)-p(u)|^2},~~
  q_N(v,b)=
  \frac{2\nu_N(b)\cdot(p_N(v)-p_N(u))}{|p_N(v)-p_N(u)|^2}.
\]
We claim that
\begin{equation}\label{26oct718}
  \sup_{b=(u,\omega)\in\mathcal B}\ 
  \sup_{v\ne u}|q_N(v,b)-q(v,b)|\leq C\delta_N.
\end{equation}
This follows from the analysis in the proof of Lemma~\ref{l:truncatedNormalRegularity}
that leads to~\eqref{e:resolvedCutKernel} and the error bound~\eqref{26oct716}.

Going back to~\eqref{26oct202}-\eqref{e:truncatedReciprocal}, we have 
\begin{equation}\label{26oct720}
\frac{1}{t_L(b)}=\max\Big\{\frac1L,\sup_{v\ne u}q(v,b)\Big\},~~
\frac{1}{t_{L,N}(b)}=\max\Big\{\frac1L,\sup_{v\ne u}q_N(v,b)\Big\}.
\end{equation}
We conclude from~\eqref{26oct710},~\eqref{26oct718} and~\eqref{26oct720} 
that
\[
  \|t_{L,N}-t_L\|_{L^\infty(\mathcal B)}\leq C\delta_N,
  \qquad
  \|F_N-F\|_{L^\infty(\mathcal B)}\leq C\delta_N.
\]
Recalling~\eqref{e:embeddingPolynomialApproximation}, we see that this is the first estimate 
in~\eqref{e:cutApproximationError}. The second estimate in~\eqref{e:cutApproximationError}
follows exactly as in the proof of the BV bound for the derivative bound 
in~\eqref{e:normalMapStability}, at the end of the proof of Lemma~\ref{l:truncatedNormalRegularity}.  
\end{proof}

\subsection{A polynomial description of the normal bundle}

In order to prove Lemma~\ref{l:endpointPolynomial}, we need a preliminary preliminary result showing the normal bundle can be described as the level set of a polynomial of~$(d+1)$  variables. 

\begin{lemma}\label{l:criticalDistancePolynomial-bis}
There is a nonzero polynomial~$A_N(x,\rho)$, $x\in\R^d$, $\rho\in\R$, such that
\begin{equation}\label{e:criticalDistanceRelation}
  A_N\bigl(x,|x-p_N(u)|^2\bigr)=0,
  \quad\hbox{for all } 
  u\in\cS,\text{ and $x\in\R^d$ s.t.}~~x-p_N(u)\perp Dp_N(u)\cT_u\cS,
\end{equation}
and with
\begin{equation}\label{e:criticalDistanceDegree}
  \deg A_N\leq CN^{P_0},
\end{equation}
with an exponent~$P_0$ that does not depend on $N$.
\end{lemma}

\begin{proof}
{\bf Step 1. A semialgebraic description of the normal.}
Let us first explain how the normal $n=\nu_N(u,\omega)$, for $(u,\omega)\in\cB$ 
can be described by a fixed number of polynomial equations and inequalities of degree 
at most $CN$ in a fixed number of variables that is independent of $N$.

Choose a fixed family of polynomial equations~$R_j(u)=0$,~~$1\le j\le \ell$, 
to represent the algebraic  manifold~$\cS$:
\begin{equation}\label{26oct1104}
\cS=\{u\in\R^a:~R_j(u)=0,~\text{for all $1\le j\le \ell$}\}.
\end{equation}
In particular, the vector $R(u)=(R_1,\dots,R_\ell(u))$ satisfies 
\begin{equation}\label{26oct722}
\cT_u\cS=\ker DR(u),~~\text{ for~$u\in \cS$.}
\end{equation}
We claim that, given $(u,\omega)\in\cB$, the identity~$n=\nu_{N}(u,\omega)$ is equivalent to the existence
of~$\xi\in\R^a$,~$\zeta\in\R^\ell$, and~$\lambda\in\R$ such that 
\begin{equation}\label{26oct721}
\begin{aligned}
  &DR(u)\xi=0,\qquad
  \omega-\lambda n=Dp_N(u)\xi,~~\lambda>0,
\end{aligned}
\end{equation}
and
\begin{equation}\label{26oct723}
[Dp_N(u)]^Tn=[DR(u)]^T\zeta,\qquad
  |n|^2=1.
\end{equation}
Indeed, due to~\eqref{26oct722}, the first condition in~\eqref{26oct721} says 
that $\xi\in\R^a$ is a tangent vector to~$\cS$
at the point $u\in\cS$. Hence, the second 
equation in~\eqref{26oct721} says that~$\omega-\lambda n$
is tangent to~$\cM_N=p_N(\cS)$. The first equation in~\eqref{26oct723} says  that~$n$ is normal
to~$p_N(\cS)$. This is because any tangent vector $v$ to $\cM_N=p_N(\cS)$ at the point $p_N(u)$ has the form $v=Dp_N(u)\eta$, with some
$\eta\in\cT_u\cS$. We know from~\eqref{26oct722} that $DR(u)\eta=0$.
Therefore, if $n$ satisfies the first condition in~\eqref{26oct723}, then
\begin{equation}
(n\cdot v)=(n\cdot Dp_N(u)\eta)=([Dp_N(u)]^Tn\cdot\eta)=([DR(u)]^T\zeta\cdot\eta)=(\zeta\cdot[DR(u)\eta])
=0.
\end{equation}
Thus,~$\lambda n$ is the orthogonal projection
of~$\omega$ onto the normal space to $\cM_N$, and the last condition in~\eqref{26oct723}
normalizes the vector $n$. 
The system~\eqref{26oct721}-\eqref{26oct723}, together with~\eqref{26oct1104},
consists of a fixed number of polynomial equations 
and inequalities,
in a fixed number of scalar variables that is independent of $N$. As the polynomials $p_N(u)$ have degrees
at most~$aN$, the polynomials forming this system have degrees at most~$c_aN$.

{\bf Step 2. Getting rid of the auxiliary variables.} %
The next step is to get rid of the auxiliary variables $\lambda>0$, $\zeta\in\R^\ell$ and~$\xi\in\R^a$  
in~\eqref{26oct721}-\eqref{26oct723}.
Let us set~$m=\dim\cS$ and consider a minor of the $\ell\times a$ matrix~$DR$ of the size $(a-m)\times(a-m)$.
On any coordinate patch where this minor is not zero, select the
corresponding~$a-m$ equations
\begin{equation}\label{26oct1202}
R_{j_1}(u)=\ldots=R_{j_{a-m}}(u)=0
\end{equation}
in~\eqref{26oct1104}.
The zero set of~\eqref{26oct1202} agrees locally with the manifold~$\cS$ because of~\eqref{26oct1104}
and non-degeneracy of the minor of the size $(a-m)$. 
On one such patch $\cQ$, and with the above choice of $R_{j_k}$, set
\[
  h(u,x)=|x-p_N(u)|^2,~u\in\cS,~x\in\R^d,
\]
and consider the following   system of~$a$ polynomial
equations in $a$ variables~$(u_1,\dots,u_a)$, with~$x\in\R^d$ fixed
\begin{equation}\label{e:stationaryDistanceSystem}
\begin{aligned}
  f_i(u,x)&=R_{j_i}(u)=0,
       &&1\leq i\leq a-m,\\
  f_{a-m+j}(u,x)
       &=(x-p_N(u))\cdot [Dp_N(u)V_j(u)]=0,
       &&1\leq j\leq m.
\end{aligned}
\end{equation}
Here, we make use of
the polynomial vector fields~$V_1(u),\ldots,V_m(u)$ that form a
basis of the tangent space~$\cT_u\cS$ on $\cQ$. The first set of equations in~\eqref{e:stationaryDistanceSystem}
says that~$u\in\cS$, while the second set  
says that $x-p_N(u)$ is orthogonal to the manifold $\cM_N$. 

{\bf Step 3. Reduction to a single polynomial.} We now show that solutions to the system~\eqref{e:stationaryDistanceSystem} lie in the zero set
of a single polynomial in $(d+1)$-variables. 
Choose an integer~$D\leq CN$ strictly larger than the degrees
in~$u$ of both all polynomials in~\eqref{e:stationaryDistanceSystem} and of~$h$. This is possible because the degree of $p_N(u)$
is not larger than $CN$, while the degrees in~$x$
are at most one for~$f_i$, and two for~$h$.

Let us take~$\eps>0$ small and perturb  \eqref{e:stationaryDistanceSystem} to the system 
\begin{equation}\label{e:perturbedDistanceSystem}
  f_i(u,x)+\varepsilon u_i^D=0,
       \qquad i=1,\ldots,a.
\end{equation}
The highest order power in $u_i$ in each of these equations is the monomial $u_i^D$. 

Consider the set $\cP$ of all the monomials of the form $u_1^{\alpha_1}\dots u_a^{\alpha_a}$ with $0\leq\alpha_i<D$. 
Their total number is~$K_D=D^a$.
For each of the $K_D$ monomials in $\cP$, repeatedly 
reduce the polynomial
\begin{equation}\label{26oct1108}
p_\alpha(x,u)=h(u,x)u^\alpha,
\end{equation}
using the rule
\begin{equation}\label{26oct912}
  u_i^D=-\varepsilon^{-1}f_i(u,x),
\end{equation}
motivated by~\eqref{e:perturbedDistanceSystem}. 
The total degree in~$u$ of each polynomial strictly decreases after each such reduction. After repeatedly
applying the reduction, starting with~$p_\alpha(x,u)$, we will eventually
obtain a linear combination
of the monomials in $\cP$, 
with coefficients depending on~$x$ and $\eps$. 
This gives a~$K_D\times K_D$ matrix~$T(x,\varepsilon)$ of the coefficients
in these linear combinations, obtained starting with each monomial in $\cP$.  
We claim that 
\begin{equation}\label{26oct1110}
\text{the entries  $T_{pj}(x,\eps)$, $1\le p,j\le K_D$, 
are polynomials in $x$ and $\eps^{-1}$ of degree~$M\le CD$.}
\end{equation}
Indeed,  before the reduction, the total degree in~$u$ of the 
polynomial $p_\alpha(x,u)$ in~\eqref{26oct1108} is at most
\[
D-1+a(D-1)=(a+1)(D-1).
\]
Thus, along each chain of replacements there are at most
$(a+1)D$ steps. Each step adds at most one to the degree in~$x$
and to the power of~$\varepsilon^{-1}$.  This proves~\eqref{26oct1110}.  

Given some $u\in\cS$, let $e(u)$ be the row of length $K_D$ containing these monomials. 
By the above construction of the matrix $T(x,\eps)$, every solution $u_\eps$ to~\eqref{e:perturbedDistanceSystem}
satisfies the identity
\[
  e(u_\eps)T(x,\varepsilon)=h(u_\eps,x)e(u_\eps).
\]
It follows that
\begin{equation}\label{26oct914}
  \det\bigl(\rho_\eps I-T(x,\varepsilon)\bigr)=0
       \quad\hbox{at }\rho_\eps=h(u_\eps,x).
\end{equation}
Since the matrix~$T(x,\eps)$ is of the size~$K_D\times K_D$, 
$K_D=D^a$, and because of~\eqref{26oct1110}, the polynomial
\[
  \mathcal R(x,\rho,\varepsilon)
    =\det\bigl(\varepsilon^M\rho I-
                         \varepsilon^M T(x,\varepsilon)\bigr)
\]
has the total degree of at most $CD^{a+1}$.  Here, $M$ is also as in~\eqref{26oct1110}. 
Let~$A(x,\rho)$
be the leading order nonzero coefficient in $\eps$ of~$\mathcal R(x,\rho,\varepsilon)$:
\begin{equation}\label{26oct1114} 
  \mathcal R(x,\rho,\varepsilon)
     =\varepsilon^q\bigl(A(x,\rho)
                             +\varepsilon A_1(x,\rho)+\cdots\bigr).
\end{equation}
The polynomial $A(x,\rho)$ is also of degree at most~$CD^{a+1}$. By~\eqref{26oct914}, we know that
\begin{equation}\label{26oct1112} 
\eps^{-q}\mathcal R(x,h(u_\eps,x),\varepsilon)=0,
\end{equation}
for all $\eps>0$ sufficiently small. 

Next, note that if $u\in\R^a$ is a solution to~\eqref{e:stationaryDistanceSystem} for a given $x\in\R^d$ 
such that the $a\times a$ matrix~$D_uf$ is non-singular, then there exists a sequence $\eps_n\to 0$ and solutions $u_{\eps_n}$ to~\eqref{e:perturbedDistanceSystem}
such that $u_{\eps_n}\to u$ as $n\to+\infty$. Passing to the limit $\eps_n\to 0$ in~\eqref{26oct1112}, and taking~\eqref{26oct1114} into account, we obtain
\begin{equation}\label{26oct1116}
  A\bigl(x,h(u,x)\bigr)=0
\end{equation}
at every such nonsingular solution.  

We claim that~\eqref{26oct1116} actually holds for any $x\in\R^d$ such that $x-p_N(u)$ is normal to~$\cM_N=p_N(\cS)$, and with $u$ in the patch under
consideration, without the non-singularity assumption on $D_uf$. 
To see this, set  
\[
w=x-p_N(u),~~x_\theta=p_N(u)+\theta w.
\]
For a given $x$, the system~\eqref{e:stationaryDistanceSystem} says that $u\in\cS$ and   $x-p_N(u)$ is normal
to $\cM_N=p_N(\cS)$ at the point $z=p_N(u)\in\cM_N$. In particular, at $\theta=0$, for $x_0=p_N(u)$, this system has an isolated solution $(u,p_N(u))$.
Therefore, the determinant of the gradient matrix $D_uf$ is not zero at $\theta=0$. As this determinant is a polynomial in $\theta$, it follows that 
it does not vanish identically in $\theta$. This, in turn, means that
it vanishes at most at a discrete set of $\theta$. We may now approximate $\theta=1$ by a sequence $\theta_k\to 1$ of non-singular $\theta_k$. For each of $\theta_k$, we have
\begin{equation}
  A\bigl(x_{\theta_k},h(u,x_{\theta_k})\bigr)=0.
\end{equation}
Passing to the limit $\theta_k\to 1$, we conclude that~\eqref{26oct1116} holds for all $u\in\cS$ and $x$ such that~$x-p_N(u)$ is normal
to $\cM_N=p_N(\cS)$ at the point $p_N(u)$.

Taking the product of the nonzero polynomials~$A$ from all
the patches gives~$A_N$ satisfying~\eqref{e:criticalDistanceRelation}--\eqref{e:criticalDistanceDegree}, finishing the proof of Lemma~\ref{l:criticalDistancePolynomial-bis}.
\end{proof}

\subsection{The semialgebraic approximation of the cut-locus (Lemma~\ref{l:endpointPolynomial})}\label{s:endpointPolynomial}

We now prove Lemma~\ref{l:endpointPolynomial} showing that the cut-locus of the manifolds~$\cM_N$ is contained in a real algebraic set of a moderately large degree.
Before delving into the proof, we study a few simple examples.
\begin{enumerate}
  \item If $\cS=\{x_2=\dots=x_d=0\}$ be a line in~$\R^d$, then its cut-locus with a cut-off time $L$ is the annulus $\{x_2^2+\dots+x_d^2=L^2\}$.

  \item If $\cS=\{x\in\R^d:~|x|^2=R^2\}$ is a sphere in $\R^d$.
    If the cutoff time $L<R$, the cut-locus is the union $\{|x|^2=(L+R)^2\}\cup\{|x|^2=(R-L)^2\}$.  If $L\ge R$ then the cut-locus is the union $\{|x|^2=(L+R)^2\}\cup\{x=0\}$.
\end{enumerate}
In either cases, it is easy to see that the cut-locus lies in a zero set of a polynomial $H(x)$. 

\begin{proof}[Proof of Lemma~\ref{l:endpointPolynomial}]
We start with the polynomial~$A_N$ obtained in Lemma~\ref{l:criticalDistancePolynomial-bis},
so that~\eqref{e:criticalDistanceRelation} holds. Our goal is to eliminate the variable $\rho$ so that
we would end up with a polynomial $H_N(x)$ of~$x\in\R^d$ alone, and such that $F_N(\cB)$ is 
contained in the 
zero level set~$\set{H_N = 0}$.
Let us remove the repeated irreducible factors in $A_N(x,\rho)$, and also remove
the factor~$\rho-L^2$ if $A_N(x,\rho)$ has one. We denote the resulting polynomial by 
$Q_N(x,\rho)$. This polynomial still satisfies~~\eqref{e:criticalDistanceRelation}
\begin{equation}\label{e:reducedDistanceRelation} 
Q_N\bigl(x,|x-p_N(u)|^2\bigr)=0,
  \quad\hbox{whenever}\quad
  u\in\cS,~~x-p_N(u)\perp Dp_N(u)\cT_u\cS,
\end{equation}
This is because, first, removing repeated factors preserves the zero
set, and, second,  removing the factor~$\rho-L^2$ also preserves
\eqref{e:criticalDistanceRelation} since the latter holds for all normal vectors.

Consider the polynomial 
\begin{equation}\label{e:endpointContainingPolynomial}
  H_N(x)=%
  Q_N(x,L^2)\,
       \operatorname{Res}_{\rho}
          \bigl(Q_N(x,\rho),\partial_\rho Q_N(x,\rho)\bigr).
\end{equation}
Here, $\operatorname{Res}_{\rho}\bigl(Q_N(x,\rho),\partial_\rho Q_N(x,\rho)\bigr)$ is the resultant of these two polynomials in $\rho$
and is itself a polynomial in $x$ that does not vanish, as long as $Q_N(x,\rho)$ 
and $\partial_\rho Q(x,\rho)$ do not have common roots
in $\rho$. Since all repeated factors have been removed, these common roots do not exist,
and the resultant is a non-zero polynomial in $x$. 
The factor $Q_N(x,L^2)$ is also not identically zero, because the factor
$\rho-L^2$ has been removed from $A_N(x,\rho)$. 
Thus, the polynomial~$H_N(x)$ is not identically equal to zero. Its
degree satisfies
\begin{equation}\label{e:endpointPolynomialDegree}
  \deg H_N\leq CN^{P_2}
\end{equation}
for a fixed exponent~$P_2$. This is because the Sylvester matrix used to compute 
the resultant has the size
at most~$2\deg Q_N$, and each entry in the matrix has degree at most
$\deg Q_N$. Thus, its determinant is a polynomial of degree at most $N^{P}$ with some $P$ that 
does not depend on $N$. 

We claim that
\begin{equation}\label{e:endpointPolynomialContainment}
  F_N(\cB)\subseteq\{H_N=0\}.
\end{equation}
We first make the following observation. If we set 
\begin{equation}\label{26oct1212}
\delta(x)=\dist(x,\cM_N)^2,~~x\in\R^d,
\end{equation}
then, we have 
\begin{equation}\label{26oct1118}
Q_N(x,\delta(x))=0,~~\hbox{ for every~$x\in\R^d$},
\end{equation}
because of  \eqref{e:reducedDistanceRelation}. 

Now suppose that
\[
x_0=p_N(u)+t _L(u,n)n\in F_N(\cB),
\]
where~$n$ is a unit normal to $\cM_N$ at the point $p_N(u)\in\cM_N$, and~$t_L(u,n)$ is its truncated cut time on~$\cM_N$.
If~$t=L$, then $\delta(x_0)=L^2$, and~$Q_N(x_0,L^2)=0$ by~\eqref{26oct1118}. Therefore, in that case we have 
$H_N(x_0)=0$.

If~$t_L(u,n)<L$, set
\[
x_\tau=p_N(u)+\tau n,~~\tau\ge 0.
\]
The vector $x_\tau-p_N(u)$ is normal to $\cM_N$ at the point $p_N(u)$, 
for every~$\tau\ge 0$, which, together 
with~\eqref{e:reducedDistanceRelation} implies that 
\[
Q_N(x_\tau,\tau^2)=0,~~\hbox{ for all $\tau\ge 0$.} 
\]
Thus, for~$\tau$ sufficiently close to~$t_L(u,n)$, both~$\rho_1=\tau^2$ and
$\rho_2=\delta(x_\tau)$ are close to~$t_L^2(u,n)$ and are roots of this
same polynomial $Q_N(x_\tau,\rho)$. On the other hand, if $H_N(x_0)\neq 0$, then 
\[
\operatorname{Res}_{\rho}Q_N(x_0,\rho),\partial_\rho Q_N(x_0,\rho)\neq 0.
\]
It follows that $\rho=\delta(x_0)$ is a simple root of $Q_N(x_0,\rho)$. As a consequence, for $\tau$ sufficiently close to
$t_L(u,n)$ the polynomial $Q(x_\tau,\rho)$ has a unique root close to $t_L^2(u,n)$. 
We conclude that 
\[
  \delta(x_\tau)=\tau^2,
\]
also for some~$\tau>t_L(u,n)$ sufficiently close to $t_L(u,n)$. 
Therefore,~$p_N(u)$ is the nearest point on $\cM_N$ even 
for some~$t>t_L(u,n)$ which is a contradiction since~$t_L(u,n)<L$. 
Therefore, we must have~$H_N(x_0)=0$, which 
proves~\eqref{e:endpointPolynomialContainment} and finishes the proof of Lemma~\ref{l:endpointPolynomial}.
 \end{proof}

\subsection{The parametrized area of the semialgebraic approximation (Lemma~\ref{l:algebraicCutArea})}\label{s:algebraicCutArea}

We now prove Lemma~\ref{l:algebraicCutArea} bounding the parametrized area of the semialgebraic approximation.

\begin{proof}[Proof of Lemma~\ref{l:algebraicCutArea}]
Let $Q_i$, $i=1,\dots,K$, 
be the charts for the manifold $\cB$, and 
$F_{N,i}$ is the restriction of~$F_N$ to $Q_i$. The cubes $Q_i'$ below contain the cubes $Q_i$. 

Let~$\pi_k:\R^d\to\R^{d-1}$ be the projection that omits the~$k$th coordinate and set 
\[
g_{k,i}=\pi_k\circ F_{N,i}:~Q_i'\to\Rm^{d-1}. %
\]
Using the Cauchy-Binet formula for the determinant of the form $\det(A^TA)$ in the case where~$A$ is a $d\times(d-1)$ matrix, we get
\begin{equation}\label{26oct1802}
  J_{d-1}F_{N,i}=\sqrt{\det(DF_{N,i}^T DF_{N,i})}
    =\Big(\sum_{k=1}^d|\det Dg_{k,i}|^2\Big)^{1/2}
    \leq\sum_{k=1}^d|\det Dg_{k,i}|.
\end{equation}
Fix some $x\in\Rm^d$ and $s>0$, and set
\[
V_i=Q_i'\cap F_{N,i}^{-1}(B(x,s)),~~1\le i\le K. 
\]
Applying 
the area formula to the Lipschitz map $g_{k,i}=\pi_k\circ F_{N,i}$ gives 
\begin{equation}\label{26oct1804} %
  \begin{aligned}
  \int_{V_i}|\det Dg_{k,i}(y)|\,dy
    & =\int_{\Rm^{d-1}}\#\{y\in V_i:g_{k,i}(y)=z\}\,dz
    \\
    &= \int_{\pi_k(B(x,s))}\#\{y\in V_i:g_{k,i}(y)=z\}\,dz,
  \end{aligned}
\end{equation}
for all $1\le i\le K$ and $1\le k\le d$. One consequence of this identity is that 
\begin{equation}
N_{ik}(z):=\#\{y\in V_i:g_{k,i}(y)=z\}<+\infty,~~\hbox{for all $1\le i\le K$ and $1\le k\le d$,}
\end{equation}
for a.e. $z\in\pi_k(B(x,s))$.
Combining~\eqref{26oct1802} and~\eqref{26oct1804}, we obtain  
\begin{equation}\label{26oct1812}
\begin{aligned}
\sum_{i=1}^K\int_{Q_i'\cap F_{N,i}^{-1}(B(x,s))}
       J_{d-1}F_{N,i}(y)dy&\le\sum_{i=1}^K\sum_{k=1}^d\int_{V_i}|\det Dg_{k,i}(y)|\,dy 
   \\
   & =\sum_{i=1}^K\sum_{k=1}^d\int_{\pi_k(B(x,s))}N_{ik}(z)\,dz.
\end{aligned}
\end{equation}
Fix $k$ and choose $z\in\pi_k(B(x,s))$ outside a set
of measure zero, so that $N_{ik}(z)<+\infty$ for every $i$
and the restriction of $H_N$ to $\pi_k^{-1}(z)$
is not identically zero.
The identity $g_{k,i}(y)=z$ for some $y\in V_i$
means that there exists
\[
z'\in\ell_z:=\{\zeta\in\Rm^d:~\pi_k(\zeta)=z\},
\]
such that 
\begin{equation}\label{26oct1808}
F_{N,i}(y)=z'.
\end{equation}
In particular, we deduce from~\eqref{26oct1808} that $H_N(z')=0$. Here, $H_N(z)$ is the polynomial constructed 
in Lemma~\ref{l:endpointPolynomial}. 
Since~$H_N$ is not identically~$0$, the restriction of~$H_N$ to the line $\ell_z$ is a nonzero  
polynomial of one real variable for almost every~$z\in\R^{d-1}$.
Hence, it has at most~$\text{deg}(H_N)\le CN^{P_2}$ roots, by~\eqref{e:endpointPolynomialDegree}. 
Thus, there are at most $CN^{P_2}$ points $z'\in\ell_z$ for which~\eqref{26oct1808} may have a solution $y\in V_i$. 

Furthermore, for a fixed~$x\in\R^d$, let~$\mathcal A(x)=\{p_N(u_1),\dots,p_N(u_J)\}$ be the collection of its nearest-points
on~$\cM_N=p_N(\cS)$. The $u_j$ are the common real zeros of the polynomials (in the variable $u\in\cS\subset\R^a$)
\[
  R_1(u),\ldots,R_\ell(u),\qquad
  |x-p_N(u)|^2-\delta(x),
\]
with $\delta(x)$ defined by~\eqref{26oct1212}. 
These polynomials have degrees at most~$CN$ in $u$. The bound on the number of isolated real zeros therefore
implies that if $\#\mathcal A(x)<+\infty$, then\begin{equation}\label{e:nearestPointMultiplicity-2}
  \#\mathcal A(x)\leq CN^a,
\end{equation}
with the constant $C$ that does not depend on $x$. 

Going back to a point
$z'\in F_N(\cB)\cap B(x,s)$ with $\pi_k(z')=z$,
we first identify the fiber $F_N^{-1}(z')$ with
the nearest-point set $\cA(z')$.
Set $t=\dist(z',\cM_N)>0$, and let $q_0$ be the
nearest point supplied by one endpoint parameter over $z'$.
Every other nearest point $q$ determines the unit normal
$n=(z'-q)/t$.
If $t=L$, its truncated cut time is $L$.
If $t<L$ and $q\ne q_0$, then, for every $h>0$,
\[
  |q+(t+h)n-q_0|^2-(t+h)^2
  =-\frac{h}{t}|q-q_0|^2<0.
\]
Thus its cut time is exactly $t$.
The point $q_0$ already has truncated cut time $t$
by its choice.
Using the diffeomorphism $G_N$ in~\eqref{26oct714},
we therefore obtain a bijection between
$F_N^{-1}(z')$ and $\cA(z')$.

Since the charts cover $\cB$ and $z'\in B(x,s)$,
the finiteness of $N_{ik}(z)$ for every $i$
implies that $F_N^{-1}(z')$ is finite.
Hence $\#\cA(z')<+\infty$.
Thus, we have
\begin{equation}\label{26oct1810}
\#\cA(z')\le CN^a.
\end{equation}
We conclude from the argument below~\eqref{26oct1808}
and~\eqref{26oct1810} that for almost every
$z\in\pi_k(B(x,s))$ we have
\begin{equation}
N_{ik}(z)\le \#\{z'\in\ell_z:~z'\in F_{N,i}(V_i)\}\cdot CN^a\le CN^{a+P_2}.
\end{equation}
Combining this estimate with~\eqref{26oct1812}
\begin{equation}\label{26oct1814}
\begin{aligned}
\sum_{i=1}^K\int_{Q_i'\cap F_{N,i}^{-1}(B(x,s))}
       J_{d-1}F_{N,i}(y)dy&\le CN^{a+P_2}s^{d-1}.
\end{aligned}
\end{equation}
This proves~\eqref{e:algebraicCutArea}. 

\end{proof}

\bibliographystyle{halpha-abbrv}
\bibliography{gautam-refs1,gautam-refs2,preprints,ayya-refs}

\newcommand{\etalchar}[1]{$^{#1}$}
\begin{thebibliography}{BDBDD24}
\expandafter\ifx\csname url\endcsname\relax
  \def\url#1{\texttt{#1}}\fi
\expandafter\ifx\csname doi\endcsname\relax
  \def\doi#1{\burlalt{doi:#1}{http://dx.doi.org/#1}}\fi
\expandafter\ifx\csname urlprefix\endcsname\relax\def\urlprefix{URL }\fi
\expandafter\ifx\csname href\endcsname\relax
  \def\href#1#2{#2}\fi
\expandafter\ifx\csname burlalt\endcsname\relax
  \def\burlalt#1#2{\href{#2}{#1}}\fi

\bibitem[And82]{Anderson82}
B.~D.~O. Anderson.
\newblock Reverse-time diffusion equation models.
\newblock {\em Stochastic Process. Appl.}, 12(3):313--326, 1982.
\newblock \doi{10.1016/0304-4149(82)90051-5}.

\bibitem[BDBDD24]{BentonDeBortoliEA24}
J.~Benton, V.~De~Bortoli, A.~Doucet, and G.~Deligiannidis.
\newblock Nearly $d$-linear convergence bounds for diffusion models via
  stochastic localization.
\newblock In B.~Kim, Y.~Yue, S.~Chaudhuri, K.~Fragkiadaki, M.~Khan, and Y.~Sun,
  editors, {\em International Conference on Learning Representations}, volume
  2024, pages 36916--36936, 2024.
\newblock
  \urlprefix\url{https://proceedings.iclr.cc/paper_files/paper/2024/file/9edce07129d2ec4c511a8783eb3d81b5-Paper-Conference.pdf}.

\bibitem[BH70]{BrambleHilbert70}
J.~H. Bramble and S.~R. Hilbert.
\newblock Estimation of linear functionals on {S}obolev spaces with application
  to {F}ourier transforms and spline interpolation.
\newblock {\em SIAM J. Numer. Anal.}, 7:112--124, 1970.
\newblock \doi{10.1137/0707006}.

\bibitem[BMR22]{BlockMrouehEA22}
A.~Block, Y.~Mroueh, and A.~Rakhlin.
\newblock Generative modeling with denoising auto-encoders and {L}angevin
  sampling, 2022, \burlalt{2002.00107}{http://arxiv.org/abs/2002.00107}.
\newblock \urlprefix\url{https://arxiv.org/abs/2002.00107}.

\bibitem[BPR06]{BasuPollackEA06}
S.~Basu, R.~Pollack, and M.-F.~c. Roy.
\newblock {\em Algorithms in real algebraic geometry}, volume~10 of {\em
  Algorithms and Computation in Mathematics}.
\newblock Springer-Verlag, Berlin, second edition, 2006.

\bibitem[Buc77]{Buchner77}
M.~A. Buchner.
\newblock Simplicial structure of the real analytic cut locus.
\newblock {\em Proc. Amer. Math. Soc.}, 64(1):118--121, 1977.
\newblock \doi{10.2307/2040994}.

\bibitem[Buc78]{Buchner78}
M.~A. Buchner.
\newblock The structure of the cut locus in dimension less than or equal to
  six.
\newblock {\em Compositio Math.}, 37(1):103--119, 1978.
\newblock \urlprefix\url{http://www.numdam.org/item?id=CM_1978__37_1_103_0}.

\bibitem[CCL{\etalchar{+}}23]{ChenChewiEA23}
S.~Chen, S.~Chewi, J.~Li, Y.~Li, A.~Salim, and A.~Zhang.
\newblock Sampling is as easy as learning the score: theory for diffusion
  models with minimal data assumptions.
\newblock In {\em The Eleventh International Conference on Learning
  Representations}, 2023.
\newblock \urlprefix\url{https://openreview.net/forum?id=zyLVMgsZ0U_}.

\bibitem[CH04]{7080175}
J.~A. Costa and A.~O. Hero.
\newblock Learning intrinsic dimension and intrinsic entropy of
  high-dimensional datasets.
\newblock In {\em 2004 12th European Signal Processing Conference}, pages
  369--372, 2004.

\bibitem[CS04]{CannarsaSinestrari04}
P.~Cannarsa and C.~Sinestrari.
\newblock {\em Semiconcave functions, {H}amilton-{J}acobi equations, and
  optimal control}, volume~58 of {\em Progress in Nonlinear Differential
  Equations and their Applications}.
\newblock Birkh\"auser Boston, Inc., Boston, MA, 2004.

\bibitem[DB22]{DeBortoli22}
V.~De~Bortoli.
\newblock Convergence of denoising diffusion models under the manifold
  hypothesis.
\newblock {\em Transactions on Machine Learning Research}, 2022.
\newblock \urlprefix\url{https://openreview.net/forum?id=MhK5aXo3gB}.
\newblock Expert Certification.

\bibitem[DMWC26]{dennehy2026diffusionmodelsrecoveraccurate}
A.~Dennehy, R.~Muthukumar, R.~Willett, and N.~Chandramoorthy.
\newblock Diffusion models recover accurate mixture weights despite score
  function insensitivity, 2026,
  \burlalt{2607.15485}{http://arxiv.org/abs/2607.15485}.
\newblock \urlprefix\url{https://arxiv.org/abs/2607.15485}.

\bibitem[Efr11]{Efron11}
B.~Efron.
\newblock Tweedie's formula and selection bias.
\newblock {\em Journal of the American Statistical Association},
  106(496):1602--1614, 2011.
\newblock \doi{10.1198/jasa.2011.tm11181}.

\bibitem[EG92]{EvansGariepy92}
L.~C. Evans and R.~F. Gariepy.
\newblock {\em Measure theory and fine properties of functions}.
\newblock Studies in Advanced Mathematics. CRC Press, Boca Raton, FL, 1992.

\bibitem[FMN13]{fefferman2013testingmanifoldhypothesis}
C.~Fefferman, S.~Mitter, and H.~Narayanan.
\newblock Testing the manifold hypothesis, 2013,
  \burlalt{1310.0425}{http://arxiv.org/abs/1310.0425}.
\newblock \urlprefix\url{https://arxiv.org/abs/1310.0425}.

\bibitem[GS78]{GluckSinger78}
H.~Gluck and D.~Singer.
\newblock Scattering of geodesic fields. {I}.
\newblock {\em Ann. of Math. (2)}, 108(2):347--372, 1978.
\newblock \doi{10.2307/1971170}.

\bibitem[GS79]{GluckSinger79}
H.~Gluck and D.~A. Singer.
\newblock Scattering of geodesic fields. {II}.
\newblock {\em Ann. of Math. (2)}, 110(2):205--225, 1979.
\newblock \doi{10.2307/1971259}.

\bibitem[GT18]{Gorban_2018}
A.~N. Gorban and I.~Y. Tyukin.
\newblock Blessing of dimensionality: mathematical foundations of the
  statistical physics of data.
\newblock {\em Philosophical Transactions of the Royal Society A: Mathematical,
  Physical and Engineering Sciences}, 376(2118):20170237, Mar. 2018.
\newblock \doi{10.1098/rsta.2017.0237}.

\bibitem[HAHM25]{hakemi2025deeperdiffusionmodelsamplify}
S.~Hakemi, N.~Akhtar, G.~M. Hassan, and A.~Mian.
\newblock Deeper diffusion models amplify bias, 2025,
  \burlalt{2505.17560}{http://arxiv.org/abs/2505.17560}.
\newblock \urlprefix\url{https://arxiv.org/abs/2505.17560}.

\bibitem[Hir94]{Hirsch94}
M.~W. Hirsch.
\newblock {\em Differential topology}, volume~33 of {\em Graduate Texts in
  Mathematics}.
\newblock Springer-Verlag, New York, 1994.
\newblock Corrected reprint of the 1976 original.

\bibitem[HIS26]{HanIyerEA26}
R.~Han, G.~Iyer, and D.~Slepčev.
\newblock Time-complexity of sampling from a multimodal distribution using
  sequential monte carlo, 2026,
  \burlalt{2508.02763}{http://arxiv.org/abs/2508.02763}.
\newblock \urlprefix\url{https://arxiv.org/abs/2508.02763}.

\bibitem[HJA20]{HoJainEA20}
J.~Ho, A.~Jain, and P.~Abbeel.
\newblock Denoising diffusion probabilistic models.
\newblock In H.~Larochelle, M.~Ranzato, R.~Hadsell, M.~Balcan, and H.~Lin,
  editors, {\em Advances in Neural Information Processing Systems}, volume~33,
  pages 6840--6851. Curran Associates, Inc., 2020.
\newblock
  \urlprefix\url{https://proceedings.neurips.cc/paper_files/paper/2020/file/4c5bcfec8584af0d967f1ab10179ca4b-Paper.pdf}.

\bibitem[Hyv05]{Hyvarinen05}
A.~Hyv{{\"a}}rinen.
\newblock Estimation of non-normalized statistical models by score matching.
\newblock {\em Journal of Machine Learning Research}, 6(24):695--709, 2005.
\newblock \urlprefix\url{http://jmlr.org/papers/v6/hyvarinen05a.html}.

\bibitem[IT01]{ItohTanaka01}
J.-i. Itoh and M.~Tanaka.
\newblock The {L}ipschitz continuity of the distance function to the cut locus.
\newblock {\em Trans. Amer. Math. Soc.}, 353(1):21--40, 2001.
\newblock \doi{10.1090/S0002-9947-00-02564-2}.

\bibitem[KAAL22]{karras2022elucidatingdesignspacediffusionbased}
T.~Karras, M.~Aittala, T.~Aila, and S.~Laine.
\newblock Elucidating the design space of diffusion-based generative models,
  2022, \burlalt{2206.00364}{http://arxiv.org/abs/2206.00364}.
\newblock \urlprefix\url{https://arxiv.org/abs/2206.00364}.

\bibitem[KHR22]{koehler2022statisticalefficiencyscorematching}
F.~Koehler, A.~Heckett, and A.~Risteski.
\newblock Statistical efficiency of score matching: The view from isoperimetry,
  2022, \burlalt{2210.00726}{http://arxiv.org/abs/2210.00726}.
\newblock \urlprefix\url{https://arxiv.org/abs/2210.00726}.

\bibitem[LB04]{NIPS2004_74934548}
E.~Levina and P.~Bickel.
\newblock Maximum likelihood estimation of intrinsic dimension.
\newblock In L.~Saul, Y.~Weiss, and L.~Bottou, editors, {\em Advances in Neural
  Information Processing Systems}, volume~17. MIT Press, 2004.
\newblock
  \urlprefix\url{https://proceedings.neurips.cc/paper_files/paper/2004/file/74934548253bcab8490ebd74afed7031-Paper.pdf}.

\bibitem[Liu08]{Liu08}
J.~S. Liu.
\newblock {\em {M}onte {C}arlo strategies in scientific computing}.
\newblock Springer Series in Statistics. Springer, New York, 2008.

\bibitem[LLT22]{LeeLuEA22}
H.~Lee, J.~Lu, and Y.~Tan.
\newblock Convergence for score-based generative modeling with polynomial
  complexity.
\newblock In S.~Koyejo, S.~Mohamed, A.~Agarwal, D.~Belgrave, K.~Cho, and A.~Oh,
  editors, {\em Advances in Neural Information Processing Systems}, volume~35,
  pages 22870--22882. Curran Associates, Inc., 2022.
\newblock \doi{10.52202/068431-1662}.

\bibitem[LLZ{\etalchar{+}}23]{li2023autodiffusiontrainingfreeoptimizationtime}
L.~Li, H.~Li, X.~Zheng, J.~Wu, X.~Xiao, R.~Wang, M.~Zheng, X.~Pan, F.~Chao, and
  R.~Ji.
\newblock Autodiffusion: Training-free optimization of time steps and
  architectures for automated diffusion model acceleration, 2023,
  \burlalt{2309.10438}{http://arxiv.org/abs/2309.10438}.
\newblock \urlprefix\url{https://arxiv.org/abs/2309.10438}.

\bibitem[LN05]{LiNirenberg05}
Y.~Li and L.~Nirenberg.
\newblock The distance function to the boundary, {F}insler geometry, and the
  singular set of viscosity solutions of some {H}amilton-{J}acobi equations.
\newblock {\em Comm. Pure Appl. Math.}, 58(1):85--146, 2005.
\newblock \doi{10.1002/cpa.20051}.

\bibitem[LRS10]{LelievreRoussetEA10}
T.~Leli\`evre, M.~Rousset, and G.~Stoltz.
\newblock {\em Free energy computations}.
\newblock Imperial College Press, London, 2010.
\newblock \doi{10.1142/9781848162488}.
\newblock A mathematical perspective.

\bibitem[LSG24]{LeeSantanaGijzen24}
H.~Lee and M.~Santana-Gijzen.
\newblock Convergence bounds for sequential {M}onte {C}arlo on multimodal
  distributions using soft decomposition, 2024,
  \burlalt{2405.19553}{http://arxiv.org/abs/2405.19553}.
\newblock \urlprefix\url{https://arxiv.org/abs/2405.19553}.

\bibitem[LSHH26]{LiShenEA26}
X.~Li, Z.~Shen, Y.-P. Hsieh, and N.~He.
\newblock When scores learn geometry: Rate separations under the manifold
  hypothesis, 2026, \burlalt{2509.24912}{http://arxiv.org/abs/2509.24912}.
\newblock \urlprefix\url{https://arxiv.org/abs/2509.24912}.

\bibitem[Lu26]{Lu26}
J.~Lu.
\newblock A mathematical introduction to diffusion models, 2026,
  \burlalt{2607.01693}{http://arxiv.org/abs/2607.01693}.
\newblock \urlprefix\url{https://arxiv.org/abs/2607.01693}.

\bibitem[LWB24]{LuWangEA24}
Y.~Lu, Z.~Wang, and G.~Bal.
\newblock Mathematical analysis of singularities in the diffusion model under
  the submanifold assumption, 2024,
  \burlalt{2301.07882}{http://arxiv.org/abs/2301.07882}.
\newblock \urlprefix\url{https://arxiv.org/abs/2301.07882}.

\bibitem[LY86]{LiYau86}
P.~Li and S.-T. Yau.
\newblock On the parabolic kernel of the {S}chr\"{o}dinger operator.
\newblock {\em Acta Math.}, 156(3-4):153--201, 1986.
\newblock \doi{10.1007/BF02399203}.

\bibitem[LZB{\etalchar{+}}22]{lu2022dpmsolverfastodesolver}
C.~Lu, Y.~Zhou, F.~Bao, J.~Chen, C.~Li, and J.~Zhu.
\newblock Dpm-solver: A fast ode solver for diffusion probabilistic model
  sampling in around 10 steps, 2022,
  \burlalt{2206.00927}{http://arxiv.org/abs/2206.00927}.
\newblock \urlprefix\url{https://arxiv.org/abs/2206.00927}.

\bibitem[Mat83]{Mather83}
J.~N. Mather.
\newblock Distance from a submanifold in {E}uclidean space.
\newblock In {\em Singularities, {P}art 2 ({A}rcata, {C}alif., 1981)},
  volume~40 of {\em Proc. Sympos. Pure Math.}, pages 199--216. Amer. Math.
  Soc., Providence, RI, 1983.
\newblock \doi{10.1090/pspum/040.2/713249}.

\bibitem[Nea01]{Neal01}
R.~M. Neal.
\newblock Annealed importance sampling.
\newblock {\em Stat. Comput.}, 11(2):125--139, 2001.
\newblock \doi{10.1023/A:1008923215028}.

\bibitem[Pid22]{Pidstrigach22}
J.~Pidstrigach.
\newblock Score-based generative models detect manifolds.
\newblock In {\em Proceedings of the 36th International Conference on Neural
  Information Processing Systems}, NIPS '22, Red Hook, NY, USA, 2022. Curran
  Associates Inc.

\bibitem[PJT19]{PaulinJasraEA19}
D.~Paulin, A.~Jasra, and A.~Thiery.
\newblock Error bounds for sequential {M}onte {C}arlo samplers for multimodal
  distributions.
\newblock {\em Bernoulli}, 25(1):310--340, 2019.
\newblock \doi{10.3150/17-bej988}.

\bibitem[RIG{\etalchar{+}}25]{roos2025imetbiasinvestigating}
N.~Roos, E.~Iakovleva, A.~Gjergji, V.~P. Pastore, and E.~Tartaglione.
\newblock How i met your bias: Investigating bias amplification in diffusion
  models, 2025, \burlalt{2512.20233}{http://arxiv.org/abs/2512.20233}.
\newblock \urlprefix\url{https://arxiv.org/abs/2512.20233}.

\bibitem[SE19]{SongErmon19}
Y.~Song and S.~Ermon.
\newblock Generative modeling by estimating gradients of the data distribution.
\newblock In H.~Wallach, H.~Larochelle, A.~Beygelzimer, F.~d\textquotesingle
  Alch\'{e}-Buc, E.~Fox, and R.~Garnett, editors, {\em Advances in Neural
  Information Processing Systems}, volume~32. Curran Associates, Inc., 2019.
\newblock
  \urlprefix\url{https://proceedings.neurips.cc/paper_files/paper/2019/file/3001ef257407d5a371a96dcd947c7d93-Paper.pdf}.

\bibitem[SFK24]{pmlr-v235-sabour24a}
A.~Sabour, S.~Fidler, and K.~Kreis.
\newblock Align your steps: Optimizing sampling schedules in diffusion models.
\newblock In R.~Salakhutdinov, Z.~Kolter, K.~Heller, A.~Weller, N.~Oliver,
  J.~Scarlett, and F.~Berkenkamp, editors, {\em Proceedings of the 41st
  International Conference on Machine Learning}, volume 235 of {\em Proceedings
  of Machine Learning Research}, pages 42947--42975. PMLR, 21--27 Jul 2024.
\newblock \urlprefix\url{https://proceedings.mlr.press/v235/sabour24a.html}.

\bibitem[SME22]{song2022denoisingdiffusionimplicitmodels}
J.~Song, C.~Meng, and S.~Ermon.
\newblock Denoising diffusion implicit models, 2022,
  \burlalt{2010.02502}{http://arxiv.org/abs/2010.02502}.
\newblock \urlprefix\url{https://arxiv.org/abs/2010.02502}.

\bibitem[SSDK{\etalchar{+}}21]{SongSohlDicksteinEA21}
Y.~Song, J.~Sohl-Dickstein, D.~P. Kingma, A.~Kumar, S.~Ermon, and B.~Poole.
\newblock Score-based generative modeling through stochastic differential
  equations.
\newblock In {\em International Conference on Learning Representations}, 2021.
\newblock \urlprefix\url{https://openreview.net/forum?id=PxTIG12RRHS}.

\bibitem[Tho72]{Thom72}
R.~Thom.
\newblock Sur le cut-locus d'une vari\'et\'e{} plong\'ee.
\newblock {\em J. Differential Geometry}, 6:577--586, 1972.
\newblock \urlprefix\url{http://projecteuclid.org/euclid.jdg/1214430644}.

\bibitem[Tog73]{Tognoli73}
A.~Tognoli.
\newblock Su una congettura di {Nash}.
\newblock {\em Annali della Scuola Normale Superiore di Pisa, Classe di
  Scienze}, 27(1):167--185, 1973.
\newblock \urlprefix\url{https://www.numdam.org/item/ASNSP_1973_3_27_1_167_0/}.

\bibitem[Tre00]{Trefethen00}
L.~N. Trefethen.
\newblock {\em Spectral methods in {MATLAB}}, volume~10 of {\em Software,
  Environments, and Tools}.
\newblock Society for Industrial and Applied Mathematics (SIAM), Philadelphia,
  PA, 2000.
\newblock \doi{10.1137/1.9780898719598}.

\bibitem[Ver18]{Vershynin18}
R.~Vershynin.
\newblock {\em High-dimensional probability}, volume~47 of {\em Cambridge
  Series in Statistical and Probabilistic Mathematics}.
\newblock Cambridge University Press, Cambridge, 2018.
\newblock \doi{10.1017/9781108231596}.
\newblock An introduction with applications in data science, With a foreword by
  Sara van de Geer.

\bibitem[Vin11]{Vincent11}
P.~Vincent.
\newblock A connection between score matching and denoising autoencoders.
\newblock {\em Neural Computation}, 23(7):1661--1674, 2011.
\newblock \doi{10.1162/NECO_a_00142}.

\bibitem[Wal77]{Wall77}
C.~T.~C. Wall.
\newblock Geometric properties of generic differentiable manifolds.
\newblock In J.~Palis and M.~do~Carmo, editors, {\em Geometry and Topology},
  pages 707--774, Berlin, Heidelberg, 1977. Springer Berlin Heidelberg.
\newblock \doi{10.1007/BFb0085382}.

\bibitem[XTJ23]{xu2023stabletargetfieldreduced}
Y.~Xu, S.~Tong, and T.~Jaakkola.
\newblock Stable target field for reduced variance score estimation in
  diffusion models, 2023,
  \burlalt{2302.00670}{http://arxiv.org/abs/2302.00670}.
\newblock \urlprefix\url{https://arxiv.org/abs/2302.00670}.

\bibitem[Yom81]{Yomdin81}
Y.~Yomdin.
\newblock On the local structure of a generic central set.
\newblock {\em Compositio Math.}, 43(2):225--238, 1981.
\newblock \urlprefix\url{http://www.numdam.org/item?id=CM_1981__43_2_225_0}.

\bibitem[Yom15]{Yomdin15}
Y.~Yomdin.
\newblock Smooth parametrizations in dynamics, analysis, diophantine and
  computational geometry.
\newblock {\em Jpn. J. Ind. Appl. Math.}, 32(2):411--435, 2015.
\newblock \doi{10.1007/s13160-015-0176-6}.

\bibitem[YZS{\etalchar{+}}25]{YangZhangEA25}
L.~Yang, Z.~Zhang, Y.~Song, S.~Hong, R.~Xu, Y.~Zhao, W.~Zhang, B.~Cui, and
  M.-H. Yang.
\newblock Diffusion models: A comprehensive survey of methods and applications,
  2025, \burlalt{2209.00796}{http://arxiv.org/abs/2209.00796}.
\newblock \urlprefix\url{https://arxiv.org/abs/2209.00796}.

\bibitem[ZLWZ25]{zhang2025collapseerrorsinduceddeterministic}
Y.~Zhang, Z.~Liao, J.~Wu, and D.~Zou.
\newblock On the collapse errors induced by the deterministic sampler for
  diffusion models, 2025,
  \burlalt{2508.16154}{http://arxiv.org/abs/2508.16154}.
\newblock \urlprefix\url{https://arxiv.org/abs/2508.16154}.

\end{thebibliography}
\end{document}